\documentclass[a4paper,11pt]{article}
\usepackage{amsfonts,amssymb,amsmath}

\usepackage[english]{babel}
 \makeatletter
 \@addtoreset{equation}{subsection}
 \makeatother

 \makeatletter
 \@addtoreset{enunciato}{section}
 \makeatother

 \newcounter{enunciato}[subsection]

 \newtheorem{ittheorem}{Theorem}
 \newtheorem{itlemma}{Lemma}
 \newtheorem{itproposition}{Proposition}
 \newtheorem{itdefinition}{Definition}
 \newtheorem{itremark}{Remark}
 \newtheorem{itclaim}{Claim}
 \newtheorem{itfact}{Fact}
 \newtheorem{itconjecture}{Conjecture}

 \newenvironment{theorem}{\addtocounter{enunciato}{1}
 \begin{ittheorem}}{\end{ittheorem}}

 \newenvironment{lemma}{\addtocounter{enunciato}{1}
 \begin{itlemma}}{\end{itlemma}}

 \newenvironment{proposition}{\addtocounter{enunciato}{1}
 \begin{itproposition}}{\end{itproposition}}

 \newenvironment{definition}{\addtocounter{enunciato}{1}
 \begin{itdefinition}}{\end{itdefinition}}

 \newenvironment{remark}{\addtocounter{enunciato}{1}
 \begin{itremark}}{\end{itremark}}

 \newenvironment{claim}{\addtocounter{enunciato}{1}
 \begin{itclaim}}{\end{itclaim}}

 \newenvironment{fact}{\addtocounter{enunciato}{1}
 \begin{itfact}}{\end{itfact}}

 \newenvironment{conjecture}{\addtocounter{enunciato}{1}
 \begin{itconjecture}}{\end{itconjecture}}

 \newcommand{\be}[1]{\begin{equation}\label{#1}}
 \newcommand{\ee}{\end{equation}}

 \newcommand{\bl}[1]{\begin{lemma}\label{#1}}
 \newcommand{\el}{\end{lemma}}

 \newcommand{\br}[1]{\begin{remark}\label{#1}}
 \newcommand{\er}{\end{remark}}

 \newcommand{\bt}[1]{\begin{theorem}\label{#1}}
 \newcommand{\et}{\end{theorem}}

 \newcommand{\bd}[1]{\begin{definition}\label{#1}}
 \newcommand{\ed}{\end{definition}}

 \newcommand{\bcl}[1]{\begin{claim}\label{#1}}
 \newcommand{\ecl}{\end{claim}}

 \newcommand{\bfact}[1]{\begin{fact}\label{#1}}
 \newcommand{\efact}{\end{fact}}

 \newcommand{\bp}[1]{\begin{proposition}\label{#1}}
 \newcommand{\ep}{\end{proposition}}

 \newcommand{\bc}[1]{\begin{corollary}\label{#1}}
 \newcommand{\ec}{\end{corollary}}

 \newcommand{\bcj}[1]{\begin{conjecture}\label{#1}}
 \newcommand{\ecj}{\end{conjecture}}

 \newcommand{\bpr}{\begin{proof}}
 \newcommand{\epr}{\end{proof}}

 \newcommand{\bprl}[1]{\begin{proofof}{\it\ref{#1}}.\,\,}
 \newcommand{\eprl}{\end{proofof}}

 \newcommand{\bi}{\begin{itemize}}
 \newcommand{\ei}{\end{itemize}}

 \newcommand{\ben}{\begin{enumerate}}
 \newcommand{\een}{\end{enumerate}}

 \newenvironment{proof}{\noindent {\em Proof}.\,\,}{\hspace*{\fill}$\halmos$\medskip}
 \newenvironment{proofof}{\noindent {\em Proof of Lemma\,\,}}{\hspace*{\fill}$\halmos$\medskip}
 \newcommand{\halmos}{\rule{1ex}{1.4ex}}

 \newcommand{\one}{{\mathchoice {1\mskip-4mu\mathrm l}
         {1\mskip-4mu\mathrm l}
         {1\mskip-4.5mu\mathrm l}
         {1\mskip-5mu\mathrm l}}}

\def \Z {{\mathbb Z}}
\def \R {{\mathbb R}}

\def \N {{\mathbb N}}

\def \ra {\rightarrow}
\def \ba {\begin{array}}
\def \ea {\end{array}}
\def \da {\downarrow}

\def \P {{\mathbb P}}
\def \E {{\mathbb E}}

\def \Id {{\it Id}}
\def \ES {{\rm E}}
\def \PS {{\rm P}}
\def \Sp {{\rm Sp}}

\def \IRW {{\hbox{\tiny\rm IRW}}}
\def \RW {{\hbox{\tiny\rm RW}}}
\def \new {{\rm new}}
\def \diag {{\rm diag}}
\def \off {{\rm off}}
\def \c {{\rm c}}
\def \onek {{1[\kappa]}}

\def \cP {{\mathcal P}}
\def \cM {{\mathcal M}}

\def \cE {{\mathcal E}}
\def \cA {{\mathcal A}}
\def \cF {{\mathcal F}}
\def \cG {{\mathcal G}}

\def\one{\rlap{\mbox{\small\rm 1}}\kern.15em 1}

\begin{document}
\title{Intermittency on catalysts:\\
symmetric exclusion}

\author{\renewcommand{\thefootnote}{\arabic{footnote}}
J.\ G\"artner
\footnotemark[1]
\\
\renewcommand{\thefootnote}{\arabic{footnote}}
F.\ den Hollander
\footnotemark[2]\, \footnotemark[3]
\\
\renewcommand{\thefootnote}{\arabic{footnote}}
G.\ Maillard
\footnotemark[4]
}

\footnotetext[1]{
Institut f\"ur Mathematik, Technische Universit\"at Berlin,
Strasse des 17.\ Juni 136, D-10623 Berlin, Germany,
{\sl jg@math.tu-berlin.de}
}
\footnotetext[2]{
Mathematical Institute, Leiden University, P.O.\ Box 9512,
2300 RA Leiden, The Netherlands,
{\sl denholla@math.leidenuniv.nl}
}
\footnotetext[3]{
EURANDOM, P.O.\ Box 513, 5600 MB Eindhoven, The Netherlands
}
\footnotetext[4]{
Institut de Math\'ematiques, \'Ecole Polytechnique F\'ed\'erale
de Lausanne, CH-1015 Lausanne, Switzerland,
{\sl gregory.maillard@epfl.ch}
}

\maketitle

\begin{abstract}

We continue our study of intermittency for the parabolic Anderson equation
$\partial u/\partial t = \kappa\Delta u + \xi u$, where $u\colon\,\Z^d\times
[0,\infty)\to\R$, $\kappa$ is the diffusion constant, $\Delta$ is the discrete
Laplacian, and $\xi\colon\,\Z^d\times [0,\infty)\to\R$ is a space-time random
medium. The solution of the equation describes the evolution of a ``reactant''
$u$ under the influence of a ``catalyst'' $\xi$.

In this paper we focus on the case where $\xi$\ is exclusion with a symmetric
random walk transition kernel, starting from equilibrium with density
$\rho\in (0,1)$. We consider the annealed Lyapunov exponents, i.e., the exponential
growth rates of the successive moments of $u$. We show that these exponents
are trivial when the random walk is recurrent, but display an interesting dependence
on the diffusion constant $\kappa$ when the random walk is transient, with
qualitatively different behavior in different dimensions. Special attention
is given to the asymptotics of the exponents for $\kappa\to\infty$, which is
controlled by moderate deviations of $\xi$ requiring a delicate expansion argument.

In G\"artner and den Hollander \cite{garhol04} the case where $\xi$ is a Poisson
field of independent (simple) random walks was studied. The two cases show
interesting differences and similarities. Throughout the paper, a comparison
of the two cases plays a crucial role.

{\it Acknowledgment.} GM was supported by a postdoctoral fellowship from the
Netherlands Organization for Scientific Research (grant 613.000.307) during
his stay at EURANDOM. The research in this paper was partially supported by
the ESF Scientific Programme ``Random Dynamics in Spatially Extended Systems''.
FdH and GM are grateful to the Pacific Institute for the Mathematical Sciences
and the Mathematics Department of the University of British Columbia, Vancouver, 
Canada, for hospitality: FdH January to August 2006, GM mid-January to mid-February 
2006 when the work on this paper was completed.

\end{abstract}

\newpage


\section{Introduction and main results}
\label{S1}

\subsection{Model}
\label{S1.1}

The parabolic Anderson equation is the partial differential equation
\be{pA}
\frac{\partial}{\partial t}u(x,t) = \kappa\Delta u(x,t) + \xi(x,t)u(x,t),
\qquad x\in\Z^d,\,t\geq 0.
\ee
Here, the $u$-field is $\R$-valued, $\kappa\in [0,\infty)$ is the diffusion
constant, $\Delta$ is the discrete Laplacian, acting on $u$ as
\be{dL}
\Delta u(x,t) = \sum_{{y\in\Z^d} \atop {\|y-x\|=1}} [u(y,t)-u(x,t)]
\ee
($\|\cdot\|$ is the Euclidian norm), while
\be{rf}
\xi = \{\xi(x,t) \colon\,x\in\Z^d,\,t \geq 0\}
\ee
is an $\R$-valued random field that evolves with time and that drives
the equation. As initial condition for (\ref{pA}) we take
\be{ic}
u(\cdot,0) \equiv 1.
\ee

In the present paper we focus on the case where $\xi$ is \emph{Symmetric Exclusion}
(SE), i.e., $\xi$ takes values in $\{0,1\}^{\Z^d}\times [0,\infty)$,
where $\xi(x,t)=1$ means that there is a particle at $x$ at time $t$ and $\xi(x,t)=0$
means that there is none, and particles move around according to a symmetric random
walk transition kernel. We choose $\xi(\cdot,0)$ according to the Bernoulli product
measure with density $\rho\in (0,1)$, i.e., initially each site has a particle with
probability $\rho$ and no particle with probability $1-\rho$, independently for
different sites. For this choice, the $\xi$-field is stationary in time.

One interpretation of (\ref{pA}) and (\ref{ic}) comes from population dynamics.
Consider a spatially homogeneous system of two types of particles, $A$ (catalyst)
and $B$ (reactant), subject to:
\begin{itemize}
\item[(i)]
$A$-particles behave autonomously, according to a prescribed stationary dynamics,
with density $\rho$;
\item[(ii)]
$B$-particles perform independent random walks with diffusion constant $\kappa$
and split into two at a rate that is equal to the number of $A$-particles present
at the same location;
\item[(iii)]
the initial density of $B$-particles is $1$.
\end{itemize}
Then
\be{uint}
\begin{array}{ll}
u(x,t) = \{&\hbox{the average number of $B$-particles at site $x$ at time $t$}\\
           &\hbox{conditioned on the evolution of the $A$-particles}\,\,\,\}.
\end{array}
\ee
It is possible to add that $B$-particles die at rate $\delta\in (0,\infty)$.
This amounts to the trivial transformation $u(x,t) \to u(x,t)e^{-\delta t}$.

In Kesten and Sidoravicius \cite{kessid03} and in G\"artner and den Hollander
\cite{garhol04}, the case was considered where $\xi$ is given by a Poisson
field of independent simple random walks. The survival versus extinction pattern
(in \cite{kessid03} for $\delta>0$) and the annealed Lyapunov exponents (in
\cite{garhol04} for $\delta=0$) were studied, in particular, their dependence
on $d$, $\kappa$ and the parameters controlling $\xi$.

Equation (\ref{pA}) is a discrete heat equation with the $\xi$-field playing
the role of a source. What makes (\ref{pA}) particularly interesting is that the
two terms in the right-hand side \emph{compete with each other}: the diffusion
induced by $\Delta$ tends to make $u$ flat, while the branching induced by $\xi$
tends to make $u$ irregular. Henceforth we call $\xi$ the ``catalyst'' and $u$ 
the ``reactant''.

\subsection{SE, Lyapunov exponents and comparison with IRW}
\label{S1.2}

Throughout the paper, we abbreviate $\Omega=\{0,1\}^{\Z^d}$ (endowed with the
product topology), and we let $p\colon\Z^d\times\Z^d\to [0,1]$ be the transition
kernel of an irreducible random walk,
\be{pdef}
\begin{aligned}
& p(x,y)=p(0,y-x) \geq 0 \,\,\,\forall\,x,y\in\Z^d, \quad
\sum_{y\in\Z^d} p(x,y)=1 \,\,\,\forall\,x\in\Z^d,\\
& p(x,x)=0 \,\,\,\forall\,x\in\Z^d, \quad
p(\cdot,\cdot) \mbox{ generates } \Z^d,
\end{aligned}
\ee
that is assumed to be \emph{symmetric},
\be{pprop}
p(x,y)=p(y,x) \,\,\,\forall\,x,y\in\Z^d.
\ee
A special case is \emph{simple} random walk
\be{SRW}
p(x,y) =
\begin{cases}
\frac{1}{2d} & \text{if } \|x-y\|=1,\\
0            & \text{otherwise}.
\end{cases}
\ee

The exclusion process is the Markov process on $\Omega$ whose generator
$L$ acts on cylindrical functions $f$ as (see Liggett \cite{lig85}, Chapter VIII)
\be{expro1}
\begin{aligned}
(Lf)(\eta) &= \sum_{x,y\in\Z^d} p(x,y)\, \eta(x)[1-\eta(y)]
\left[f\left(\eta^{x,y}\right)-f(\eta)\right]
= \sum_{\{x,y\}\subset\Z^d} p(x,y)\,
\left[f\left(\eta^{x,y}\right)-f(\eta)\right],
\end{aligned}
\ee
where the latter sum runs over unoriented bonds $\{x,y\}$ between any pair of sites
$x,y\in\Z^d$, and
\be{expro2}
\eta^{x,y}(z) =
\begin{cases}
\eta(z)  &\text{if } z \neq x,y,\\
\eta(y)  & \text{if } z=x,\\
\eta(x)  & \text{if } z=y.
\end{cases}
\ee
The first line of (\ref{expro1}) says that a particle at site $x$ jumps to a vacancy at site
$y$ at rate $p(x,y)$, the second line says that the states of $x$ and $y$ are interchanged
along the bond $\{x,y\}$ at rate $p(x,y)$. For $\rho \in [0,1]$, let $\nu_\rho$ be the
Bernoulli product measure on $\Omega$ with density $\rho$. This is an invariant measure for SE.
Under (\ref{pdef}--\ref{pprop}), $(\nu_\rho)_{\rho\in [0,1]}$ are the only extremal equilibria
(see Liggett \cite{lig85}, Chapter VIII, Theorem 1.44). We denote by $\P_\eta$ the law of
$\xi$ starting from $\eta\in\Omega$ and
write $\P_{\nu_\rho} = \int_\Omega \nu_{\rho}(d\eta)\,\P_\eta$.

In the \emph{graphical representation} of SE, space is drawn sidewards, time is
drawn upwards, and for each pair of sites $x,y\in\Z^d$ links are drawn between
$x$ and $y$ at Poisson rate $p(x,y)$. The configuration at time $t$ is obtained
from the one at time $0$ by transporting the local states along paths that move
upwards with time and sidewards along links (see Fig.\ 1).

We will frequently use the following property, which is immediate from the
graphical representation:
\be{graph}
\E_{\,\eta}\left(\xi(y,t)\right) = \sum_{x\in\Z^d} \eta(x)\, p_t(x,y),
\qquad \eta\in\Omega,\,y\in\Z^d,\,t\geq 0.
\ee
Similar expressions hold for higher order correlations. Here, $p_t(x,y)$ is the
probability that the random walk with transition kernel $p(\cdot,\cdot)$ and step
rate $1$ moves from $x$ to $y$ in time $t$. The graphical representation shows that
the evolution is invariant under time reversal and, in particular, the equilibria
$(\nu_\rho)_{\rho\in [0,1]}$ are \emph{reversible}. This fact will turn out to be very 
important later on.


\vskip 0.8truecm

\setlength{\unitlength}{0.27cm}

\begin{picture}(20,10)(-9,1)

  \put(0,0){\line(22,0){22}}
  \put(0,11){\line(22,0){22}}

  \put(2,0){\line(0,12){12}}
  \put(5,0){\line(0,12){12}}
  \put(8,0){\line(0,12){12}}
  \put(11,0){\line(0,12){12}}
  \put(14,0){\line(0,12){12}}
  \put(17,0){\line(0,12){12}}
  \put(20,0){\line(0,12){12}}

  \qbezier[15](2.1,2)(3.5,2)(4.9,2)
  \qbezier[15](5.1,4)(6.5,4)(7.9,4)
  \qbezier[15](8.1,7)(9.5,7)(10.9,7)
  \qbezier[15](2.1,8)(3.5,8)(4.9,8)
  \qbezier[15](11.1,2.5)(12.5,2.5)(13.9,2.5)
  \qbezier[15](17.1,4)(18.5,4)(19.9,4)
  \qbezier[15](11.1,4.5)(12.5,4.5)(13.9,4.5)
  \qbezier[15](14.1,7.5)(15.5,7.5)(16.9,7.5)

  \put(10.7,-1.2){$x$}
  \put(8.4,11.4){$y$}
  \put(-1,-.3){$0$}
  \put(-1,10.7){$t$}

  \put(12,2.7){$\rightarrow$}
  \put(12,4.7){$\leftarrow$}
  \put(9,7.3){$\leftarrow$}
  \put(11.3,.8){$\uparrow$}
  \put(14.3,3){$\uparrow$}
  \put(11.3,5.5){$\uparrow$}
  \put(8.3,8.8){$\uparrow$}

  \put(11,0){\circle*{.35}}
  \put(8,11){\circle*{.35}}
  \put(23,0){$\Z^d$}

  \put(-1.5,-4){\small
               Fig.\ 1: Graphical representation. The dashed lines are links.
               \normalsize}
               \put(-1.5,-5.5){\small
               The arrows represent a path from $(x,0)$ to $(y,t)$.
               \normalsize}

\end{picture}

\vskip 2.2truecm


By the Feynman-Kac formula, the solution of (\ref{pA}) and (\ref{ic}) reads
\be{fey-kac1}
u(x,t) = \ES_{\,x}\left(\exp\left[\int_0^t ds\,\,
\xi\left(X^\kappa(s),t-s\right)\right]\right),
\ee
where $X^\kappa$ is simple random walk on $\Z^d$ with step rate $2d\kappa$ and $\ES_{\,x}$
denotes expectation with respect to $X^\kappa$ given $X^\kappa(0)=x$. We will often write
$\xi_t(x)$ and $X_t^\kappa$ instead of $\xi(x,t)$ and $X^\kappa(t)$, respectively.

For $p\in\N$ and $t>0$, define
\be{lyapdef}
\Lambda_p(t) = \frac{1}{pt} \log \E_{\,\nu_\rho}\left(u(0,t)^p\right).
\ee
Then
\be{fey-kac2}
\Lambda_p(t) = \frac{1}{pt} \log \E_{\,\nu_\rho}\bigg(\ES_{\,0,\dots,0}
\bigg(\exp\bigg[\int_0^t ds \sum_{q=1}^{p}
\xi\big(X_q^\kappa(s),s\big)\bigg]\bigg)\bigg),
\ee
where $X_q^\kappa$, $q=1,\dots,p$, are $p$ independent copies of $X^\kappa$,
$\ES_{\,0,\dots,0}$ denotes expectation w.r.t.\ $X_q^\kappa, q=1,\dots,p$, given
$X_1^\kappa(0)=\cdots=X_p^\kappa(0)=0$, and the time argument $t-s$ in
(\ref{fey-kac1}) is replaced by $s$ in (\ref{fey-kac2}) via the reversibility
of $\xi$ starting from $\nu_\rho$. If the last quantity admits a limit as
$t\to\infty$, then we define
\be{lyap2}
\lambda_p = \lim_{t\to\infty} \Lambda_p(t)
\ee
to be the $p$-th \emph{annealed Lyapunov exponent}.

{}From H\"older's inequality applied to (\ref{lyapdef}) it follows that
$\Lambda_p(t)\geq\Lambda_{p-1}(t)$ for all $t>0$ and $p\in\N\setminus\{1\}$.
Hence $\lambda_p\geq\lambda_{p-1}$ for all $p\in\N\setminus\{1\}$. We say
that the system is \emph{$p$-intermittent} if $\lambda_p>\lambda_{p-1}$.
In the latter case the system is $q$-intermittent for all $q>p$ as well
(cf. G\"artner and Molchanov \cite{garmol90}, Section 1.1). We say that the
system is \emph{intermittent} if it is $p$-intermittent for all
$p\in\N\setminus\{1\}$. Intermittent means that the $u$-field develops sparse
high peaks dominating the moments in such a way that each moment is dominated
by its own collection of peaks (see G\"artner and K\"onig \cite{garkon04},
Section 1.3, and den Hollander \cite{garhol04}, Section 1.2).

Let $(\tilde\xi_{t})_{t\geq 0}$ be the process of \emph{Independent Random Walks} (IRW)
with step rate 1, transition kernel $p(\cdot,\cdot)$ and state space $\Omega$. Let 
$\E_{\eta}^{\IRW}$ denote expectation w.r.t.\ $(\tilde\xi_{t})_{t\geq 0}$ starting 
from $\tilde\xi_{0}=\eta$, and write $\E_{\,\nu_\rho}^{\IRW}=\int_\Omega\nu_\rho(d\eta)\,
\E_{\eta}^{\IRW}$. Throughout the paper we will make use of the following inequality 
comparing SE and IRW. The proof of this inequality is given in Appendix \ref{A} and 
uses a lemma due to Landim \cite{lan92}.

\bp{IRW-comp}
For any $K\colon\,\Z^d\times[0,\infty)\to \R$ such that either $K\geq 0$ or $K\leq 0$,
any $t\geq 0$ such that $\sum_{z\in\Z^d}\int_0^t ds\,|K(z,s)|<\infty$ and any $\eta\in\Omega$,
\be{IRW-comp-ineq}
\E_{\,\eta}\Bigg(\exp\Bigg[\sum_{z\in \Z^d}\int_0^t ds\,\,
K(z,s)\, \xi_{s}(z)\Bigg]\Bigg)
\leq \E_{\,\eta}^{\IRW}\Bigg(\exp\Bigg[\sum_{z\in \Z^d}\int_0^t ds\,\,
K(z,s)\, \tilde\xi_{s}(z)\Bigg]\Bigg).
\ee
\ep

\noindent
This powerful inequality will allow us to obtain bounds that are more easily computable.

\subsection{Main theorems}
\label{S1.3}

Our first result states that the Lyapunov exponents exist and behave nicely as a
function of $\kappa$. We write $\lambda_p(\kappa)$ to exhibit the dependence on
$\kappa$, suppressing $d$ and $\rho$.

\bt{Lyaexist}
Let $d\geq 1$, $\rho\in (0,1)$ and $p\in\N$.\\
(i) For all $\kappa\in [0,\infty)$, the limit in {\rm (\ref{lyap2})} exists and is finite.\\
(ii) On $[0,\infty)$, $\kappa\to\lambda_p(\kappa)$ is continuous, non-increasing and convex.
\et

Our second result states that the Lyapunov exponents are trivial for recurrent
random walk but are non-trivial for transient random walk (see Fig.\ 2).

\bt{Lyalow}
Let $d\geq 1$, $\rho\in (0,1)$ and $p\in\N$.\\
(i) If $p(\cdot,\cdot)$ is recurrent, then $\lambda_p(\kappa)=1$ for all $\kappa
\in [0,\infty)$.\\
(ii) If $p(\cdot,\cdot)$ is transient, then $\rho<\lambda_p(\kappa)<1$ for all
$\kappa\in [0,\infty)$. Moreover, $\kappa\mapsto\lambda_p(\kappa)$ is strictly
decreasing with $\lim_{\kappa\to\infty}\lambda_p(\kappa)=\rho$.
\et


\vspace{0.5cm}

\setlength{\unitlength}{0.3cm}

\begin{picture}(20,8)(-7,0)

  \put(0,0){\line(8,0){8}}
  \put(0,0){\line(0,7){7}}
  {\thicklines
   \qbezier(0,5)(3,5)(6,5)
  }
  \put(-.8,-1.3){$0$}
  \put(-1,4.8){$1$}
  \put(9,-.3){$\kappa$}
  \put(-1,8){$\lambda_p(\kappa)$}

  \put(15,0){\line(8,0){8}}
  \put(15,0){\line(0,7){7}}
  {\thicklines
   \qbezier(15,4)(18,2.5)(21,2.2)
  }
  \qbezier[60](15,2)(18,2)(22,2)
  \qbezier[60](15,5)(18,5)(22,5)
  \put(14.2,-1.3){$0$}
  \put(14,4.8){$1$}
  \put(14,1.8){$\rho$}
  \put(24,-.3){$\kappa$}
  \put(14,8){$\lambda_p(\kappa)$}
  \put(0,5){\circle*{.45}}
  \put(15,4){\circle*{.45}}

  \put(-1,-3.5){\small
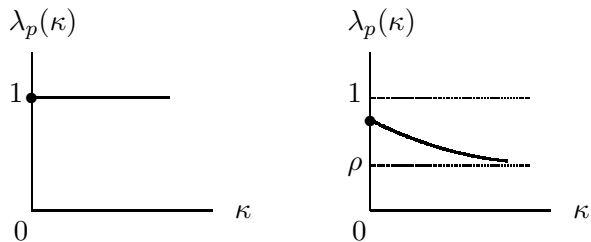
             Fig.\ 2: Qualitative picture of $\kappa\mapsto\lambda_p(\kappa)$
             for recurrent, respectively,
             \normalsize}
  \put(-1,-4.8){\small
              transient random walk.
              \normalsize}

\end{picture}

\vskip 2truecm


Our third result shows that for transient random walk the system is intermittent
at $\kappa=0$.

\bt{Lyaint}
Let $d\geq 1$ and $\rho\in (0,1)$. If $p(\cdot,\cdot)$ is transient, then $p\mapsto\lambda_p(0)$
is strictly increasing.
\et

Our fourth and final result identifies the behavior of the Lyapunov exponents
for large $\kappa$ when $d\geq 4$ and $p(\cdot,\cdot)$ is simple random walk
(see Fig.\ 3).

\bt{Lyahighlim}
Assume {\rm (\ref{SRW})}. Let $d\geq 4$, $\rho\in (0,1)$ and $p\in\N$ . Then
\be{limlamb}
\lim_{\kappa\to\infty} 2d\kappa[\lambda_p(\kappa)-\rho] =
\rho(1-\rho)G_d
\ee
with $G_d$ the Green function at the origin of simple random walk on $\Z^d$.
\et


\vskip 1truecm

\setlength{\unitlength}{0.25cm}

\begin{picture}(20,10)(-10,1)

  \put(0,-2){\line(18,0){18}}
  \put(0,-2){\line(0,14){14}}
  {\thicklines
   \qbezier(0,8)(2,6.6)(4,5.3)
   \qbezier(0,6)(2,4.8)(4,3.8)
   \qbezier(0,4)(2,3)(4,2.25)
  }
  {\thicklines
  \qbezier[15](11,0.8)(13,0.5)(16,0.4)
  }
  \qbezier[60](0,0)(9,0)(17,0)
  \qbezier[60](0,10.5)(9,10.5)(17,10.5)
  \put(-.8,-3.3){$0$}
  \put(-1.2,-.3){$\rho$}
  \put(-1.2,10.2){$1$}
  \put(0,8){\circle*{.45}}
  \put(0,6){\circle*{.45}}
  \put(0,4){\circle*{.45}}
  \put(-4,8){$p=3$}
  \put(-4,6){$p=2$}
  \put(-4,4){$p=1$}
  \put(7,2.5){{\bf ?}}
  \put(19,-2.3){$\kappa$}
  \put(-1.5,12.7){$\lambda_p(\kappa)$}
  \put(-4,-6){\small
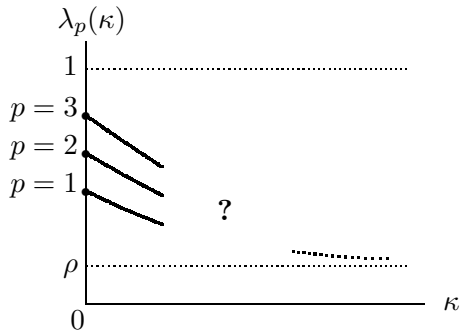
               Fig.\ 3: Qualitative picture of $\kappa\mapsto\lambda_p(\kappa)$
               for $p=1,2,3$ for simple
               \normalsize}
  \put(-4,-7.5){\small
               random walk in $d\geq 4$. The dotted line moving down represents
               \normalsize}
  \put(-4,-9){\small
                the asymptotics given by the r.h.s.\ of (\ref{limlamb}).
               \normalsize}
\end{picture}

\vskip 4truecm


\subsection{Discussion}
\label{S1.4}

Theorem \ref{Lyaexist} gives general properties that need no further comment.
We will see that they in fact hold for \emph{any stationary, reversible and
bounded $\xi$}.

The intuition behind Theorem \ref{Lyalow} is the following. If the catalyst is
driven by a recurrent random walk, then it suffers from ``traffic jams'', i.e.,
with not too small a probability there is a large region around the origin that
the catalyst fully occupies for a long time. Since with not too small a
probability the simple random walk (driving the reactant) can stay inside
this large region for the same amount of time, the average growth rate of the
reactant at the origin is maximal. This phenomenon may be expressed by saying
that \emph{for recurrent random walk clumping of the catalyst dominates the
growth of the moments}. For transient random walk, on the other hand, clumping
of the catalyst is present (the growth rate of the reactant is $>\rho$), but it
is not dominant (the growth rate of the reactant is $<1$). As the diffusion
constant $\kappa$ of the reactant increases, the effect of the clumping of the
catalyst gradually diminishes and the growth rate of the reactant gradually
decreases to the density of the catalyst.

Theorem \ref{Lyaint} shows that if the reactant stands still and the catalyst
is driven by a transient random walk, then the system is intermittent. Apparently,
the successive moments of the reactant, which are equal to the exponential moments
of the occupation time of the origin by the catalyst (take (\ref{fey-kac1}) with
$\kappa=0$), are sensitive to successive degrees of clumping. By continuity,
\emph{intermittency persists for small} $\kappa$.

Theorem \ref{Lyahighlim} shows that, when the catalyst is driven by simple random walk,
\emph{all} Lyapunov exponents decay to $\rho$ as $\kappa\to\infty$ \emph{in the same
manner} when $d \geq 4$. The case $d=3$ remains open. We conjecture:

\bcj{cjd=3}
Assume {\rm (\ref{SRW})}. Let $d=3$, $\rho\in (0,1)$ and $p\in\N$ . Then
\be{limlamb*}
\lim_{\kappa\to\infty} 2d\kappa[\lambda_p(\kappa)-\rho] =
\rho(1-\rho)G_d + [2d\rho(1-\rho)]^2\cP
\ee
with
\be{Rpdef}
\cP = \sup_{{f \in H^1(\R^3)} \atop {\|f\|_2=1}}
\left[\,\left\|\left(-\Delta_{\R^3}\right)^{-1/2}\,f^2\right\|_2^2
- \left\|\nabla_{\R^3} f\right\|_2^2\,\right]
\in (0,\infty),
\ee
where $\nabla_{\R^3}$ and $\Delta_{\R^3}$ are the continuous gradient and
Laplacian, $\|\cdot\|_2$ is the $L^2(\R^3)$-norm, $H^1(\R^3)=\{f\colon\,
\R^3\to\R\colon\,f,\nabla_{\R^3}f\in L^2(\R^3)\}$, and
\be{Gint}
\left\|\left(-\Delta_{\R^3}\right)^{-1/2}\,f^2\right\|_2^2
= \int_{\R^3} dx\,f^2(x) \int_{\R^3} dy\,f^2(y)\,\,\frac{1}{4\pi\|x-y\|}.
\ee
\ecj
In section \ref{S1.5} we will explain how this conjecture arises in analogy with the
case of IRW studied in G\"artner and den Hollander \cite{garhol04}.
If Conjecture \ref{cjd=3} holds true, then in $d=3$ \emph{intermittency persists for large}
$\kappa$. It would still remain open whether the same is true for $d\geq 4$. To decide the
latter, we need a finer asymptotics for $d\geq 4$. A large diffusion constant of the reactant
prevents the solution $u$ to easily localize around the regions where the catalyst clumps, but
it is not clear whether this is able to destroy intermittency for $d\geq 4$.

We further conjecture:

\bcj{cj1}
In $d=3$, the system is intermittent for all $\kappa\in [0,\infty)$.
\ecj

\bcj{cj2}
In $d\geq 4$, there exists a strictly increasing sequence $0<\kappa_2<\kappa_3
<\ldots$ such that for $p=2,3,\ldots$ the system is $p$-intermittent if and only
if $0\leq\kappa<\kappa_p$.
\ecj

\noindent
In words, we conjecture that in $d=3$ the curves in Fig.\ 3 never merge, whereas for $d\geq 4$
the curves merge successively.

\subsection{Heuristics behind Theorem \ref{Lyahighlim} and Conjecture \ref{cjd=3}}
\label{S1.5}

The \emph{heuristics} behind Theorem \ref{Lyahighlim} and Conjecture \ref{cjd=3} is
the following. Consider the case $p=1$. Scaling time by $\kappa$ in (\ref{fey-kac2}),
we have $\lambda_1(\kappa)=\kappa \lambda_1^*(\kappa)$ with
\be{lambda*scal}
\lambda_1^*(\kappa) = \lim_{t\to\infty} \Lambda_1^*(\kappa;t)
\quad \mbox{and} \quad \Lambda_1^*(\kappa;t)
= \frac{1}{t} \log \E_{\,\nu_\rho,0}
\left(\exp\left[\frac{1}{\kappa} \int_0^t ds\,\,
\xi\Big(X(s),\frac{s}\kappa\Big)\right]\right),
\ee
where $X=X^1$ and we abbreviate
\be{Edefs}
\E_{\,\nu_\rho,0} = \E_{\,\nu_\rho}\ES_{\,0}.
\ee
For large $\kappa$, the $\xi$-field in (\ref{lambda*scal}) evolves slowly and
therefore does not manage to cooperate with the $X$-process in determining the
growth rate. Also, the prefactor $1/\kappa$ in the exponent is small. As a result,
the expectation over the $\xi$-field can be computed via a \emph{Gaussian
approximation} that becomes sharp in the limit as $\kappa\to\infty$, i.e.,
\be{Gappr1}
\begin{aligned}
&\Lambda_1^*(\kappa;t) - \frac{\rho}{\kappa}
= \frac{1}{t} \log \E_{\,\nu_\rho,0}
\left(\exp\left[\frac{1}{\kappa} \int_0^t ds\,\,
\Big[\xi\Big(X(s),\frac{s}\kappa\Big)-\rho\Big]\right]\right)\\
&\qquad \approx \frac{1}{t} \log
\ES_{\,0}\left(\exp\left[\frac{1}{2\kappa^2}
\int_0^t ds \int_0^t du\,\,
\E_{\,\nu_\rho}\bigg(\Big[\xi\Big(X(s),\frac{s}\kappa\Big)-\rho\Big]
\Big[\xi\Big(X(u),\frac{u}\kappa\Big)-\rho\Big]\bigg)\right]\right).
\end{aligned}
\ee
(In essence, what happens here is that the asymptotics for $\kappa\to\infty$ is driven
by moderate deviations of the $\xi$-field, which fall in the Gaussian regime.) The exponent
in the r.h.s.\ of (\ref{Gappr1}) equals
\be{Gappr2}
\frac{1}{\kappa^2} \int_0^t ds \int_s^t du\,\,
\E_{\,\nu_\rho}\bigg(\Big[\xi\Big(X(s),\frac{s}\kappa\Big)-\rho\Big]
\Big[\xi\Big(X(u),\frac{u}\kappa\Big)-\rho\Big]\bigg).
\ee
Now, for $x,y\in\Z^d$ and $b \geq a \geq 0$ we have
\be{Gappr3}
\begin{aligned}
\E_{\,\nu_\rho}\bigg(\big[\xi(x,a)-\rho\big]\big[\xi(y,b)-\rho\big]\bigg)
&= \E_{\,\nu_\rho}\bigg(\big[\xi(x,0)-\rho\big]
\big[\xi(y,b-a)-\rho\big]\bigg)\\
&= \int_\Omega \nu_\rho(d\eta)\,\big[\eta(x)-\rho\big]
\E_{\,\eta}\bigg(\big[\xi(y,b-a)-\rho\big]\bigg)\bigg)\\
&= \sum_{z\in\Z^d} p_{b-a}(z,y)\,\int_\Omega \nu_\rho(d\eta)\,
\big[\eta(x)-\rho\big]\big[\eta(z)-\rho\big]\\
&= \rho(1-\rho)\,p_{b-a}(x,y),
\end{aligned}
\ee
where the first equality uses the stationarity of $\xi$, the third equality
uses (\ref{graph}) from the graphical representation, and the fourth equality
uses that $\nu_\rho$ is Bernoulli. Substituting (\ref{Gappr3}) into (\ref{Gappr2}),
we get that the r.h.s.\ of (\ref{Gappr1}) equals
\be{Gappr5}
\frac{1}{t} \log
\ES_{\,0}\left(\exp\left[\frac{\rho(1-\rho)}{\kappa^2}\int_0^t ds \int_s^t du\,\,
p_{\frac{u-s}\kappa}(X(s),X(u))\right]\right).
\ee
This is precisely the integral that was investigated in G\"artner and den Hollander 
\cite{garhol04} (see Sections 5--8 and equations (1.5.4--1.5.11) of that paper).
Therefore the limit
\be{limfinal}
\lim_{\kappa\to\infty} \kappa[\lambda_1(\kappa)-\rho]
= \lim_{\kappa\to\infty}
\kappa^2\lim_{t\to\infty}\left[\Lambda_1^*(\kappa;t)-\frac{\rho}{\kappa}\right]
= \lim_{\kappa\to\infty}
\kappa^2\lim_{t\to\infty} \mbox{{\rm (\ref{Gappr5})}}
\ee
can be read off from \cite{garhol04} and yields (\ref{limlamb}) for $d \geq 4$ and 
(\ref{limlamb*}) for $d=3$. A similar heuristics applies for $p>1$.

The r.h.s.\ of (\ref{limlamb}), which is valid for $d \geq 4$, is obtained from the 
above computations by moving the expectation in (\ref{Gappr5}) into the exponent. 
Indeed,
\be{intavcal1}
\ES_{\,0}\left(p_{\frac{u-s}\kappa}(X(s),X(u))\right)
= \sum_{x,y\in\Z^d} p_{2ds}(0,x)p_{2d(u-s)}(x,y)p_{\frac{u-s}\kappa}(x,y)
= p_{2d(u-s)(1+\frac{1}{2d\kappa})}(0,0)
\ee
and hence
\be{intavcal2}
\int_0^t ds \int_s^t du\,\,\ES_{\,0}\left(p_{\frac{u-s}\kappa}(X(s),X(u))\right)
= \int_0^t ds \int_0^{t-s} dv\,p_{2dv(1+\frac{1}{2d\kappa})}(0,0)
\sim t\, \frac{1}{2d(1+\frac{1}{2d\kappa})}G_d.
\ee
Thus we see that the result in Theorem \ref{Lyahighlim} comes from a second order
asymptotics on $\xi$ and a first order asymptotics on $X$. Despite this simple fact,
it turns out to be hard to make the above heuristics rigorous. For $d=3$, on the
other hand, we expect the first order asymptotics on $X$ to fail, leading to the
more complicated behavior in (\ref{limlamb*}) .

\medskip\noindent
{\bf Remark 1:} In (\ref{pA}), the $\xi$-field may be multiplied by a coupling
constant $\gamma\in (0,\infty)$. This produces no change in Theorems \ref{Lyaexist},
\ref{Lyalow}(i) and \ref{Lyaint}. In Theorem \ref{Lyalow}(ii), $(\rho,1)$ becomes
$(\gamma\rho,\gamma)$, while in the r.h.s.\ of Theorem \ref{Lyahighlim} and Conjecture
\ref{cjd=3}, $\rho(1-\rho)$ gets multiplied by $\gamma^2$. Similarly, if the simple
random walk in Theorem \ref{Lyahighlim} is replaced by a random walk with transition
kernel $p(\cdot,\cdot)$ satisfying (\ref{pdef}--\ref{pprop}), then we expect that in
(\ref{limlamb}) and (\ref{limlamb*}) $G_d$ becomes the Green function at the origin of
this random walk and a factor $1/\sigma^4$ appears in front of the last term in the
r.h.s.\ of (\ref{limlamb*}) with $\sigma^2$ the variance of $p(\cdot,\cdot)$.

\medskip\noindent
{\bf Remark 2:} In G\"artner and den Hollander \cite{garhol04} the catalyst was $\gamma$
times a Poisson field with density $\rho$ of independent simple random walks stepping at
rate $2d\theta$, where $\gamma,\rho,\theta\in (0,\infty)$ are parameters. It was found
that the Lyapunov exponents are infinite in $d=1,2$ for all $p$ and in $d\geq 3$ for
$p\geq 2d\theta/\gamma G_d$, irrespective of $\kappa$ and $\rho$. In $d \geq 3$ for 
$p<2d\theta/\gamma G_d$, on the other hand, the Lyapunov exponents are finite for all 
$\kappa$, and exhibit a dichotomy similar to the one expressed by Theorem \ref{Lyahighlim} 
and Conjecture \ref{cjd=3}. Apparently, in this regime the two types of catalyst are 
qualitatively similar. Remarkably, the \emph{same} asymptotic behavior for large $\kappa$ 
was found (with $\rho\gamma^2$ replacing $\rho(1-\rho)$ in (\ref{limlamb})), and the 
\emph{same} variational formula as in (\ref{Rpdef}) was seen to play a central role in 
$d=3$. [Note: In \cite{garhol04} the symbols $\nu,\rho,G_d$ were used instead of 
$\rho,\theta,G_d/2d$.]

\subsection{Outline}
\label{S1.6}

In Section \ref{S2} we derive a variational formula for $\lambda_p$ from which
Theorem \ref{Lyaexist} follows immediately. The arguments that will be used to
derive this variational formula apply to an arbitrary bounded, stationary and
\emph{reversible} catalyst. Thus, the properties in Theorem \ref{Lyaexist} are
quite general. In Section \ref{S3} we do a range of estimates, either directly
on (\ref{fey-kac2}) or on the variational formula for $\lambda_p$ derived in
Section \ref{S2}, to prove Theorems \ref{Lyalow} and \ref{Lyaint}. Here, the
special properties of SE, in particular, its space-time correlation structure
expressed through the graphical representation (see Fig.$\;1$), are crucial.
These results hold for an arbitrary random walk subject to (\ref{pdef}--\ref{pprop}).
Finally, in Section \ref{S4} we prove Theorem \ref{Lyahighlim}, which is restricted
to simple random walk. The analysis consists of a long series of estimates, taking
up half of the paper and, in essence, showing that the problem reduces to understanding 
the asymptotic behavior of (\ref{Gappr5}). This reduction is important, because it 
explains why there is some degree of \emph{universality} in the behavior for 
$\kappa\to\infty$ under different types of catalysts: apparently, the Gaussian 
approximation and the two-point correlation function in space and time determine 
the asymptotics (recall the heuristic argument in Section \ref{S1.5}).


\section{Lyapunov exponents: general properties}
\label{S2}

In this section we prove Theorem \ref{Lyaexist}. In Section
\ref{S2.1} we formulate a large deviation principle for the occupation time of the
origin in SE due to Landim \cite{lan92}, which will be needed in Section \ref{S3.2}.
In Section \ref{S2.2} we extend the line of thought in \cite{lan92} and derive a
variational formula for $\lambda_p$ from which Theorem \ref{Lyaexist} will follow
immediately.

\subsection{Large deviations for the occupation time of the origin}
\label{S2.1}

Kipnis \cite{kip87}, building on techniques developed by Arratia \cite{arr85},
proved that the \emph{occupation time of the origin up to time $t$},
\be{occtime}
T_t = \int_0^t \xi(0,s)\, ds,
\ee
satisfies a strong law of large numbers and a central limit theorem. Landim
\cite{lan92} subsequently proved that $T_t$ satisfies a \emph{large deviation
principle}, i.e.,
\be{LDPSE}
\begin{aligned}
\limsup_{t\to\infty} \frac{1}{t} \log \P_{\nu_\rho}
\left(T_t/t \in F \right) &\leq - \inf_{\alpha \in F} \Psi_d(\alpha),
\qquad F\subseteq [0,1] \text{ closed},\\
\liminf_{t\to\infty} \frac{1}{t} \log \P_{\nu_\rho}
\left(T_t/t \in G \right) &\geq - \inf_{\alpha \in G} \Psi_d(\alpha),
\qquad G\subseteq [0,1] \text{ open},
\end{aligned}
\ee
with $\Psi_d\colon\,[0,1]\to [0,\infty)$ a rate function that is given by the
following formulas. Define the Dirichlet form associated with the generator $L$
of SE given in (\ref{expro1}),
\be{Diridef}
\cE(f) = \left(-L f,f\right)_{L^2(\nu_\rho)}
= \int_\Omega \nu_{\rho}(d \eta)\,\, \frac12 \sum_{\{x,y\}\subset\Z^d} p(x,y)\,
\left[f\left(\eta^{x,y}\right)- f(\eta)\right]^2,
\quad f\in L^2(\nu_\rho).
\ee
Next, define $I_d\colon\,\cM_1(\Omega)\to[0,\infty]$ by
\be{rdSEb}
I_d(\mu) = \lim\limits_{\varepsilon\da 0}
\inf\limits_{{\nu\in B(\mu,\varepsilon)} \atop {\nu\ll\nu_\rho}}
\cE\left(\sqrt{\frac{d\nu}{d\nu_\rho}}\,\right),
\ee
where $\cM_1(\Omega)$ is the set of probability measures on $\Omega$ (endowed
with the Prokhorov metric), and $B(\mu,\varepsilon)$ is the open ball of radius
$\varepsilon$ centered at $\mu$. Then
\be{rfSEa}
\Psi_d(\alpha) = \inf_{{\mu\in\cM_1(\Omega)} \atop {\int_\Omega \eta(0)\mu(d\eta)=\alpha}}
I_d(\mu).
\ee

The function $I_d$ is convex and lower semi-continuous. Since $\Psi_d$ is a minimum
of $I_d$ under a linear constraint, it inherits these properties (see e.g.\ den
Hollander \cite{hol00}, Theorem III.32).

\medskip\noindent
{\bf Remark:} The function $\tilde I_d(\mu)$ defined by $\cE(\sqrt{d\mu/d\nu_\rho})$
if $\mu\ll\nu_\rho$ and $\infty$ otherwise is not lower semi-continuous, and therefore
does not qualify as a large deviation rate function. This is why the regularization in
(\ref{rdSEb}) is necessary (see also Deuschel and Stroock \cite{deustr89}, Section 5.3).

\medskip
Landim \cite{lan92} showed that $\Psi_d$ has a unique zero at $\rho$ and satisfies
the quadratic lower bound
\be{rflb}
\Psi_d(\alpha) \geq \frac{1}{2G_d} \left(\sqrt{\alpha}-\sqrt{\rho}\right)^2
\ee
with $G_d$ the Green function at the origin of the random walk with transition
kernel $p(\cdot,\cdot)$. This bound was obtained with the help of Proposition
\ref{IRW-comp} with $K(z,s)=\delta_0(z)$, which implies that the occupation time
for SE is stochastically smaller than the occupation time for IRW with the same
density (see \cite{lan92}, Proposition 4.1). For the latter the rate function can
be computed and equals the lower bound in (\ref{rflb}) for $\alpha \geq \rho$. The
same lower bound holds for $\alpha \leq \rho$, which can be seen by interchanging
the role of the states $0$ and $1$.

Thus, as seen from (\ref{rflb}), for transient random walk the rate function $\Psi_d$ 
is non-trivial. For recurrent random walk $\Psi_d$ turns out to be zero, so a different 
scaling is needed in (\ref{LDPSE}) to get a non-trivial large deviation principle. This is
carried out in Landim \cite{lan92} for $d=1$ and in Chang, Landim and Lee
\cite{chalanlee04} for $d=2$ (when $p(\cdot,\cdot)$ has positive and finite
variance). We shall, however, not need this refinement.

\subsection{Variational formula for $\lambda_p(\kappa)$: proof of Theorem
\ref{Lyaexist}}
\label{S2.2}

Return to (\ref{fey-kac2}). In this section we show that, by considering $\xi$ and
$X^\kappa_1,\dots,X^\kappa_p$ as a joint random process and exploiting the reversibility 
of $\xi$, we can use the spectral theorem to express the Lyapunov exponents in terms of 
a variational formula. From the latter it will follow that $\kappa\mapsto\lambda_p(\kappa)$ 
is continuous, non-increasing and convex on $[0,\infty)$.

Define
\be{comproc}
Y(t) = \left(\xi(t),X^\kappa_1(t),\dots,X^\kappa_p(t)\right), \quad t\geq 0,
\ee
and
\be{pot}
V(\eta,x_1,\dots,x_p) = \sum_{i=1}^p \eta(x_i), \quad \eta\in\Omega,
\,x_1,\dots,x_p\in\Z^d.
\ee
Then we may write (\ref{fey-kac2}) as
\be{combLamb}
\Lambda_p(t) = \frac{1}{pt} \log \E_{\,\nu_\rho,0,\dots,0}
\left(\exp\left[\int_0^t V(Y(s))ds\right]\right).
\ee
The random process $Y=(Y(t))_{t\geq 0}$ takes values in $\Omega \times (\Z^d)^p$ and
has generator
\be{gen}
G^\kappa = L+\kappa\sum_{i=1}^p \Delta_i
\ee
in $L^2(\nu_\rho\otimes m^p)$ (endowed with the inner product $(\cdot,\cdot)$), with
$L$ given by (\ref{expro1}), $\Delta_i$ the discrete Laplacian acting on the $i$-th
spatial coordinate, and $m$ the counting measure on $\Z^d$. Let
\be{combgen}
G^\kappa_V = G^\kappa + V.
\ee
By (\ref{pprop}), this is a \emph{self-adjoint} operator. Our claim is that $\lambda_p$ equals
$\frac{1}{p}$ times the upper boundary of the spectrum of $G^\kappa_V$.

\bp{LyaRR}
$\lambda_p = \frac{1}{p}\mu_p$ with $\mu_p=\sup\Sp\left(G^\kappa_V\right)$.
\ep

\bpr
The proof is standard. Let $(\cP_t)_{t\geq 0}$ denote the semigroup generated by
$G^\kappa_V$.

\medskip\noindent
\underline{Upper bound}:
Let $Q_{t\log t} = [-t\log t,t \log t]^d \cap \Z^d$. By a standard large deviation
estimate for simple random walk, we have
\be{LDRW}
\begin{aligned}
&\E_{\,\nu_\rho,0,\dots,0}
\left(\exp\left[\int_0^t V(Y(s))ds\right]\right)\\
&\qquad = \E_{\,\nu_\rho,0,\dots,0}
\left(\exp\left[\int_0^t V(Y(s))ds\right]
\one\left\{X^\kappa_i(t)\in Q_{t\log t} \mbox{ for } i=1,\dots,p\right\}\right)
+ R_t
\end{aligned}
\ee
with $\lim_{t\to\infty}\frac{1}{t}\log R_t = -\infty$. Thus it suffices to focus
on the term with the indicator.

Estimate, with the help of the spectral theorem (Kato \cite{ka76}, Section VI.5),
\be{lb1}
\begin{aligned}
&\E_{\,\nu_\rho,0,\dots,0}
\left(\exp\left[\int_0^t V(Y(s))ds\right]
\one\left\{X^\kappa_i(t)\in Q_{t\log t} \mbox{ for } i=1,\dots,p\right\}\right)\\
&\quad\leq \sum_{x_1,\dots,x_p \in Q_{t\log t}} \E_{\,\nu_\rho,x_1,\dots,x_p}
\left(\exp\left[\int_0^t V(Y(s))ds\right]
\one\left\{X^\kappa_i(t)\in Q_{t\log t} \mbox{ for } i=1,\dots,p\right\}\right)\\
&\quad= \left(\one_{(Q_{t\log t})^p}, \cP_t \one_{(Q_{t\log t})^p}\right)
= \int_{(-\infty,\mu_p]} e^{\mu t}\,
d\|E_\mu \one_{(Q_{t\log t})^p}\|^2_{L^2(\nu_\rho\otimes m^p)}\\
&\quad\leq e^{\mu_p t}\, \|\one_{(Q_{t\log t})^p}\|^2_{L^2(\nu_\rho\otimes m^p)},
\end{aligned}
\ee
where $\one_{(Q_{t\log t})^p}$ is the indicator function of $(Q_{t\log t})^p \subset(\Z^d)^p$
and $(E_\mu)_{\mu \in \R}$ denotes the spectral family of orthogonal projection operators
associated with $G^\kappa_V$. Since $\|\one_{(Q_{t\log t})^p}\|^2_{L^2(\nu_\rho
\otimes m^p)}=|Q_{t\log t}|^p$ does not increase exponentially fast, it follows from
(\ref{lyap2}), (\ref{combLamb}) and
(\ref{LDRW}--\ref{lb1}) that $\lambda_p \leq \frac{1}{p}\mu_p$.

\medskip\noindent
\underline{Lower bound}:
For every $\delta>0$ there exists an $f_\delta\in L^2(\nu_\rho\otimes m^p)$ such
that
\be{rcspec}
(E_{\mu_p}-E_{\mu_p-\delta}) f_\delta \neq 0
\ee
(see Kato \cite{ka76}, Section VI.2; the spectrum of $G^\kappa_V$ coincides with
the set of $\mu$'s for which $E_{\mu+\delta}-E_{\mu-\delta}\neq 0$ for all
$\delta>0$). Approximating $f_\delta$ by bounded functions, we may without loss
of generality assume that $0 \leq f_\delta \leq 1$. Similarly, approximating
$f_\delta$ by bounded functions with finite support in the spatial variables, we
may assume without loss of generality that there exists a finite
$K_\delta\subset\Z^d$ such that
\be{fdelbds}
0 \leq f_\delta \leq \one_{(K_\delta)^p}.
\ee
First estimate
\be{ub1}
\begin{aligned}
&\E_{\,\nu_\rho,0,\dots,0}
\left(\exp\left[\int_0^t V(Y(s))ds\right]\right)\\
&\quad\geq \sum_{x_1,\dots,x_p \in K_\delta}
\E_{\,\nu_\rho,0,\dots,0}\left(\exp\left[\int_1^t V(Y(s))ds\right]
\one\{X^\kappa_1(1)=x_1,\dots,X^\kappa_p(1)=x_p\}\right)\\
&\quad\geq \sum_{x_1,\dots,x_p \in K_\delta}
\E_{\,\nu_\rho,0,\dots,0}\Bigg(\one\{X^\kappa_1(1)=x_1,\dots,X^\kappa_p(1)=x_p\}\\
&\qquad\qquad\qquad\qquad\qquad \times \E_{\,\xi(1),x_1,\dots,x_p}
\left(\exp\left[\int_0^{t-1} V(Y(s))ds\right]\right)\Bigg)\\
&\quad= \sum_{x_1,\dots,x_p \in K_\delta}
p_1^\kappa(0,x_1) \dots p_1^\kappa(0,x_p)\,
\E_{\,\nu_\rho,x_1,\dots,x_p}\left(\exp\left[\int_0^{t-1} V(Y(s))ds\right]\right)\\
&\quad\geq C_\delta^p\,
\sum_{x_1,\dots,x_p \in K_\delta}
\E_{\,\nu_\rho,x_1,\dots,x_p}\left(\exp\left[\int_0^{t-1} V(Y(s))ds\right]\right),
\end{aligned}
\ee
where $p_t^\kappa(x,y)=\P_x(X^{\kappa}(t)=y)$ and $C_\delta=\min_{x\in K_\delta}
p_1^\kappa(0,x)>0$. The equality in (\ref{ub1}) uses that $\nu_\rho$ is invariant
for the exclusion dynamics.
Next estimate
\be{ub1ext}
\begin{aligned}
\mbox{r.h.s.\ (\ref{ub1})}
&\geq C_\delta^p\,
\int_\Omega \nu_\rho(d\eta) \sum_{x_1,\dots,x_p \in \Z^d}
f_\delta(\eta,x_1,\dots,x_p)\\
&\qquad\qquad\qquad \times \E_{\,\eta,x_1,\dots,x_p}
\left(\exp\left[\int_0^{t-1} V(Y(s))ds\right]\right) f_\delta(Y(t-1))\\
&= C_\delta^p\, (f_\delta,\cP_{t-1} f_\delta)
\geq \frac{C_\delta^p}{|K_\delta|^p}\,
\int_{(\mu_p-\delta,\mu_p]} e^{\mu (t-1)}\,
d\|E_\mu f_\delta\|^2_{L^2(\nu_\rho\otimes m^p)}\\
&\geq C_\delta^p\,
e^{(\mu_p-\delta)(t-1)}\,
\|(E_{\mu_p}-E_{\mu_p-\delta})f_\delta\|^2_{L^2(\nu_\rho\otimes m^p)},
\end{aligned}
\ee
where the first inequality uses (\ref{fdelbds}). Combine (\ref{ub1}--\ref{ub1ext})
with (\ref{rcspec}), and recall (\ref{combLamb}), to get
$\lambda_p\geq \frac{1}{p}(\mu_p-\delta)$. Let $\delta\da 0$, to obtain $\lambda_p
\geq \frac{1}{p}\mu_p$.
\epr

The Rayleigh-Ritz formula for $\mu_p$ applied to Proposition \ref{LyaRR} gives
(recall (\ref{expro1}), (\ref{pot}) and (\ref{gen}--\ref{combgen})):

\bp{RR}
For all $p\in\N$,
\be{RRform}
\lambda_p = \frac{1}{p}\mu_p = \frac{1}{p}\,\,
\sup_{\|f\|_{L^2(\nu_\rho\otimes m^p)}=1} (G^\kappa_V f,f)
\ee
with
\be{decompgen}
(G^\kappa_V f,f) = A_1(f) - A_2(f) - \kappa A_3(f),
\ee
where
\be{Dirwrite}
\begin{aligned}
A_1(f) &= \int_\Omega \nu_\rho(d\eta) \sum_{z_1,\dots,z_p\in\Z^d}
\left[\sum_{i=1}^p \eta(z_i)\right] f(\eta,z_1,\dots,z_p)^2,\\[0.3cm]
A_2(f) &= \int_\Omega \nu_\rho(d\eta) \sum_{z_1,\dots,z_p\in\Z^d}
\frac12 \sum_{\{x,y\}\subset\Z^d}
p(x,y)\,[f(\eta^{x,y},z_1,\dots,z_p)-f(\eta,z_1,\dots,z_p)]^2,\\[0.3cm]
A_3(f) &= \int_\Omega \nu_\rho(d\eta) \sum_{z_1,\dots,z_p\in\Z^d}
\frac12 \sum_{i=1}^p \sum_{{y_i\in\Z^d} \atop {\|y_i-z_i\|=1}}
[f(\eta,z_1,\dots,z_p)|_{z_i\to y_i}-f(\eta,z_1,\dots,z_p)]^2,
\end{aligned}
\ee
and $z_i\to y_i$ means that the argument $z_i$ is replaced by $y_i$.
\ep

We are now ready to give the proof of Theorem \ref{Lyaexist}.

\bpr
The existence of $\lambda_p$ was established in Proposition \ref{LyaRR}. By
(\ref{decompgen}--\ref{Dirwrite}), the r.h.s.\ of (\ref{RRform}) is a supremum
over functions that are linear and non-increasing in $\kappa$. Consequently,
$\kappa\mapsto\lambda_p(\kappa)$ is lower semi-continuous, convex and
non-increasing on $[0,\infty)$ (and, hence, also continuous).
\epr

The variational formula in Proposition \ref{RR} is useful to deduce qualitative properties
of $\lambda_p$, as demonstrated above. Unfortunately, it is not clear how to deduce from it
more detailed information about the Lyapunov exponents. To achieve the latter, we resort
in Sections \ref{S3} and \ref{S4} to different techniques, only occasionally making use of
Proposition \ref{RR}.


\section{Lyapunov exponents: recurrent vs.\ transient random walk}
\label{S3}

In this section we prove Theorems \ref{Lyalow} and \ref{Lyaint}. In Section
\ref{S3.1} we consider recurrent random walk, in Section \ref{S3.2} transient
random walk.

\subsection{Recurrent random walk: proof of Theorem \ref{Lyalow}(i)}
\label{S3.1}

The key to the proof of Theorem \ref{Lyalow}(i) is the following.
\bl{trafic-jam}
If $p(\cdot,\cdot)$ is recurrent, then for any finite box $Q\subset\Z^d$,
\be{stirineq3*}
\lim_{t\to\infty} \frac{1}{t} \log \P_{\nu_\rho}\Bigl(\xi(x,s) = 1\,\,
\forall\,s \in [0,t]\,\,\forall\,x\in Q\Bigr) = 0.
\ee
\el

\bpr
In the spirit of Arratia \cite{arr85}, Section 3, we argue as follows.
Let
\be{stirdef}
H_t^Q = \left\{x\in\Z^d\colon\,\text{there is a path from } (x,0)
\text{ to } Q \times [0,t] \text{ in the graphical representation}\right\}.
\ee


\vskip 0.3truecm

\setlength{\unitlength}{0.3cm}

\begin{picture}(20,10)(-5,1)

  \put(0,0){\line(22,0){22}}
  \put(0,11){\line(22,0){22}}

  \put(2,0){\line(0,12){12}}
  \put(5,0){\line(0,12){12}}
  \put(8,0){\line(0,12){12}}
  \put(11,0){\line(0,12){12}}
  \put(14,0){\line(0,12){12}}
  \put(17,0){\line(0,12){12}}
  \put(20,0){\line(0,12){12}}

  \qbezier[15](2.1,2)(3.5,2)(4.9,2)
  \qbezier[15](5.1,4)(6.5,4)(7.9,4)
  \qbezier[15](8.1,7)(9.5,7)(10.9,7)
  \qbezier[15](2.1,8)(3.5,8)(4.9,8)
  \qbezier[15](11.1,2.5)(12.5,2.5)(13.9,2.5)
  \qbezier[15](17.1,4)(18.5,4)(19.9,4)
  \qbezier[15](11.1,4.5)(12.5,4.5)(13.9,4.5)
  \qbezier[15](14.1,7.5)(15.5,7.5)(16.9,7.5)

  \put(1.7,-1.2){$x$}
  \put(-1,-.3){$0$}
  \put(-1,10.7){$t$}
  \put(10,-1.5){$[$}
  \put(20.9,-1.5){$]$}
  \put(15.2,-1.5){$Q$}
  \put(16.8,-1.5){$\longrightarrow$}
  \put(12.5,-1.5){$\longleftarrow$}

  \put(3,2.5){$\rightarrow$}
  \put(6,4.5){$\rightarrow$}
  \put(9,7.5){$\rightarrow$}
  \put(2.5,.6){$\uparrow$}
  \put(5.5,2.6){$\uparrow$}
  \put(8.5,5.1){$\uparrow$}

  \put(2,0){\circle*{.35}}
  \put(11,7){\circle*{.35}}
  \put(23,0){$\Z^d$}

  \put(-1,-5){\small
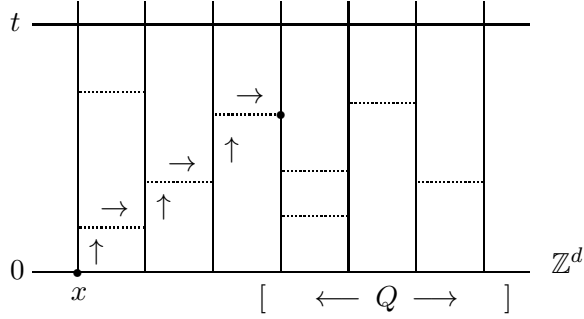
               Fig.\ 4: A path from $(x,0)$ to $Q\times [0,t]$ (recall Fig.\ 1).
             \normalsize}

\end{picture}

\vskip 2.4truecm


\noindent
Note that $H^Q_0=Q$ and that $t\mapsto H^Q_t$ is non-decreasing. Denote by $\cP$ and
$\cE$, respectively, probability and expectation associated with the graphical
representation. Then
\be{stirrep}
\P_{\nu_\rho}\Bigl(\xi(x,s) = 1\,\,
\forall\,s \in [0,t]\,\,\forall\,x\in Q\Bigr)
= (\cP\otimes\nu_\rho)\left(H_t^Q \subseteq \xi(0)\right),
\ee
where $\xi(0)=\{x\in\Z^d\colon\,\xi(x,0)=1\}$ is the set of initial locations of the
particles. Indeed, (\ref{stirrep}) holds because if $\xi(x,0)=0$ for some $x \in H^Q_t$,
then this 0 will propagate into $Q$ prior to time $t$ (see Fig.\ 4).

By Jensen's inequality,
\be{stirineq1}
(\cP\otimes\nu_\rho)\left(H_t^Q \subseteq \xi(0)\right)
= \cE\left(\rho^{|H_t^Q|}\right) \geq \rho^{\cE|H_t^Q|}.
\ee
Moreover, $H_t^Q\subseteq\cup_{y\in Q} H_t^{\{y\}}$, and hence
\be{stirineq2}
\cE|H_t^Q| \leq |Q|\, \cE|H_t^{\{0\}}|.
\ee
Furthermore, we have
\be{range1}
\cE|H_t^{\{0\}}|=\E_{\,0}^{\,p(\cdot,\cdot)}R_t,
\ee
where $R_t$ is the range after time $t$ of the random walk with transition kernel
$p(\cdot,\cdot)$ driving $\xi$ and $\E_{\,0}^{\,p(\cdot,\cdot)}$ denotes expectation
w.r.t.\;this random walk starting from $0$. Indeed, by time reversal, the probability
that there is a path from $(x,0)$ to $\{0\}\times [0,t]$ in the graphical representation
is equal to the probability that the random walk starting from 0 hits $x$ prior to time
$t$. It follows from (\ref{stirrep}--\ref{range1}) that
\be{stirineq3}
\frac{1}{t} \log \P_{\nu_\rho}\Bigl(\xi(x,s) = 1\,\,
\forall\,s \in [0,t]\,\,\forall\,x\in Q\Bigr)
\geq - |Q|\log\left(\frac{1}{\rho}\right)\left\{\frac{1}{t}
\E_{\,0}^{\,p(\cdot,\cdot)}R_t\right\}.
\ee
Finally, since $\lim_{t\to\infty}\frac{1}{t}\E_{\,0}^{\,p(\cdot,\cdot)}R_t=0$ when
$p(\cdot,\cdot)$ is recurrent (see Spitzer \cite{sp76}, Chapter 1, Section 4), we get
(\ref{stirineq3*}).
\epr

We are now ready to give the proof of Theorem \ref{Lyalow}(i).

\bpr
Since $p\mapsto\lambda_p$ is non-decreasing and $\lambda_p \leq 1$ for all $p\in\N$,
it suffices to give the proof for $p=1$. For $p=1$, (\ref{fey-kac2}) gives
\be{L1a}
\Lambda_1(t) = \frac{1}{t} \log
\E_{\,\nu_\rho,0}\left(\exp\left[\int_0^t
\xi\left(X^\kappa(s),s\right)ds\right]\right).
\ee
By restricting $X^\kappa$ to stay inside a finite box $Q\subset\Z^d$ up to
time $t$ and requiring $\xi$ to be 1 throughout this box up to time $t$, we
obtain
\be{L1b}
\begin{aligned}
&\E_{\,\nu_\rho,0}\left(\exp\left[\int_0^t
\xi(X^\kappa(s),s)\, ds\right]\right)\\
&\geq e^t \,\P_{\nu_\rho}\Bigl(\xi(x,s) = 1\,\,
\forall\,s \in [0,t]\,\,\forall\,x\in Q\Bigr)
\,\PS_0\Bigl(X^\kappa(s)\in Q\,\,\forall\,s\in [0,t] \Bigr).
\end{aligned}
\ee
For the first factor, we apply (\ref{stirineq3*}). For the second factor, we have
\be{L1d}
\lim_{t\to\infty} \frac{1}{t} \log
\PS_0\Bigl(X^\kappa(s)\in Q\,\,\forall\,s\in [0,t] \Bigr)
= -\lambda^\kappa(Q)
\ee
with $\lambda^\kappa(Q)>0$ the principal Dirichlet eigenvalue on $Q$ of $-\kappa\Delta$,
the generator of the simple random walk $X^\kappa$.
Combining (\ref{stirineq3*}) and (\ref{L1a}--\ref{L1d}), we arrive at
\be{resSE}
\lambda_1 = \lim_{t\to\infty} \Lambda_1(t) \geq 1-\lambda^\kappa(Q).
\ee
Finally, let $Q\to\Z^d$ and use that $\lim_{Q\to\Z^d} \lambda^\kappa(Q)=0$ for any $\kappa$,
to arrive at $\lambda_1\geq 1$. Since, trivially, $\lambda_1\leq 1$, we get $\lambda_1=1$.
\epr

\subsection{Transient random walk: proof of Theorems \ref{Lyalow}(ii) and \ref{Lyaint}}
\label{S3.2}

Theorem \ref{Lyalow}(ii) is proved in Sections \ref{S3.2.1} and \ref{S3.2.2}--\ref{S3.2.4},
Theorem \ref{Lyaint} in Section \ref{S3.2.5}. In Section \ref{S3.2.6} we make a link between
Section \ref{S2.1} and Proposition \ref{RR} for $\kappa=0$ that provides further background
for Theorem \ref{Lyaint}. Throughout the present section we assume that the random walk kernel 
$p(\cdot,\cdot)$ is transient.

\subsubsection{Proof of the lower bound in Theorem \ref{Lyalow}(ii)}
\label{S3.2.1}

\bp{lem-kap0}
$\lambda_p(\kappa)>\rho$ for all $\kappa\in[0,\infty)$ and $p\in\N$.
\ep

\bpr
Since $p\mapsto\lambda_p(\kappa)$ is non-decreasing for all $\kappa$, it suffices to
give the proof for $p=1$. For every $\epsilon>0$ there exists a function
$\phi_\epsilon\colon\Z^d\to\R$ such that
\be{phicond}
\sum_{x\in\Z^d} \phi_\epsilon(x)^2=1 \quad \mbox{ and } \quad
\sum_{{x,y\in\Z^d} \atop {\|x-y\|=1}}\,
[\phi_\epsilon(x)-\phi_\epsilon(y)]^2 \leq \epsilon^2.
\ee
Let
\be{fepsdef}
f_\epsilon(\eta,x) = \frac{1+\epsilon\eta(x)}{[1+(2\epsilon+\epsilon^2)\rho]^{1/2}}\,
\phi_\epsilon(x), \qquad \eta \in \Omega,\,x\in\Z^d.
\ee
Then
\be{normcheck}
\|f_\epsilon\|_{L^2(\nu_\rho \otimes m)}^2
= \int_\Omega \nu_\rho(d\eta)\,\sum_{x\in\Z^d}
\frac{[1+\epsilon\eta(x)]^2}{1+(2\epsilon+\epsilon^2)\rho}\,\phi_\epsilon(x)^2
= \sum_{x\in\Z^d} \phi_\epsilon(x)^2
= 1.
\ee
Therefore we may use $f_\epsilon$ as a test function in (\ref{RRform}) in Proposition
\ref{RR}. This gives
\be{lbtest}
\lambda_1 = \mu_1 \geq \frac{1}{1+(2\epsilon+\epsilon^2)\rho}\, (I-II-\kappa\,III)
\ee
with
\be{termI}
I = \int_\Omega \nu_\rho(d\eta)\sum_{z\in\Z^d}
\eta(z)\, [1+\epsilon\eta(z)]^2\, \phi_\epsilon(z)^2
= (1+2\epsilon+\epsilon^2)\rho\,
\sum_{z\in\Z^d} \phi_\epsilon(z)^2
= (1+2\epsilon+\epsilon^2)\rho
\ee
and
\be{termII}
\begin{aligned}
II &= \int_\Omega \nu_\rho(d\eta) \sum_{z\in\Z^d} \frac14 \sum_{x,y\in\Z^d}
p(x,y)\,\epsilon^2 [\eta^{x,y}(z)-\eta(z)]^2\phi_\epsilon(z)^2\\
&= \frac12\int_\Omega \nu_\rho(d\eta)\,\,\sum_{x,y\in\Z^d} p(x,y)\,
\epsilon^2 [\eta(x)-\eta(y)]^2\, \phi_\epsilon(x)^2\\
&= \epsilon^2\rho(1-\rho)\,\sum_{{x,y\in\Z^d} \atop {x\neq y}}
p(x,y)\,\phi_\epsilon(x)^2
\leq  \epsilon^2\rho(1-\rho)
\end{aligned}
\ee
and
\be{termIII}
\begin{aligned}
III &= \frac12 \int_\Omega \nu_\rho(d\eta) \sum_{{x,y\in\Z^d} \atop {\|x-y\|=1}}
\big\{[1+\epsilon\eta(x)]\phi_\epsilon(x)
-[1+\epsilon\eta(y)]\phi_\epsilon(y)\big\}^2\\
&= \frac12 \sum_{{x,y\in\Z^d} \atop {\|x-y\|=1}} \big\{[1+(2\epsilon+\epsilon^2)\rho]
[\phi_\epsilon(x)^2+\phi_\epsilon(y)^2]
-2(1+\epsilon\rho)^2\phi_\epsilon(x)\phi_\epsilon(y) \big\}\\
&= \frac12[1+(2\epsilon+\epsilon^2)\rho] \sum_{{x,y\in\Z^d} \atop {\|x-y\|=1}}
[\phi_\epsilon(x)-\phi_\epsilon(y)]^2 + \epsilon^2 \rho(1-\rho)
\sum_{{x,y\in\Z^d} \atop {\|x-y\|=1}}
\phi_\epsilon(x)\phi_\epsilon(y)\\
&\leq \frac12 [1+(2\epsilon+\epsilon^2)\rho]\epsilon^2
+ 2d\epsilon^2 \rho(1-\rho).
\end{aligned}
\ee
In the last line we use that $\phi_\epsilon(x)\phi_\epsilon(y)\leq \frac{1}{2}\phi_\epsilon(x)^2
+\frac{1}{2}\phi_\epsilon(y)^2$. Combining (\ref{lbtest}--\ref{termIII}), we find
\be{lbtestres}
\lambda_1 = \mu_1 \geq \rho\,\frac{1+2\epsilon+O(\epsilon^2)}{1+2\epsilon\rho+O(\epsilon^2)}.
\ee
Because $\rho\in(0,1)$, it follows that for $\epsilon$ small enough the r.h.s.\
is strictly larger than $\rho$.
\epr

\subsubsection{Proof of Theorem \ref{Lyaint}}
\label{S3.2.5}

\bpr
For $\kappa=0$, (\ref{fey-kac2}) reduces to
\be{lyap1}
\Lambda_p(t) = \frac{1}{pt} \log
\E_{\,\nu_\rho}\left(\exp\left[p \int_0^t \xi(0,s) ds\right]\right)
= \frac{1}{pt} \log
\E_{\,\nu_\rho}\left(\exp\left[pT_t\right]\right)
\ee
(recall (\ref{occtime})). In order to compute $\lambda_p(0)=\lim_{t\to\infty}
\Lambda_p(t)$, we may use the large deviation principle for $(T_t)_{t\geq 0}$
cited in Section \ref{S2.1} due to Landim \cite{lan92}. Indeed, by applying
Varadhan's Lemma (see e.g.\ den Hollander \cite{hol00}, Theorem III.13) to
(\ref{lyap1}), we get
\be{Varocc}
\lambda_p(0) = \frac{1}{p}
\max_{\alpha \in [0,1]} \big[p\alpha - \Psi_d(\alpha)\big]
\ee
with $\Psi_d$ the rate function given in (\ref{rfSEa}). Since $\Psi_d$ is lower
semi-continuous, (\ref{Varocc}) has at least one maximizer $\alpha_p$:
\be{inc5}
\lambda_p(0) = \alpha_p - \frac{1}{p} \Psi_d(\alpha_p).
\ee
By Proposition \ref{lem-kap0} for $\kappa=0$, we have $\lambda_p(0)>\rho$.
Hence $\alpha_p>\rho$ (because $\Psi_d(\rho)=0$). Since $p(\cdot,\cdot)$ is
transient, we may use the quadratic lower bound in (\ref{rflb}) to see that
$\Psi_d(\alpha_p)>0$. Therefore we get from (\ref{Varocc}--\ref{inc5}) that
\be{inc9}
\lambda_{p+1}(0)
\geq \frac{1}{p+1}\left[\alpha_p(p+1)-\Psi_d(\alpha_p)\right]
= \alpha_p - \frac{1}{p+1} \Psi_d(\alpha_p)
>\alpha_p-\frac1p \Psi_d(\alpha_p)
= \lambda_p(0).
\ee
Since $p$ is arbitrary, this completes the proof of Theorem \ref{Lyaint}.
\epr

\subsubsection{Proof of the upper bound in Theorem \ref{Lyalow}(ii)}
\label{S3.2.2}

\bp{lem_kap0*}
$\lambda_p(\kappa)<1$ for all $\kappa\in[0,\infty)$ and $p\in\N$.
\ep

\bpr
By Theorem \ref{Lyaint}, which was proved in Section \ref{S3.2.5}, we know that
$p\mapsto\lambda_p(0)$ is strictly increasing. Since $\lambda_p(0) \leq 1$ for
all $p\in\N$, it therefore follows that $\lambda_p(0)<1$ for all $p\in\N$.
Moreover, by Theorem \ref{Lyaexist}(ii), which was proved in Section \ref{S2.2},
we know that $\kappa\mapsto\lambda_p(\kappa)$ is non-increasing. It therefore
follows that $\lambda_p(\kappa)<1$ for all $\kappa\in[0,\infty)$ and $p\in\N$.
\epr

\subsubsection{Proof of the asymptotics in Theorem \ref{Lyalow}(ii)}
\label{S3.2.3}

The proof of the next proposition is somewhat delicate.

\bp{lem-kaplim}
$\lim_{\kappa \to\infty}\lambda_p(\kappa)=\rho$ for all $p\in\N$.
\ep

\bpr
We give the proof for $p=1$. The generalization to arbitrary $p$ is straightforward
and will be explained at the end. We need a cube $Q=[-R,R]^d\cap\Z^d$ of length $2R$,
centered at the origin and $\delta\in(0,1)$. Limits are taken in the order
\be{lim-ord}
t\to\infty,\,\, \kappa\to\infty,\,\, \delta\downarrow 0,\,\, Q\uparrow\Z^d.
\ee
The proof proceeds in 4 steps, each containing a lemma.

\noindent
\underline{Step 1:} Let $X^{\kappa,Q}$ be simple random walk on $Q$ obtained from
$X^\kappa$ by suppressing jumps outside of $Q$. Then $(\xi_t,X_t^{\kappa,Q})_{t\geq 0}$ is
a Markov process on $\Omega\times Q$ with self-adjoint generator in $L^2(\nu_\rho\otimes m_Q)$,
where $m_Q$ is the counting measure on $Q$.

\bl{kaplim-lem1}
For all $Q$ finite (centered and cubic) and $\kappa \in [0,\infty)$,
\be{kaplim-3}
\E_{\,\nu_\rho,0}\left(\exp\left[\int_0^t ds\, \xi(X_s^\kappa,s)\right]\right)
\leq e^{o(t)}\,
\E_{\,\nu_\rho,0}\left(\exp\left[\int_0^t ds\, \xi\big(X_s^{\kappa,Q},s\big)\right]\right),\,\,
t\to\infty.
\ee
\el

\bpr
We consider the partition of $\Z^d$ into cubes $Q_z=2Rz+Q$, $z\in\Z^d$. The Lyapunov exponent
$\lambda_1(\kappa)$ associated with $X^\kappa$ is given by the variational formula
(\ref{RRform}--\ref{Dirwrite}) for $p=1$. It can be estimated from above by splitting the sums
over $\Z^d$ in (\ref{Dirwrite}) into separate sums over the individual cubes $Q_z$ and
suppressing in $A_3(f)$ the summands on pairs of lattice sites belonging to different cubes.
The resulting expression is easily seen to coincide with the original variational expression
(\ref{RRform}), except that the supremum is restricted in addition to functions $f$ with
spatial support contained in $Q$. But this is precisely the Lyapunov exponent $\lambda_1^Q(\kappa)$
associated with $X^{\kappa,Q}$. Hence, $\lambda_1(\kappa)\leq \lambda_1^Q(\kappa)$, and this
implies (\ref{kaplim-3}).
\epr

\noindent
\underline{Step 2:} For large $\kappa$ the random walk $X^{\kappa,Q}$ moves fast through
the finite box $Q$ and therefore samples it in a way that is close to the uniform distribution.

\bl{kaplim-lem2}
For all $Q$ finite and $\delta\in (0,1)$, there exist $\varepsilon=\varepsilon(\kappa,\delta,Q)$
and $N_0=N_0(\delta,\varepsilon)$, satisfying $\lim_{\kappa\to\infty}\varepsilon(\kappa,\delta,Q)=0$
and $\lim_{\delta,\varepsilon\downarrow 0}N_0(\delta,\varepsilon)=N_0>1$, such that
\be{kaplim-7}
\begin{aligned}
\E_{\,\nu_\rho,0}\left(\exp\left[\int_0^t ds\, \xi\big(X_s^{\kappa,Q},s\big)\right]\right)
&\leq o(1)+ \exp\bigg[\bigg(\Big(1+\frac{1+\varepsilon}{1-\delta}\Big)\delta N_0|Q|
+\frac{\delta+\varepsilon}{1-\delta}\bigg)(t+\delta)\bigg]\\
&\qquad\qquad\times
\E_{\,\nu_\rho}\bigg(\exp\bigg[\int_0^{t+\delta} ds\,
\frac{1}{|Q|}\sum_{y\in Q}\xi(y,s)\bigg]\bigg),\,\,t\to\infty.
\end{aligned}
\ee
\el

\bpr
We split time into intervals of length $\delta>0$. Let $I_k$ be the indicator of the event
that $\xi$ has a jump time in $Q$ during the {\it k}-th time interval. If $I_k=0$, then
$\xi_s=\xi_{(k-1)\delta}$ for all $s\in[(k-1)\delta,k\delta)$. Hence,
\be{kaplim-13}
\int_{(k-1)\delta}^{k\delta} ds\, \xi_s\big(X_s^{\kappa,Q}\big)
\leq \int_{(k-1)\delta}^{k\delta}ds\, \xi_{(k-1)\delta}\big(X_s^{\kappa,Q}\big)
+\delta I_k
\ee
and, consequently, we have for all $x\in\Z^d$ and $k=1,\dots,\lceil t/\delta\rceil$,
\be{kaplim-9}
\begin{aligned}
\ES_{\,x}\bigg(\exp\bigg[\int_0^\delta ds\, \xi_{(k-1)\delta+s}\big(X_s^{\kappa,Q}\big)\bigg]\bigg)
&\leq e^{\delta I_k}\,
\ES_{\,x}\bigg(\exp\bigg[\int_0^\delta ds\, \eta\big(X_s^{\kappa,Q}\big)\bigg]\bigg),
\end{aligned}
\ee
where we abbreviate $\xi_{(k-1)\delta}=\eta$. Next, we do a Taylor expansion and use
the Markov property of $X^{\kappa,Q}$, to obtain ($s_0=0$)
\be{kaplim-21}
\begin{aligned}
&\ES_{\,x}\left(\exp\left[\int_0^\delta ds\, \eta\big(X_s^{\kappa,Q}\big)\right]\right)
=\sum_{n=0}^{\infty} \left(\prod_{l=1}^{n}\int_{s_{l-1}}^\delta ds_l\right)
\ES_{x}\left(\prod_{m=1}^{n} \eta\big(X_{s_m}^{\kappa,Q}\big)\right)\\
&\quad\leq \sum_{n=0}^{\infty} \left(\prod_{l=1}^{n}\int_{s_{l-1}}^\delta ds_l\right)
\left(\prod_{m=1}^n \max_{x\in Q}\ES_{x}\left(\eta\big(X_{s_m-s_{m-1}}^{\kappa,Q}\big)\right)\right)\\
&\quad\leq \sum_{n=0}^{\infty} \left\{\int_0^\delta ds\, \max_{x\in Q}
\ES_{\,x}\left(\eta\big(X_s^{\kappa,Q}\big)\right)\right\}^n
\leq\exp\bigg[\frac{1}{1-\delta}\int_0^\delta ds\, \max_{x\in Q}
\ES_{\,x}\Big(\eta\big(X_s^{\kappa,Q}\big)\Big)\bigg]\\
&\quad\leq\exp\Bigg[\frac{1}{1-\delta}\sum_{y\in Q}\eta(y) \int_0^\delta ds\, \max_{x\in Q}
\ES_{\,x}\Big(\delta_y\big(X_s^{\kappa,Q}\big)\Big)\Bigg],
\end{aligned}
\ee
where we use that $\max_{x\in Q}
\E_{\,x}\left(\int_0^\delta ds\, \eta\big(X_s^{\kappa,Q}\big)\right)\leq \delta$.
Now, let $p_s^{\kappa,Q}(\cdot,\cdot)$ denote the transition kernel of $X^{\kappa,Q}$.
Note that
\be{kaplim-25}
\lim_{\kappa\to\infty} p_s^{\kappa,Q}(x,y) = \frac{1}{|Q|}
\quad \text{for all } s>0,\,\,  Q \text{ finite} \text{ and } x,y\in Q.
\ee
Hence
\be{kaplim-26}
\lim_{\kappa\to\infty}\ES_x \Big(\delta_y\big(X_s^{\kappa,Q}\big)\Big)
=\frac{1}{|Q|} \quad \text{for all } s>0,\,\,  Q \text{ finite} \text{ and } x,y\in Q.
\ee
Therefore, by the Lebesgue dominated convergence theorem, we have
\be{kaplim-27}
\lim_{\kappa\to\infty}\int_0^\delta ds\, \max_{x\in Q}\ES_{\,x}\Big(
\delta_y\big(X_s^{\kappa,Q}\big)\Big) = \delta\, \frac{1}{|Q|}
\quad \text{for all } \delta>0,\,\,  Q \text{ finite}\text{ and } y\in Q.
\ee
This implies that the expression in the exponent in the r.h.s.\ of (\ref{kaplim-21})
converges to
\be{kaplim-28}
\frac{\delta}{1-\delta}\frac1{|Q|} \sum_{y\in Q}\eta(y),
\ee
uniformly in $\eta\in\Omega$. Combining the latter with (\ref{kaplim-21}), we see that
there exists some $\varepsilon=\varepsilon(\kappa,\delta,Q)$, satisfying
$\lim_{\kappa\to\infty}\varepsilon(\kappa,\delta,Q)=0$, such that for all $x\in Q$,
\be{kaplim-29}
\ES_{\,x}\bigg(\exp\bigg[\int_0^\delta ds\, \eta\big(X_s^{\kappa,Q}\big)\bigg]\bigg)
\leq \exp\left[\frac{1+\varepsilon}{1-\delta}\, \delta
\, \frac{1}{|Q|}\sum_{y\in Q}\eta(y)\right]
\quad \text{for all } \delta \in (0,1) \text{ and } Q \text{ finite}.
\ee
Next, similarly as in (\ref{kaplim-13}), we have
\be{kaplim-31}
\delta\, \frac{1}{|Q|}\sum_{y\in Q} \xi_{(k-1)\delta}(y)
\leq \int_{(k-1)\delta}^{k\delta}ds\, \frac{1}{|Q|}\sum_{y\in Q} \xi_s(y)
+\delta I_k.
\ee
Applying the Markov property to $X^{\kappa,Q}$, and using (\ref{kaplim-13}) and
(\ref{kaplim-29}-\ref{kaplim-31}), we find that
\be{kaplim-35}
\begin{aligned}
\E_{\,\nu_\rho,0}\left(\exp\left[\int_0^t ds\, \xi\big(X_s^{\kappa,Q},s\big)\right]\right)
&\leq \E_{\,\nu_\rho}\Bigg(\exp\bigg[
\Big(1+\frac{1+\varepsilon}{1-\delta}\Big)\delta\, N_{t+\delta}
+\frac{\delta+\epsilon}{1-\delta}(t+\delta)\bigg]\\
&\qquad\qquad\qquad\times
\exp\Bigg[\int_0^{t+\delta}ds\,
\frac1{|Q|}\sum_{y\in Q} \xi_s(y)\Bigg]\Bigg),
\end{aligned}
\ee
where $N_{t+\delta}$ is the total number of jumps that $\xi$ makes inside $Q$ up 
to time $t+\delta$. The second term in the r.h.s.\ of (\ref{kaplim-35}) equals the 
second term in the r.h.s.\ of (\ref{kaplim-7}). The first term will be negligible on 
an exponential scale for $\delta\downarrow 0$, because, as can be seen from the graphical 
representation, $N_{t+\delta}$ is stochastically smaller that the total number of jumps 
up to time $t+\delta$ of a Poisson process with rate $|Q\cup\partial Q|$. Indeed, 
abbreviating
\be{abJdef}
a = \left(1+\frac{1+\varepsilon}{1-\delta}\right)\delta, \quad
b = \frac{\delta+\varepsilon}{1-\delta}, \quad
M_{t+\delta} = \int_0^{t+\delta}ds\,\,\frac{1}{|Q|} \sum_{y\in Q} \xi_s(y),
\ee
we estimate, for each $N$,
\be{kaplim-39}
\begin{aligned}
&{\rm r.h.s.} \,\,(\ref{kaplim-35})
=\E_{\,\nu_\rho}\left(e^{aN_{t+\delta}+b(t+\delta)+M_{t+\delta}}\right)\\
&\qquad\leq e^{(b+1)(t+\delta)}\, \E_{\,\nu_\rho}\Big(e^{aN_{t+\delta}}\,
1\{N_{t+\delta}\geq N|Q|(t+\delta)\}\Big)
+e^{(aN|Q|+b)(t+\delta)}\, \E_{\,\nu_\rho}\left(e^{M_{t+\delta}}\right).
\end{aligned}
\ee
For $N\geq N_0=N_0(a,b)$, the first term tends to zero as $t\to\infty$ and
can be discarded. Hence
\be{kaplim-41}
{\rm r.h.s.} \,\,(\ref{kaplim-35})
\leq e^{(aN_0|Q|+b)(t+\delta)}\, \E_{\,\nu_\rho}\left(e^{bM_{t+\delta}}\right),
\ee
which is the desired bound in (\ref{kaplim-7}). Note that $a\downarrow 0$, $b\downarrow 1$
as $\delta, \varepsilon \downarrow 0$ and hence $N_0(a,b)\downarrow N_0>1$.
\epr

\noindent
\underline{Step 3:} By combining Lemmas \ref{kaplim-lem1}--\ref{kaplim-lem2},
we now know that for any $Q$ finite,
\be{kaplim-45}
\lim_{\kappa\to\infty}\lambda_1(\kappa)
\leq \lim_{t\to\infty}\frac1t\log\E_{\,\nu_\rho}
\Bigg(\exp\Bigg[\int_0^t ds\, \frac1{|Q|}\sum_{y\in Q}\xi_s(y)\Bigg]\Bigg),
\ee
where we have taken the limits $\kappa\to\infty$ and $\delta\downarrow 0$.
According to Proposition \ref{IRW-comp} (with $K(z,s)=(1/|Q|)1_Q(z)$),
\be{kaplim-47}
\E_{\,\nu_\rho}
\Bigg(\exp\Bigg[\int_0^t ds\, \frac1{|Q|}\sum_{y\in Q}\xi_s(y)\Bigg]\Bigg)
\leq \E_{\,\nu_\rho}^{\IRW}
\Bigg(\exp\Bigg[\int_0^t ds\, \frac1{|Q|}\sum_{y\in Q}\tilde\xi_s(y)\Bigg]\Bigg),
\ee
where $(\tilde\xi_{t})_{t\geq 0}$ is the process of Independent Random Walks on
$\Z^d$ with step rate 1 and transition kernel $p(\cdot,\cdot)$, and $\E_{\,\nu_\rho}^{\IRW} 
=\int_\Omega \nu_\rho(d\eta)\,\E_{\,\eta}^{\IRW}$. The r.h.s.\ can be computed
and estimated as follows. Write
\be{pgen}
(\Delta^{(p)}f)(x) = \sum_{y\in\Z^d} p(x,y)[f(y)-f(x)], \qquad x\in\Z^d,
\ee
to denote the generator of the random walk with step rate 1 and transition kernel 
$p(\cdot,\cdot)$.

\bl{kaplim-lem3}
For all $Q$ finite,
\be{kaplim-49}
{\rm r.h.s.}\,\,(\ref{kaplim-47})
\leq e^{\rho t}\,\exp\left[\int_0^t ds\, \frac1{|Q|}\sum_{x\in Q}w^Q(x,s)\right],
\ee
where $w^Q\colon\,\Z^d\times[0,\infty)\to \R$ is the solution of the Cauchy problem
\be{kaplim-50}
\frac{\partial w^Q}{\partial t}(x,t) = \Delta^{(p)} w^Q(x,t)
+ \left\{\frac1{|Q|} 1_{Q}(x)\right\}[w^Q(x,t)+1],
\qquad
w^Q(\cdot,0)\equiv 0,
\ee
which has the representation
\be{kaplim-51}
w^Q(x,t)=\E_{\, x}^{\RW}\left(\exp\left[\int_0^t ds\, \frac1{|Q|}1_{Q}(Y_s)\right]\right)
-1\geq 0,
\ee
where $Y=(Y_t)_{t\geq 0}$ is the single random walk with step rate 1 and transition kernel 
$p(\cdot,\cdot)$, and $\E_x^{\RW}$ denotes the expectation w.r.t.\ to $Y$ starting from
$Y_0=x$.
\el

\bpr
Let
\be{kaplim-51.3}
A_\eta = \{x\in\Z^d\colon\, \eta(x)=1\}, \qquad \eta\in\Omega.
\ee
Then
\be{kaplim-51.5}
\begin{aligned}
{\rm r.h.s.} \,\,(\ref{kaplim-47})
&= \int_\Omega \nu_\rho(d\eta)\, \E_{\,\eta}^{\IRW}\left(\exp\bigg[
\int_0^t ds\, \frac1{|Q|} \sum_{x\in A_\eta} \sum_{y\in Q}
1_y(\tilde\xi_{s,x})\bigg]\right)\\
&= \int_\Omega \nu_\rho(d\eta)\, \prod_{x\in A_\eta}\E_{\,x}^{\RW}
\left(\exp\bigg[\int_0^t ds\, \frac1{|Q|}\, 1_Q(Y_s)
\bigg]\right),
\end{aligned}
\ee
where $\tilde\xi_{s,x}$ is the position at time $s$ of the random walk starting from
$\tilde\xi_{0,x}=x$ (in the process of Independent Random Walks
$\tilde\xi=(\tilde\xi_t)_{t\geq 0}$). Let
\be{kaplim-51.9}
v^Q(x,t)
=\E_{\,x}^{\RW}
\left(\exp\bigg[\int_0^t ds\, \frac1{|Q|}\, 1_Q(Y_s)
\bigg]\right).
\ee
By the Feynman-Kac formula, $v^Q(x,t)$ is  the solution of the Cauchy problem
\be{kaplim-51.11}
\frac{\partial v^Q}{\partial t}(x,t) = \Delta^{(p)} v^Q(x,t)
+ \left\{\frac1{|Q|} 1_Q(x)\right\}\, v^Q(x,t),
\qquad v^Q(\cdot,0) \equiv 1.
\ee
Now put
\be{kaplim-51.13}
w^Q(x,t)=v^Q(x,t)-1.
\ee
Then (\ref{kaplim-51.11}) can be rewritten as (\ref{kaplim-50}). Combining
(\ref{kaplim-51.5}--\ref{kaplim-51.9}) and (\ref{kaplim-51.13}), we get
\be{kaplim-51.15}
\begin{aligned}
{\rm r.h.s.} \,\,(\ref{kaplim-47})
&= \int_\Omega \nu_\rho(d\eta)\, \prod_{x\in A_\eta}
\big(1+w^Q(x,t)\big)
= \int_\Omega \nu_\rho(d\eta)\, \prod_{x\in \Z^d}
\big(1+\eta(x)\, w^Q(x,t)\big)\\
&= \prod_{x\in \Z^d}
\big(1+\rho\, w^Q(x,t)\big)
\leq \exp\left[\rho\sum_{x\in \Z^d}
w^Q(x,t)\right],
\end{aligned}
\ee
where we use that $\nu_\rho$ is the Bernoulli product measure with density $\rho$.
Summing (\ref{kaplim-50}) over $\Z^d$, we have
\be{kaplim-51.17}
\frac{\partial}{\partial t}\Bigg(\sum_{x\in\Z^d} w^Q(x,t)\Bigg)
=\sum_{x\in Q} \frac1{|Q|}\, w^Q(x,t)+1.
\ee
Integrating (\ref{kaplim-51.17}) w.r.t.\ time, we get
\be{kaplim-51.19}
\sum_{x\in\Z^d}w^Q(x,t)
=\int_0^t ds\, \sum_{x\in Q} \frac1{|Q|}\, w^Q(x,s)+t.
\ee
Combining (\ref{kaplim-51.15}) and(\ref{kaplim-51.19}), we get the claim.
\epr

\noindent
\underline{Step 4:} The proof is completed by showing the following:
\bl{kaplim-lem4}
\be{kaplim-53}
\lim_{Q\uparrow\Z^d}\lim_{t\to\infty}\frac1t \int_0^t ds\,
\frac1{|Q|}\sum_{x\in Q} w^Q(x,s)=0.
\ee
\el

\bpr
Let $\cG$ denote the Green operator acting on functions $V\colon\,\Z^d\to[0,\infty)$ as
\be{kaplim-54}
(\cG V)(x)=\sum_{y\in\Z^d}G(x,y)V(y),\qquad x\in\Z^d,
\ee
where $G(x,y)=\int_0^\infty dt\, p_t(x,y)$ denotes the Green kernel on $\Z^d$.
We have
\be{kaplim-55}
\bigg\|\cG\Big(\frac1{|Q|}1_{|Q|}\Big)\bigg\|_{\infty}
=\sup_{x\in\Z^d}\sum_{y\in Q} G(x,y) \frac1{|Q|}.
\ee
The r.h.s.\ tends to zero as $Q\uparrow\Z^d$, because $G(x,y)$ tends to
zero as $\|x-y\|\to\infty$. Hence Lemma 8.2.1 in G\"artner and den Hollander
\cite{garhol04} can be applied to (\ref{kaplim-51}) for $Q$ large enough, to yield
\be{kaplim-57}
\sup_{{x\in\Z^d}\atop{s\geq 0}} w^Q(x,s) \leq \varepsilon(Q) \downarrow 0\quad
\text{as}\quad Q\uparrow\Z^d,
\ee
which proves (\ref{kaplim-53}).
\epr

Combine (\ref{kaplim-45}--\ref{kaplim-47}), (\ref{kaplim-49}) and (\ref{kaplim-53})
to get the claim in Proposition \ref{lem-kaplim}.

This completes the proof of Proposition \ref{lem-kaplim} for $p=1$. The generalization to
arbitrary $p$ is straightforward and runs as follows. Return to (\ref{fey-kac2}). Separate
the $p$ terms under the sum with the help of H\"older's inequality with weights $1/p$.
Next, use (\ref{kaplim-3}) for each of the $p$ factors, leading to $\frac1p\log$ of the
r.h.s.\ of (\ref{kaplim-3}) with an extra factor $p$ in the exponent. Then proceed as before,
which leads to Lemma \ref{kaplim-lem3} but with $w^Q$ the solution of (\ref{kaplim-50}) with
$\frac{p}{|Q|}1_Q(x)$ between braces. Then again proceed as before, which leads to
(\ref{kaplim-51.15}) but with an extra factor $p$ in the r.h.s.\ of (\ref{kaplim-51.19}).
The latter gives a factor $e^{p\rho t}$ replacing $e^{\rho t}$ in (\ref{kaplim-49}).
Now use Lemma \ref{kaplim-lem4} to get the claim.
\epr

\subsubsection{Proof of the strict monotonicity in Theorem \ref{Lyalow}(ii)}
\label{S3.2.4}

By Theorem \ref{Lyaexist}(ii), $\kappa\mapsto\lambda_p(\kappa)$ is convex.
Because of Proposition \ref{lem-kap0} and Proposition \ref{lem-kaplim}, it must be
strictly decreasing. This completes the proof of Theorem \ref{Lyalow}(ii).

\subsubsection{Relation between Proposition \ref{RR} and (\ref{Varocc})}
\label{S3.2.6}

Let $\kappa=0$ and $p=1$. The generalization to arbitrary $p$ is straightforward.

For $\kappa=0$ and $p=1$, (\ref{RRform}--\ref{Dirwrite}) in Proposition \ref{RR}
read
\be{varcomp1}
\lambda_1(0) = \sup_{\|f\|_{L^2(\nu_\rho\otimes m)}=1} \sum_{z\in\Z^d}
\left[\int_\Omega \nu_\rho(d\eta)\eta(z)\,f(\eta,z)^2
- \cE(f(\cdot,z))\right],
\ee
where we recall (\ref{Diridef}). Split the supremum into two parts,
\be{varcomp2}
\lambda_1(0) = \sup_{\|g\|_{L^2(m)}=1}\,\,
\sup_{\|f_z\|_{L^2(\nu_\rho)}=1\,\forall\,z\in\Z^d}
\sum_{z\in\Z^d} g^2(z)
\left[\int_\Omega \nu_\rho(d\eta)\,\eta(z)\,f_z(\eta)^2
- \cE(f_z)\right],
\ee
where $f_z(\eta)=f(\eta,z)/g(z)$ with $g(z)^2=\int_\Omega\nu_\rho(d\eta)f(\eta,z)^2$.
The second supremum in (\ref{varcomp2}), which runs over a family of functions indexed
by $z$, can be brought under the sum,
\be{varcomp3}
\lambda_1(0) =  \sup_{\|g\|_{L^2(m)}=1} \sum_{z\in\Z^d} g^2(z)
\sup_{\|f_z\|_{L^2(\nu_\rho)}=1}
\left[\int_\Omega \nu_\rho(d\eta)\,\eta(z)\,f_z(\eta)^2
- \cE(f_z)\right].
\ee
By the shift-invariance of $\nu_\rho$, we may replace $\eta(z)$ by $\eta(0)$ under the
second supremum in (\ref{varcomp3}), in which case the latter no longer depends on
$z$, and we get
\be{varcomp4}
\begin{aligned}
\lambda_1(0) &= \sup_{\|f\|_{L^2(\nu_\rho)}=1}
\left[\int_\Omega \nu_\rho(d\eta)\,\eta(0)\,f(\eta)^2
- \cE(f)\right]\\
&= \sup_{\|f\|_{L^2(\nu_\rho)}=1}
\left[\int_\Omega \nu_\rho(d\eta)\,\eta(0)\,f(\eta)^2
- I_d(f^2 \nu_\rho)\right].
\end{aligned}
\ee
Here note that the smoothing in (\ref{rdSEb}) can be removed under the supremum in
(\ref{varcomp4}) (recall the remark made below (\ref{rfSEa})). But the r.h.s.\
of (\ref{varcomp4}) is precisely the r.h.s.\ of (\ref{Varocc}) for $p=1$,
where we recall (\ref{rdSEb}--\ref{rfSEa}) and put $f(\eta)^2=(d\mu/d\nu_\rho)(\eta)$.


\section{Lyapunov exponents: transient simple random walk}
\label{S4}

This section is devoted to the proof of Theorem \ref{Lyahighlim}, where $d \geq 4$ and
$p(\cdot,\cdot)$ is simple random walk given by (\ref{SRW}), i.e., $\xi$ is simple 
symmetric exclusion (SSE). \emph{The proof is long and technical}, taking up half of the
present paper. In Sections \ref{S4.1}--\ref{S4.6}, we give the proof for $p=1$. In Section 
\ref{S4.7}, we indicate how to extend the proof to arbitrary $p$.

\subsection{Outline}
\label{S4.0}

In Section \ref{S4.1}, we do an appropriate scaling in $\kappa$. In Section \ref{S4.2},
we introduce an SSE+RW generator and an auxiliary exponential martingale. In Section
\ref{S4.3}, we compute upper and lower bounds for the l.h.s.\ of (\ref{limlamb}) in
terms of certain key quantities, and we complete the proof of Theorem \ref{Lyahighlim}
(for $p=1$) subject to two propositions, whose proof is given in Sections \ref{S4.5}--\ref{S4.6}.
In Section \ref{S4.4}, we list some preparatory facts that are needed as we go along. 

As before, we write $X^\kappa_s,\xi_s(x)$ instead of $X^\kappa(s),\xi(x,s)$. We abbreviate
\be{sumsdefs}
1[\kappa] = 1 + \frac{1}{2d\kappa},
\ee
and write $\{a,b\}$ to denote the unoriented bond between nearest-neighbor sites $a,b\in\Z^d$
(recall (\ref{SRW})--(\ref{expro1})). Three parameters will be important: $t$, $\kappa$ and $T$.
We will take limits in the following order:
\be{limsorders}
t\to\infty, \quad \kappa\to\infty, \quad T\to\infty.
\ee

\subsection{Scaling}
\label{S4.1}

We have $X_t^\kappa = X_{\kappa t}$, $t\geq 0$, where $X=(X_t)_{t\geq 0}$ is simple
random walk with step rate $2d$, being independent of $(\xi_t)_{t\geq 0}$. We therefore have
\be{lb1.15}
\E_{\,\nu_\rho,0}\left(\exp\left[\int_0^t ds\,\,
\xi_s\left(X_s^\kappa\right)\right]\right)
= \E_{\,\nu_\rho,0}\left(\exp \left[\frac{1}{\kappa}\,\int_0^{\kappa t} ds\,\,
\xi_{\frac{s}{\kappa}}\left(X_s\right)\right]\right).
\ee
Define the scaled Lyapunov exponent (recall (\ref{fey-kac2}--\ref{lyap2}))
\be{lscal}
\lambda_1^*(\kappa) = \lim_{t\to\infty} \Lambda_1^*(\kappa;t)
\quad \mbox{with} \quad
\Lambda_1^*(\kappa;t) = \frac{1}{t} \log
\E_{\,\nu_\rho,0}\left(\exp\left[\frac{1}{\kappa}\int_0^t ds\,\,
\xi_{\frac{s}{\kappa}}\left(X_s\right)\right]\right).
\ee
Then $\lambda_1(\kappa)=\kappa\lambda_1^*(\kappa)$. Therefore, in what follows
we will focus on the quantity
\be{lb2.3}
\lambda_1^*(\kappa) - \frac{\rho}{\kappa}
= \lim_{t\to\infty}\frac{1}{t} \log\E_{\,\nu_\rho,0}
\left(\exp\left[\frac{1}{\kappa}\int_0^t ds\,
\left(\xi_{\frac{s}{\kappa}}(X_s)-\rho\right)\right]\right)
\ee
and compute its asymptotic behavior for large $\kappa$. We must show that
\be{limscalfin}
\lim_{\kappa\to\infty} 2d\kappa^2\left[\lambda_1^*(\kappa)-\frac{\rho}{\kappa}\right]
= \rho(1-\rho)G_d.
\ee

\subsection{SSE+RW generator and an auxiliary exponential martingale}
\label{S4.2}

For $t\geq 0$, let
\be{lb2.5}
Z_t = (\xi_{\frac{t}{\kappa}},X_t)
\ee
and denote by $\P_{\eta,x}$ the law of $Z$ starting from $Z_0=(\eta,x)$.
Then $Z=(Z_t)_{t\geq 0}$ is a Markov process on $\Omega\times\Z^d$ with
generator
\be{lb2.7}
\cA = \frac{1}{\kappa} L + \Delta
\ee
(acting on the Banach space of bounded continuous functions on $\Omega\times \Z^d$,
equipped with the supremum norm). Let $(\cP_t)_{t\geq 0}$ be the semigroup generated
by $\cA$. The following lemma will be crucial to rewrite the expectation in the r.h.s.\
of (\ref{lb2.3}) in a more manageable form.

\bl{martlem}
Fix $\kappa>0$ and $r>0$. For all $t\geq 0$ and all bounded continuous functions
$\psi\colon\,\Omega\times\Z^d\to\R$ such that $\psi$ and $\exp{[(r/\kappa)\psi]}$ belong
to the domain of $\cA$, define
\be{martdef}
\begin{aligned}
M_t^r &= \frac{r}{\kappa}\left[\psi(Z_t)-\psi(Z_0)-
\int_0^t ds\, \cA \psi(Z_s)\right],\\
N_t^r &=\exp\left[M_t^r-\int_0^t ds\,
\left[\left(e^{-\frac{r}{\kappa}\psi}\cA e^{\frac{r}{\kappa}\psi}\right)
- \cA \left(\frac{r}{\kappa}\psi\right)\right](Z_s)\right].
\end{aligned}
\ee
Then:\\
(i) $M^r=(M_t^r)_{t\geq 0}$ is a $\P_{\eta,x}$-martingale for all $(\eta,x)$.\\
(ii) For $t\geq 0$, let $\cP_t^{\new}$ be the operator defined by
\be{lb2i.1}
(\cP_t^{\new}f)(\eta,x)
= e^{-\frac{r}{\kappa}\psi(\eta,x)}\,\E_{\,\eta,x}\left(\exp\left[-\int_0^t ds\,
\left(e^{-\frac{r}{\kappa}\psi}\cA e^{\frac{r}{\kappa}\psi}\right)(Z_s)
\right]\left(e^{\frac{r}{\kappa}\psi}f\right)(Z_t)\right)
\ee
for bounded continuous $f\colon\,\Omega\times\Z^d\to\R$. Then $(\cP_t^{\new})_{t\geq 0}$
is a strongly continuous semigroup with generator
\be{lb2i.2}
(\cA^{\new}f)(\eta,x)
= \left[e^{-\frac{r}{\kappa}\psi}\cA\left(e^{\frac{r}{\kappa}\psi}f
\right) - \left(e^{-\frac{r}{\kappa}\psi}\cA e^{\frac{r}{\kappa}\psi}
\right)f\right](\eta,x).
\ee
(iii) $N^r=(N_t^r)_{t\geq 0}$ is a $\P_{\eta,x}$-martingale for all $(\eta,x)$.\\
(iv) Define a new path measure $\P_{\eta,x}^{\new}$ by putting
\be{newpathdef}
\frac{d\,\P_{\eta,x}^{\new}}{d\,\P_{\eta,x}}\left((Z_s)_{0\leq s \leq t}\right)
= N_t^r, \qquad t \geq 0.
\ee
Then, under $\P_{\eta,x}^{\new}$, $(Z_t)_{t\geq 0}$ is a Markov process with
semigroup $(\cP_t^{\new})_{t\geq 0}$.
\el

\bpr
The proof is standard.

\noindent
(i) This follows from the fact that $\cA$ is a Markov generator and $\psi$
belongs to its domain (see Liggett \cite{lig85}, Chapter I, Section 5).

\noindent
(ii) Let $\eta\in\Omega$, $x\in\Z^d$ and $f\colon\,\Omega\times\Z^d\to\R$ bounded
measurable. Rewrite (\ref{lb2i.1}) as
\be{lb2i.7}
\begin{aligned}
(\cP_{t}^{\new}f)(\eta,x)
&= \E_{\,\eta,x}\left(\exp\left[\frac{r}{\kappa}\psi(Z_t)
-\frac{r}{\kappa}\psi(Z_0)-\int_0^t ds\,
\left(e^{-\frac{r}{\kappa}\psi}\cA e^{\frac{r}{\kappa}\psi}\right)(Z_s)\right]
f(Z_t)\right)\\
&= \E_{\,\eta,x}\left(N_t^r f(Z_t)\right).
\end{aligned}
\ee
This gives
\be{lb2i.5}
(\cP_0^{\new}f)(\eta,x) = f(\eta,x)
\ee
and
\be{lb2i.9}
\begin{aligned}
(\cP_{t_1+t_2}^{\new}f)(\eta,x)
&=\E_{\,\eta,x}\left(N_{t_1+t_2}^r f(Z_{t_1+t_2})\right)
=\E_{\,\eta,x}\left(N_{t_1}^r \frac{N_{t_1+t_2}^r}{N_{t_1}^r}f(Z_{t_1+t_2})\right)\\
&=\E_{\,\eta,x}\bigg(N_{t_1}^r\E_{Z_{t_1}}\Big(N_{t_2}^r f(Z_{t_2})\Big)\bigg)
= \Big(\cP_{t_1}^{\new}(\cP_{t_2}^{\new}f)\Big)(\eta,x),
\end{aligned}
\ee
where we use the Markov property of $Z$ at time $t_1$ (under $\P_{\eta,x}$) together
with the fact that $N_{t_1+t_2}^r/N_{t_1}^r$ only depends on $Z_t$ for $t\in [t_1,t_1+t_2]$.
Equations (\ref{lb2i.5}--\ref{lb2i.9}) show that $(\cP_t^\new)_{t\geq 0}$ is a semigroup
which is easily seen to be strongly continuous.

Taking the derivative of (\ref{lb2i.1}) in the norm w.r.t.\ $t$ at $t=0$, we get (\ref{lb2i.2}).
Next, if $f\equiv 1$, then (\ref{lb2i.2}) gives $\cA^\new 1 = 0$. This last equality implies that
\be{lb2i.17}
\frac{1}{\lambda}\left(\lambda\Id - \cA^\new\right)1 = 1
\quad\forall\,\lambda>0.
\ee
Since $\lambda\Id - \cA^\new$ is invertible,
we get
\be{lb2i.19}
\left(\lambda\Id - \cA^\new\right)^{-1} 1 = \frac{1}{\lambda}
\quad\forall\,\lambda>0,
\ee
i.e.,
\be{lb2i.21}
\int_0^\infty dt\, e^{-\lambda t}\,\cP_t^\new 1 = \frac{1}{\lambda}
\quad\forall\,\lambda>0.
\ee
Inverting this Laplace transform, we see that
\be{lb2i.22}
\cP_t^\new 1 = 1 \quad\forall\,t\geq 0.
\ee

\noindent
(iii) Fix $t\geq 0$ and $h>0$. Since $N_t^r$ is $\cF_t$-measurable,
with $\cF_t$ the $\sigma$-algebra generated by $(Z_s)_{0\leq s\leq t}$, we have
\be{lb2i.23}
\begin{aligned}
&\E_{\,\eta,x}\left(N_{t+h}^r\big\vert\cF_t\right)\\
&= N_t^r\, \E_{\,\eta,x}\left(\exp\left[M_{t+h}^r- M_t^r -\int_t^{t+h}ds\,
\left[\left(e^{-\frac{r}{\kappa}\psi}\cA e^{\frac{r}{\kappa} \psi}\right)
- \cA \left(\frac{r}{\kappa}\psi\right)\right](Z_s)\right]
\bigg\vert\,\cF_t\right).
\end{aligned}
\ee
Applying the Markov property of $Z$ at time $t$, we get
\be{lb2i.25}
\begin{aligned}
\E_{\,\eta,x}\left(N_{t+h}^r\mid\cF_t\right)
&= N_t^r\, \E_{Z_t}\left(\exp\left[\frac{r}{\kappa}\psi(Z_h)
- \frac{r}{\kappa}\psi(Z_0) -\int_0^{h}ds\,
\left(e^{-\frac{r}{\kappa}\psi}\cA
e^{\frac{r}{\kappa}\psi}\right)(Z_s)\right]\right)\\
&= N_t^r\left(\cP_h^\new 1\right)(Z_t)
= N_t^r,
\end{aligned}
\ee
where the third equality uses (\ref{lb2i.22}).

\noindent
(iv) This follows from (iii) via a calculation similar to (\ref{lb2i.9}).
\epr

\subsection{Proof of Theorem \ref{Lyahighlim}}
\label{S4.3}

In this section we compute upper and lower bounds for the r.h.s.\ of (\ref{lb2.3})
in terms of certain key quantities (Proposition \ref{lbholder} below). We then state
two propositions for these quantities (Propositions \ref{mainLem1}--\ref{remainLem1}
below), from which Theorem \ref{Lyahighlim} will follow. The proof of these two
propositions is given in Sections \ref{S4.5}--\ref{S4.6}.

For $T>0$, let $\psi\colon\,\Omega\times\Z^d$
be defined by
\be{Gap1.57}
\psi(\eta,x)
=\int_{0}^{T}ds \left(\cP_s \phi\right)(\eta,x)
\quad \mbox{with} \quad \phi(\eta,x)=\eta(x)-\rho,
\ee
where $(\cP_t)_{t\geq 0}$ is the semigroup generated by $\cA$ (recall (\ref{lb2.7})).
We have
\be{Gap1.65}
\psi(\eta,x)
=\int_0^T ds\, \E_{\,\eta,x}\left(\phi(Z_s)\right)
=\int_0^T ds\, \E_{\,\eta}\, \sum_{y\in\Z^d}
p_{2ds}(y,x) \Big(\xi_{\frac{s}{\kappa}}(y)-\rho\Big),
\ee
where $p_t(x,y)$ is the probability that simple random walk with step rate 1 moves
from $x$ to $y$ in time $t$ (recall that we assume (\ref{SRW})). Using (\ref{graph}),
we obtain the representation
\be{Gap1.69}
\psi(\eta,x) = \int_0^{T} ds\, \sum_{z\in\Z^d}
p_{2ds\onek}(z,x)\big[\eta(z)-\rho\big],
\ee
where $\onek$ is given by (\ref{sumsdefs}). Note that $\psi$ depends on $\kappa$ and $T$.
We suppress this dependence. Similarly,
\be{Gap1.71}
-\cA\psi = \int_0^T ds\,(-\cA\cP_s \phi) = \phi- \cP_T \phi,
\ee
with
\be{Gap1.73}
(\cP_T \phi)(\eta,x)
= \E_{\,\eta,x}\left(\phi(Z_T)\right)
= \E_{\,\eta,x}\left(\xi_{\frac{T}{\kappa}}(X_{T})-\rho\right)
= \sum_{z\in\Z^d} p_{2dT\onek}(z,x)\,[\eta(z)-\rho].
\ee

The auxiliary function $\psi$ will play a key role throughout the remaining sections.
The integral in (\ref{Gap1.57}) is a regularization that is useful when dealing with
central limit type behavior of Markov processes (see e.g.\ Kipnis \cite{kip87}).
Heuristically, $T=\infty$ corresponds to $-\cA\psi=\phi$. Later we will let
$T\to\infty$.

The following proposition serves as the starting point of our asymptotic analysis.

\bp{lbholder}
For any $\kappa,T>0$,
\be{Gap1.27}
\lambda_1^\ast(\kappa)-\frac{\rho}{\kappa}
\,\,{\leq \atop \geq}\,\,I_1^{r,q}(\kappa,T) + I_2^{r,q}(\kappa,T),
\ee
where
\be{I12defs}
\begin{aligned}
I_1^{r,q}(\kappa,T) &=
\frac1{2q}\limsup_{t\to\infty}\frac{1}{t}\log\E_{\,\nu_\rho,0}\bigg(\exp\bigg[
\frac{2q}r\int_0^t ds\, \left[
\left(e^{-\frac{r}{\kappa}\psi} \cA e^{\frac{r}{\kappa}\psi}\right)
-\cA\left(\frac{r}{\kappa}\psi\right)\right](Z_s)\bigg]\bigg),\\
I_2^{r,q}(\kappa,T) &=
\frac1{2q}\limsup_{t\to\infty}\frac{1}{t}\log\E_{\,\nu_\rho,0}\bigg(\exp\bigg[
\frac{2q}{\kappa}\int_0^t ds\, \left(\cP_T\phi\right)(Z_s) \bigg]\bigg),
\end{aligned}
\ee
and $1/r + 1/q=1$ for any $r,q>1$ in the first inequality and any $q<0<r<1$
in the second inequality.
\ep

\bpr
Recall (\ref{lb2.3}). From the first line of (\ref{martdef}) and (\ref{Gap1.71})
it follows that
\be{Gap1.43}
\frac1r M_t^r+\frac{1}{\kappa}\psi(Z_0)-\frac{1}{\kappa}\psi(Z_t)
= \frac{1}{\kappa}\int_0^t ds\,[(-\cA)\psi](Z_s)
= \frac{1}{\kappa}\int_0^t ds\,\phi(Z_s)
-\frac{1}{\kappa}\int_0^t ds\, \left(\cP_T\phi\right)(Z_s).
\ee
Hence
\be{Gap1.47}
\begin{aligned}
&\E_{\,\nu_\rho,0}\left(\exp\left[\frac{1}{\kappa}\int_0^tds\,
\phi(Z_s)\right]\right)\\
&\qquad\qquad\qquad=\E_{\,\nu_\rho,0}
\left(\exp\left[\frac1r M_t^r+\frac{1}{\kappa}\psi(Z_0)
-\frac{1}{\kappa}\psi(Z_t)
+\frac{1}{\kappa}\int_0^t ds\, (\cP_T\phi)(Z_s)\right]\right)\\
&\qquad\qquad\qquad= \E_{\,\nu_\rho,0}\Big(\exp\Big[ U_t^r + \frac1r V_t^r\Big]\Big)
\end{aligned}
\ee
with
\be{Gap1.33}
U_t^r =\frac1r\int_0^t ds\, \left[
\left(e^{-\frac{r}{\kappa}\psi} \cA e^{\frac{r}{\kappa}\psi}\right)
-\cA\left(\frac{r}{\kappa}\psi\right)\right](Z_s)
+\frac{1}{\kappa}\Big(\psi(Z_0)-\psi(Z_t)\Big)
+\frac{1}{\kappa}\int_0^t ds\, \left(\cP_T\phi\right)(Z_s)
\ee
and
\be{Gap1.57*}
V_t^r =M_t^r -\int_0^t ds\, \left[
\left(e^{-\frac{r}{\kappa}\psi} \cA e^{\frac{r}{\kappa}\psi}\right)
-\cA\left(\frac{r}{\kappa}\psi\right)\right](Z_s).
\ee
By H\"older's inequality, with $r,q>1$ such that $1/r+1/q=1$,
it follows from (\ref{Gap1.47}) that
\be{Gap1.15}
\begin{aligned}
\E_{\,\nu_\rho}\left(\exp\left[\frac{1}{\kappa}\int_0^t ds\,
\phi(Z_s)\right]\right)
&\leq \Big(\E_{\,\nu_\rho,0}\big(\exp\big[V_t^r\big]\big)\Big)^{1/r}
\Big(\E_{\,\nu_\rho,0}\big(\exp\big[qU_t^r\big]\big)\Big)^{1/q}\\
&=\Big(\E_{\,\nu_\rho,0}\big(\exp\big[qU_t^r\big]\big)\Big)^{1/q},
\end{aligned}
\ee
where the second line of (\ref{Gap1.15}) comes from the fact that $N_t^r=\exp[V_t^r]$
is a martingale, by Lemma \ref{martlem}(iii). Similarly, by the reverse of H\"older's
inequality, with $q<0<r<1$ such that $1/r+1/q=1$, it follows from (\ref{Gap1.47}) that
\be{Gap1.16}
\begin{aligned}
\E_{\,\nu_\rho}\left(\exp\left[\frac{1}{\kappa}\int_0^t ds\,
\phi(Z_s)\right]\right)
&\geq \Big(\E_{\,\nu_\rho,0}\big(\exp\big[V_t^{r}\big]\big)\Big)^{1/r}
\Big(\E_{\,\nu_\rho,0}\big(\exp\big[qU_t^{r}\big]\big)\Big)^{1/q}\\
&=\Big(\E_{\,\nu_\rho,0}\big(\exp\big[qU_t^{r}\big]\big)\Big)^{1/q}.
\end{aligned}
\ee
The middle term in the r.h.s.\ of (\ref{Gap1.33}) can be discarded,
because (\ref{Gap1.69}) shows that $-\rho T\leq \psi\leq (1-\rho)T$. Apply the Cauchy-Schwarz
inequality to the r.h.s.\ of (\ref{Gap1.15}--\ref{Gap1.16}) to separate the other two terms
in the r.h.s.\ of (\ref{Gap1.33}).
\epr

Note that in the r.h.s.\ of (\ref{I12defs}) the prefactors of the logarithms and the prefactors
in the exponents are \emph{both positive for the upper bound and both negative for the lower
bound}. This will be important later on.

The following two propositions will be proved in Sections \ref{S4.5}-\ref{S4.6}, respectively.
Abbreviate
\be{limorders}
\limsup_{t,\kappa,T\to\infty}
=\limsup_{T\to\infty}\limsup_{\kappa\to\infty}\limsup_{t\to\infty}.
\ee

\bp{mainLem1}
If $d\geq 3$, then for any $\alpha\in\R$ and $r>0$,
\be{main1-3}
\limsup_{t,\kappa,T\to\infty}\frac{\kappa^2}{t}
\log\E_{\,\nu_\rho,0}\left(\exp\left[\frac\alpha{r}\int_0^t ds\,
\left[e^{-\frac{r}{\kappa} \psi} \cA e^{\frac{r}{\kappa}\psi}
-\cA\Big(\frac{r}{\kappa}\psi\Big)\right](Z_s)\right]\right)
\leq\alpha r\,\rho(1-\rho)\frac{1}{2d}G_d.
\ee
\ep

\bp{remainLem1}
If $d\geq 4$, then for any $\alpha\in\R$,
\be{remain1-3}
\limsup_{t,\kappa,T\to\infty}
\frac{\kappa^2}{t}\log \E_{\,\nu_\rho,0}\bigg(\exp\bigg[
\frac{\alpha}{\kappa}\int_0^t ds\,\left(\cP_T\phi\right)(Z_s)\bigg]\bigg)\leq 0.
\ee
\ep

Picking $\alpha=2q$ in Proposition \ref{mainLem1}, we see that the first term in
the r.h.s.\ of (\ref{Gap1.27}) satisfies the bounds
\be{I1bds}
\begin{aligned}
\limsup_{T\to\infty}\limsup_{\kappa\to\infty} \kappa^2 I_1^{r,q}(\kappa,T)
&\leq r\,\rho(1-\rho)\frac{1}{2d}G_d &\mbox{if } r>1,\\
\liminf_{T\to\infty}\liminf_{\kappa\to\infty} \kappa^2 I_1^{r,q}(\kappa,T)
&\geq r\,\rho(1-\rho)\frac{1}{2d}G_d &\mbox{if } r<1.
\end{aligned}
\ee
Letting $r$ tend to 1, we obtain
\be{I1bds*}
\lim_{T\to\infty}\lim_{\kappa\to\infty} \kappa^2 I_1^{r,q}(\kappa,T)
= \rho(1-\rho)\frac{1}{2d}G_d.
\ee
Picking $\alpha=2q$ in Proposition \ref{remainLem1}, we see that the second
term in the r.h.s.\ of (\ref{Gap1.27}) satisfies
\be{I2bds}
\limsup_{T\to\infty}\limsup_{\kappa\to\infty} \kappa^2 I_2^{r,q}(\kappa,T) = 0
\qquad \mbox{if } d \geq 4.
\ee
Combining (\ref{I1bds*}--\ref{I2bds}), we see that we have completed the proof
of Theorem \ref{Lyahighlim} for $d\geq 4$.

In order to prove Conjecture \ref{cjd=3}, we would have to extend Proposition
\ref{remainLem1} to $d=3$ and show that it contributes the second term in the
r.h.s.\ of (\ref{remain1-3}) rather than being negligible.

\subsection{Preparatory facts and notation}
\label{S4.4}

In order to estimate $I_1^{r,q}(\kappa,T)$ and $I_2^{r,q}(\kappa,T)$, we need a number
of preparatory facts. These are listed in Lemmas \ref{fdifbdlem}--\ref{approxLaplem}
below.

It follows from (\ref{Gap1.69}) that
\be{psidifext}
\psi(\eta,b)-\psi(\eta,a) = \int_0^T ds \sum_{z\in\Z^d}\big(p_{2ds\onek}(z,b)-p_{2ds\onek}(z,a)\big)
[\eta(z)-\rho]
\ee
and
\be{Tayl.3}
\begin{aligned}
\psi\big(\eta^{a,b},x\big)-\psi(\eta,x)
&= \int_0^T ds\,
\sum_{z\in\Z^d}
p_{2ds\onek}(z,x)\left[\eta^{a,b}(z)-\eta(z)\right]\\
&= \int_0^T ds\,
\left(p_{2ds\onek}(b,x)-p_{2ds\onek}(a,x)\right)
\left[\eta(a)-\eta(b)\right],
\end{aligned}
\ee
where we recall the definitions of $1[\kappa]$ and $\eta^{a,b}$ in (\ref{sumsdefs}) and
(\ref{expro2}), respectively. We need bounds on both these differences.

\bl{fdifbdlem}
For any $\eta\in\Omega$, $a,b,x\in\Z^d$ and $\kappa,T>0$,
\be{fdifest0}
\left|\psi\big(\eta,b\big)-\psi(\eta,a)\right| \leq 2T,
\ee
\be{fdifest1}
\left|\psi\big(\eta^{a,b},x\big)-\psi(\eta,x)\right|
\leq 2G_d < \infty,
\ee
and
\be{fdifest2}
\sum_{\{a,b\}}\Big(\psi\big(\eta^{a,b},x\big)-\psi(\eta,x)\Big)^2
\leq \frac1{2d}G_d < \infty,
\ee
where $G_d$ is the Green function at the origin of simple random walk.
\el

\bpr
The bound in (\ref{fdifest0}) is immediate from (\ref{psidifext}). By (\ref{Tayl.3}),
we have
\be{Tayl.7}
\left|\psi\big(\eta^{a,b},x\big)-\psi(\eta,x)\right|
\leq \int_0^T ds\,
\left|p_{2ds\onek}(b,x)
-p_{2ds\onek}(a,x)\right|.
\ee
Using the bound $p_t(x,y)\leq p_t(0,0)$ (which is immediate from the Fourier representation
of the transition kernel), we get
\be{Tayl.21}
\left|\psi\big(\eta^{a,b},x\big)-\psi(\eta,x)\right|
\leq 2\int_0^\infty ds\,p_{2ds\onek}(0,0) \leq 2G_d.
\ee
Again by (\ref{Tayl.3}), we have
\be{Tayl.23}
\begin{aligned}
&\sum_{\{a,b\}}\Big(\psi\big(\eta^{a,b},x\big)-\psi(\eta,x)\Big)^2
= \sum_{\{a,b\}} [\eta[(a)-\eta(b)]^2
\left(\int_0^T ds\,\Big(p_{2ds\onek}(b,x)-p_{2ds\onek}(a,x)\Big)\right)^2\\
&\qquad\leq 2\int_0^T du\, \int_u^T dv\, \sum_{\{a,b\}}
\Big(p_{2du\onek}(b,x)-p_{2du\onek}(a,x)\Big)
\Big(p_{2dv\onek}(b,x)-p_{2dv\onek}(a,x)\Big)\\
&\qquad=-2\int_0^T du\, \int_u^T dv\, \sum_{a\in\Z^d}p_{2du\onek}(a,x)
\Big[\Delta_1 p_{2dv\onek}(a,x)\Big]\\
&\qquad = -\frac{2}{\onek}\int_0^T du\, \int_u^T dv\,
\sum_{a\in\Z^d} p_{2du\onek}(a,x)
\bigg[\frac{\partial}{\partial v} p_{2dv\onek}(a,x)\bigg]\\
&\qquad = -\frac{2}{\onek}\int_0^T du\,
\sum_{a\in\Z^d} p_{2du\onek}(a,x) \big(p_{2dT\onek}(a,x)-p_{2du\onek}(a,x)\big)\\
&\qquad \leq\frac{2}{\onek}\int_0^T du\,
\sum_{a\in\Z^d}p_{2du\onek}^2(a,x)\\
&\qquad \leq\frac{2}{\onek}\int_0^\infty du\,
p_{4du\onek}(0,0)
=\frac{1}{2d(\onek)^2}G_d(0)
\leq \frac1{2d} G_d,
\end{aligned}
\ee
where $\Delta_1$ denotes the discrete Laplacian acting on the first coordinate, and
in the fifth line we use that $(\partial/\partial t)p_t =(1/2d)\Delta_1 p_t$.
\epr

For $x\in\Z^d$, let $\tau_x\colon \Omega\to\Omega$ be the $x$-shift on $\Omega$
defined by
\be{shift}
\big(\tau_x \eta\big)(z)=\eta(z+x), \qquad \eta\in\Omega,\,z\in\Z^d.
\ee

\bl{spec-rep}
For any bounded measurable $W\colon \Omega\times \Z^d \to \R$,
\be{spec-rep-3}
\begin{aligned}
&\limsup_{t\to\infty}\frac1t \log \E_{\,\nu_\rho,0}\bigg(\exp\bigg[
\int_0^t ds\,\, W\big(\xi_{\frac{s}{\kappa}},X_s\big)\bigg]\bigg)\\
&\qquad\leq \limsup_{t\to\infty}\frac1t \log \E_{\,\nu_\rho}\bigg(\exp\bigg[
\int_0^t ds\,\, W\big(\xi_{\frac{s}{\kappa}},0\big)\bigg]\bigg),
\end{aligned}
\ee
provided
\be{spec-rep-2}
W(\eta,x)=W(\tau_x\eta,0)\quad\forall\, \eta\in\Omega,\,\, x\in\Z^d.
\ee
\el

\bpr
The proof uses arguments similar to those in Sections \ref{S2.2} and \ref{S3.2.6}.
Recall (\ref{lb2.5}). Proposition \ref{RR} with $p=1$, applied to the self-adjoint
operator $G_W^\kappa=\frac1\kappa L +\Delta+W$ (instead of $G_V^\kappa$ in
(\ref{gen}--\ref{combgen})), gives
\be{er1-17}
\lim_{t\to\infty}\frac1t \log \E_{\,\nu_\rho,0}\left(\exp\bigg[
\int_0^t ds\, W(Z_s)\bigg]\right)
=\sup_{\|f\|_{L^2(\nu_\rho\otimes m)}=1}
\bigg(B_1(f)-\frac1\kappa B_2(f)-B_3(f)\bigg)
\ee
with
\be{er1-19}
\begin{aligned}
B_1(f) &= \int_\Omega \nu_\rho(d\eta)\sum_{z\in\Z^d} W(\eta,z)\, f(\eta,z)^2,\\[0.3cm]
B_2(f) &= \int_\Omega \nu_\rho(d\eta)\sum_{z\in\Z^d}\frac12 \sum_{\{x,y\}\subset\Z^d}
p(x,y)[f(\eta^{x,y},z)-f(\eta,z)]^2,\\[0.3cm]
B_3(f) &= \int_\Omega \nu_\rho(d\eta)\sum_{z\in\Z^d}\frac12 \sum_{{y\in\Z^d} \atop {\|y-z\|=1}}
[f(\eta,y)-f(\eta,z)]^2.
\end{aligned}
\ee
An upper bound is obtained by dropping $B_3(f)$, i.e., the part associated with the simple
random walk $X$. After that, split the supremum into two parts,
\be{er1-21}
\begin{aligned}
&\sup_{\|f\|_{L^2(\nu_\rho\otimes m)}=1}
\Big(B_1(f)- B_2(f)\Big)\\
&=\sup_{\|g\|_{L^2(m)}=1}\,\,
\sup_{\|f_z\|_{L^2(\nu_\rho)}=1\,\forall\,z\in\Z^d}
\sum_{z\in\Z^d} g(z)^2 \int_\Omega \nu_\rho(d\eta)\\
&\qquad\qquad\qquad\times
\bigg(W(\eta,z)\, f_z(\eta)^2-\frac12\sum_{\{x,y\}\subset\Z^d}
p(x,y)[f_z(\eta^{x,y})-f_z(\eta)]^2\bigg),
\end{aligned}
\ee
where $f_z(\eta)=f(\eta,z)/g(z)$ with $g(z)^2=\int_\Omega \nu_\rho(d\eta) f(\eta,z)^2$.
The second supremum in (\ref{er1-21}), which runs over a family of functions indexed by $z$,
can be brought under the sum. This gives
\be{er1-25}
\begin{aligned}
&{\rm r.h.s.\ (\ref{er1-21})} =
\sup_{\|g\|_{L^2(m)}=1} \sum_{z\in\Z^d} g(z)^2
\sup_{\|f_z\|_{L^2(\nu_\rho)}=1}
\int_\Omega \nu_\rho(d\eta)\\
&\qquad\qquad\qquad\qquad\qquad\times
\bigg(W(\eta,z)\, f_z(\eta)^2-\frac12\sum_{\{x,y\}\subset\Z^d}
p(x,y)[f_z(\eta^{x,y})-f_z(\eta)]^2\bigg).
\end{aligned}
\ee
By (\ref{spec-rep-2}) and the shift-invariance of $\nu_\rho$, we may replace $z$ by $0$ under
the second supremum in (\ref{er1-25}), in which case the latter no longer depends on $z$,
and we get
\be{er1-29}
\begin{aligned}
{\rm r.h.s.\ (\ref{er1-25})}
&=\sup_{\|f\|_{L^2(\nu_\rho)}=1}
\int_\Omega \nu_\rho(d\eta) \bigg[W(\eta,0)\, f(\eta)^2
- \frac12\sum_{\{x,y\}\subset\Z^d}
p(x,y)[f(\eta^{x,y})-f(\eta)]^2\bigg]\\
&=\lim_{t\to\infty}\frac1t \log \E_{\,\nu_\rho}\left(\exp\bigg[
\int_0^t ds\,\, W\Big(\xi_{\frac{s}{\kappa}},0\Big)
\bigg]\right),
\end{aligned}
\ee
where the second equality comes from the analogue of Proposition \ref{RR} with self-adjoint
operator $\frac1\kappa L+W(\cdot,0)$ (instead of $G_V^\kappa$).
\epr

\bl{varform}
For any $\rho\in (0,1)$,
\be{varform3}
\max_{\beta\in[0,1]} \Big[\gamma \beta
-\frac{1}{2G_d}\big(\sqrt{\beta}-\sqrt{\rho}\big)^2\Big]
= \frac{\rho\gamma}{1-2G_d\gamma}\sim\rho\gamma
\quad\text{as}\quad \gamma \downarrow 0.
\ee
\el

\bpr
A straightforward computation shows that the maximum in (\ref{varform3}) is attained at
\be{varform7}
\beta=\frac{\rho}{\big(1-2G_d\gamma\big)^2}\in[0,1]
\ee
for small enough $\gamma$. Substitution yields the claim.
\epr

\bl{approxLaplem}
There exists $C>0$ such that, for all $t\geq 0$ and $x,y\in\Z^d$,
\be{approxLap1}
p_t(x,y) \leq \frac{C}{(1+t)^{\frac{d}{2}}}.
\ee
\el

\bpr
This is a standard fact. Indeed, we can decompose the transition kernel of simple random walk
with step rate $1$ as
\be{pprod}
p_{dt}(x,y) = \prod_{j=1}^d p_t^{(1)}(x^j,y^j),
\qquad x=(x^1,\dots,x^d),\,y=(y^1,\dots,y^d),
\ee
where $p_t^{(1)}(x,y)$ is the transition kernel of 1-dimensional simple random walk with
step rate $1$. In Fourier representation,
\be{approxLap1.15}
\begin{aligned}
p_t^{(1)}(x,y) &= \frac{1}{2\pi}\int_{-\pi}^{\pi} dk\,
e^{ik\cdot (y-x)}\, e^{-t\hat\varphi(k)},\,\,\,
\hat\varphi(k)=1-\cos k.
\end{aligned}
\ee
The bound in (\ref{approxLap1}) follows from (\ref{pprod}) and
\be{approxLap1.17}
p_t^{(1)}(x,y) \leq p_t^{(1)}(0,0)
= \frac{1}{2\pi}\int_{-\pi}^{\pi} dk\,e^{-t\hat\varphi(k)}
\leq \frac{C}{(1+t)^{\frac{1}{2}}},
\quad t \geq 0,\,x,y\in\Z^d.
\ee
\epr

\subsection{Proof of Proposition \ref{mainLem1}}
\label{S4.5}

The proof of Proposition \ref{mainLem1} is given in Section \ref{main1} subject to four
lemmas. The latter will be proved in Sections \ref{main2}--\ref{main5}, respectively.
All results are valid for $d\geq 3$.

\subsubsection{Proof of Proposition \ref{mainLem1}}
\label{main1}

\bl{mainLem2}
Uniformly in $\eta\in\Omega$ and $x\in\Z^d$, as $\kappa\to\infty$,
\be{main2-3}
\left[\left(e^{-\frac{r}{\kappa} \psi} \cA\,
e^{\frac{r}{\kappa}\psi}\right)
-\cA\,\left(\frac{r}{\kappa}\psi\right)\right](\eta,x)
= \frac{r^2}{2\kappa^2}\sum_{e\colon\, \|e\|=1}
\Big(\psi(\eta,x+e)-\psi(\eta,x)\Big)^2
+ O\bigg(\Big(\frac{1}{\kappa}\Big)^3\bigg).
\ee
\el

\bl{mainLem3}
For any $\kappa,T>0$, $\alpha\in\R$ and $r>0$,
\be{main3-3}
\begin{aligned}
&\limsup_{t\to\infty}\frac{1}{t}\log \E_{\,\nu_\rho,0}
\Bigg(\exp\Bigg[\frac{\alpha r}{2\kappa^2}\int_0^t ds\, \sum_{e\colon\,\|e\|=1}
\bigg(\psi\Big(\xi_{\frac{s}\kappa},X_s+e\Big)
-\psi\Big(\xi_{\frac{s}\kappa},X_s\Big)\bigg)^2\Bigg]\Bigg)\\
&\leq \limsup_{t\to\infty}\frac{1}{2t}\log \E_{\,\nu_\rho}
\Bigg(\exp\Bigg[\frac{\alpha r}{\kappa^2}\int_0^t ds\,
\sum_{z\in\Z^d} K_\diag^{\kappa,T}(z)\Big(\xi_{\frac{s}\kappa}(z)-\rho\Big)^2\Bigg]\Bigg)\\
&\quad+ \limsup_{t\to\infty}\frac{1}{2t}\log \E_{\,\nu_\rho}
\Bigg(\exp\Bigg[\frac{\alpha r}{\kappa^2}\int_0^t ds\,
\sum_{{z_1,z_2\in\Z^d}\atop{z_1\neq z_2}} K_\off^{\kappa,T}(z_1,z_2)
\Big(\xi_{\frac{s}\kappa}(z_1)-\rho\Big)\Big(\xi_{\frac{s}\kappa}(z_2)-\rho\Big)\Bigg]\Bigg),
\end{aligned}
\ee
where
\be{main3-5}
\begin{aligned}
K_\diag^{\kappa,T}(z)&=\sum_{e\colon\,\|e\|=1}\Big(\chi(z+e)-\chi(z)\Big)^2,\\
K_\off^{\kappa,T}(z_1,z_2)&=\sum_{e:\|e\|=1}
\Big(\chi(z_1+e)-\chi(z_1)\Big)\Big(\chi(z_2+e)-\chi(z_2)\Big),\quad z_1\neq z_2,
\end{aligned}
\ee
with
\be{main3-7}
\chi(z)=\int_0^T du\, p_{2du\onek}(0,z).
\ee
\el

\bl{mainLem4}
For any $\alpha\in\R$ and $r>0$,
\be{main4-3}
\limsup_{t,\kappa,T\to\infty}\frac{\kappa^2}{t}\log \E_{\,\nu_\rho}
\Bigg(\exp\Bigg[\frac{\alpha r}{\kappa^2}\int_0^t ds\,
\sum_{z\in\Z^d} K_\diag^{\kappa,T}(z)\Big(\xi_{\frac{s}\kappa}(z)-\rho\Big)^2\Bigg]\Bigg)
\leq \alpha r\,\rho(1-\rho)\frac{1}{d}G_d.
\ee
\el

\bl{mainLem5}
For any $\alpha\in\R$ and $r>0$,
\be{main5-3}
\limsup_{t,\kappa,T\to\infty}\frac{\kappa^2}{t}\log \E_{\,\nu_\rho}
\Bigg(\exp\Bigg[\frac{\alpha r}{\kappa^2}\int_0^t ds\,
\sum_{{z_1,z_2\in\Z^d}\atop{z_1\neq z_2}} K_\off^{\kappa,T}(z_1,z_2)
\Big(\xi_{\frac{s}\kappa}(z_1)-\rho\Big)\Big(\xi_{\frac{s}\kappa}(z_2)-\rho\Big)\Bigg]\Bigg)
\leq0.
\ee
\el

\noindent
Combining Lemmas \ref{mainLem2}--\ref{mainLem5}, we obtain the claim in Proposition
\ref{mainLem1}.

\subsubsection{Proof of Lemma \ref{mainLem2}}
\label{main2}

Lemma \ref{mainLem2} is immediate from (\ref{lb2.7}) and the following two lemmas.

\bl{mainLem2.1}
Uniformly in $\eta\in\Omega$ and $x\in\Z^d$, as $\kappa\to\infty$,
\be{main2.1-3}
\frac{1}{\kappa}\left[\left(e^{-\frac{r}{\kappa} \psi} L\,
e^{\frac{r}{\kappa}\psi}\right)
-L\,\left(\frac{r}{\kappa}\psi\right)\right](\eta,x)
= O\left(\frac1{\kappa^3}\right).
\ee
\el
\bl{mainLem2.2}

Uniformly in $\eta\in\Omega$ and $x\in\Z^d$, as $\kappa\to\infty$,
\be{main2.2-3}
\left[\left(e^{-\frac{r}{\kappa} \psi} \Delta\,
e^{\frac{r}{\kappa}\psi}\right)
-\Delta\,\left(\frac{r}{\kappa}\psi\right)\right](\eta,x)
= \frac{r^2}{2\kappa^2}\sum_{e\colon\,\|e\|=1}
\Big(\psi(\eta,x+e)-\psi(\eta,x)\Big)^2
+ O\left(\frac{1}{\kappa^3}\right).
\ee
\el

\bprl{mainLem2.1}
By (\ref{SRW}--\ref{expro1}), we have
\be{main2.1-5}
\begin{aligned}
&\left[\left(e^{-\frac{r}{\kappa} \psi} L\,
e^{\frac{r}{\kappa}\psi}\right)
-L\,\left(\frac{r}{\kappa}\psi\right)\right](\eta,x)\\
&\qquad\qquad =\frac{1}{2d} \sum_{\{a,b\}}
\left(e^{\frac{r}{\kappa}[\psi(\eta^{a,b},x)-\psi(\eta,x)]}
-1-\frac{r}{\kappa}[\psi(\eta^{a,b},x)-\psi(\eta,x)]\right).
\end{aligned}
\ee
Taylor expansion of the r.h.s.\ of (\ref{main2.1-5}) gives
that uniformly in $\eta\in\Omega$ and $x\in\Z^d$,
\be{main2.1-7}
\frac{1}{\kappa}\left[\left(e^{-\frac{r}{\kappa} \psi} L\,
e^{\frac{r}{\kappa}\psi}\right)
-L\,\left(\frac{r}{\kappa}\psi\right)\right](\eta,x)
= \frac{r^2}{4d\kappa^3} \sum_{\{a,b\}}
\Big(\psi\big(\eta^{a,b},x\big)-\psi(\eta,x)\Big)^2
e^{o(1)}\\
= O\left(\frac1{\kappa^3}\right),
\ee
where we use (\ref{fdifest1}--\ref{fdifest2}).
\eprl

\bprl{mainLem2.2}
By (\ref{dL}), we have
\be{main2.2-5}
\begin{aligned}
&\left[\left(e^{-\frac{r}{\kappa} \psi} \Delta\,
e^{\frac{r}{\kappa}\psi}\right)
-\Delta\,\left(\frac{r}{\kappa}\psi\right)\right](\eta,x)\\
&\qquad= \sum_{e\colon \|e\|=1}
\left(e^{\frac{r}{\kappa}[\psi(\eta,x+e)-\psi(\eta,x)]}
-1-\frac{r}{\kappa}[\psi(\eta,x+e)-\psi(\eta,x)]\right).
\end{aligned}
\ee
Taylor expansion of the r.h.s.\ of (\ref{main2.2-5}) gives that uniformly
in $\eta\in\Omega$ and $x\in\Z^d$,
\be{main2.2-7}
\left[\left(e^{-\frac{r}{\kappa} \psi} \Delta\,
e^{\frac{r}{\kappa}\psi}\right)
-\Delta\,\left(\frac{r}{\kappa}\psi\right)\right](\eta,x)
= \frac{r^2}{2\kappa^2} \sum_{e\colon \|e\|=1}
\Big(\psi(\eta,x+e)-\psi(\eta,x)\Big)^2+R_{\kappa,T}(\eta,x)
\ee
with
\be{main2.2-9}
\left|R_{\kappa,T}(\eta,x)\right|
\leq \frac{r^3}{6\kappa^3}\sum_{e\colon \|e\|=1}
\Big|\psi(\eta,x+e)-\psi(\eta,x)\Big|^3e^{o(1)}
\leq \frac{8 d r^3}{3\kappa^3} T^3e^{o(1)},
\ee
where we use (\ref{fdifest0}). Combining (\ref{main2.2-7}--\ref{main2.2-9}), we
arrive at (\ref{main2.2-3}).
\eprl

\subsubsection{Proof of Lemmas \ref{mainLem3}}
\label{main3}

\bpr
By (\ref{Gap1.69}), we have for all $\eta\in\Omega$ and $x\in\Z^d$,
\be{main3-7*}
\begin{aligned}
&\sum_{e\colon\,\|e\|=1}
\Big(\psi\big(\eta,x+e\big)
-\psi\big(\eta,x\big)\Big)^2\\
&\qquad=\sum_{e\colon\,\|e\|=1}\,\sum_{z_1,z_2\in\Z^d}\int_0^T du\, \int_0^T dv\,
\Big(p_{2du\onek}(z_1,x+e)-p_{2du\onek}(z_1,x)\Big)\\
&\qquad\qquad\qquad\times
\Big(p_{2dv\onek}(z_2,x+e)-p_{2dv\onek}(z_2,x)\Big)
\big[\eta(z_1)-\rho\big]\big[\eta(z_2)-\rho\big]\\
&\qquad=\sum_{z_1,z_2\in\Z^d} K^{\kappa,T}(z_1,z_2)\,
\big[\eta(z_1+x)-\rho\big]\big[\eta(z_2+x)-\rho\big],
\end{aligned}
\ee
where $K^{\kappa,T}\colon\,\Z^d\times\Z^d\mapsto\R$ is given by
\be{KTkappa}
K^{\kappa,T}(z_1,z_2)=\sum_{e:\|e\|=1}
\Big(\chi(z_1+e)-\chi(z_1)\Big)\Big(\chi(z_2+e)-\chi(z_2)\Big).
\ee
Therefore, for all $\kappa,T>0$,
\be{main3-9}
\begin{aligned}
&\limsup_{t\to\infty}\frac{1}{t}\log\E_{\,\nu_\rho,0}
\Bigg(\exp\Bigg[\frac{\alpha r}{2\kappa^2}\int_0^t ds\, \sum_{e:\|e\|=1}
\bigg(\psi\Big(\xi_{\frac{s}\kappa},X_s+e\Big)
-\psi\Big(\xi_{\frac{s}\kappa},X_s\Big)\bigg)^2\Bigg]\Bigg)\\
&\qquad=
\limsup_{t\to\infty}\frac{1}{t}\log\E_{\,\nu_\rho,0}
\Bigg(\exp\Bigg[\frac{\alpha r}{2\kappa^2}\int_0^t ds\, \sum_{z_1,z_2\in\Z^d}
K^{\kappa,T}(z_1,z_2)\,\\
&\qquad\qquad\qquad\qquad\qquad\qquad\qquad\times
\Big(\xi_\frac{s}\kappa(z_1+X_s)-\rho\Big)
\Big(\xi_\frac{s}\kappa(z_2+X_s)-\rho\Big)\Bigg]\Bigg)\\
&\qquad\leq
\limsup_{t\to\infty}\frac{1}{t}\log\E_{\,\nu_\rho}
\Bigg(\exp\Bigg[\frac{\alpha r}{2\kappa^2}\int_0^t ds\, \sum_{z_1,z_2\in\Z^d}
K^{\kappa,T}(z_1,z_2)\, \Big(\xi_\frac{s}\kappa(z_1)-\rho\Big)
\Big(\xi_\frac{s}\kappa(z_2)-\rho\Big)\Bigg]\Bigg),
\end{aligned}
\ee
where in the last line we use Lemma \ref{spec-rep} with
\be{main3-11}
W(\eta,x)=\frac{\alpha r}{2\kappa^2}\sum_{z_1,z_2\in\Z^d}
K^{\kappa,T}(z_1,z_2)\, \big[\eta(z_1+x)-\rho\big]
\big[\eta(z_2+x)-\rho\big],
\ee
which satisfies $W(\eta,x)=W(\tau_x\eta,0)$ as required in (\ref{spec-rep-2}). Splitting
the sum in the r.h.s.\ of (\ref{main3-9}) into its diagonal and off-diagonal part and
using the Cauchy-Schwarz inequality, we arrive at (\ref{main3-3}).
\epr

\subsubsection{Proof of Lemma \ref{mainLem4}}
\label{main4}

The proof of Lemma \ref{mainLem4} is based on the following two lemmas.
Recall (\ref{occtime}).

\bl{Klim}
For any $T>0$ there exists $C_T>0$, satisfying $\lim_{T\to\infty} C_T=0$,
such that
\be{Klim3}
\lim_{\kappa\to\infty}\|K_\diag^{\kappa,T}\|_1 = \frac{1}{d}G_d +C_T.
\ee
\el

\bl{ldpart}
For any $T>0$, $\alpha\in\R$ and $r>0$,
\be{ldpart-3}
\limsup_{t,\kappa\to\infty}\frac{\kappa^2}{t}\log \E_{\,\nu_\rho}
\bigg(\exp\bigg[\frac{\alpha r}{\kappa}(1-2\rho)\|K_\diag^{\kappa,T}\|_1
T_{t/\kappa}\bigg]\bigg)
\leq\alpha r \rho(1-2\rho)\lim_{\kappa\to\infty}  \|K_\diag^{\kappa,T}\|_1.
\ee
\el

\noindent
Before giving the proofs of Lemmas \ref{Klim}--\ref{ldpart}, we first
prove Lemma \ref{mainLem4}.

\bprl{mainLem4}
By Jensen's inequality, we have
\be{main4-7}
\begin{aligned}
&\E_{\,\nu_\rho}
\Bigg(\exp\Bigg[\frac{\alpha r}{\kappa^2}\int_0^t ds\,
\sum_{z\in\Z^d} K_\diag^{\kappa,T}(z)\Big(\xi_{\frac{s}\kappa}(z)-\rho\Big)^2\Bigg]\Bigg)\\
&\qquad\leq
\sum_{z\in\Z^d} \frac{K_\diag^{\kappa,T}(z)}{\|K_\diag^{\kappa,T}\|_1}\E_{\,\nu_\rho}
\Bigg(\exp\Bigg[\frac{\alpha r}{\kappa^2}\|K_\diag^{\kappa,T}\|_1\int_0^t ds\,
\Big(\xi_\frac{s}\kappa(z)-\rho\Big)^2\Bigg]\Bigg)\\
&\qquad=
\E_{\,\nu_\rho}
\Bigg(\exp\Bigg[\frac{\alpha r}{\kappa^2}\|K_\diag^{\kappa,T}\|_1\int_0^t ds\,
\Big(\xi_\frac{s}\kappa(0)-\rho\Big)^2\Bigg]\Bigg)\\
&\qquad=
\exp\Big[\frac{\alpha r}{\kappa^2}\rho^2\|K_\diag^{\kappa,T}\|_1 t\Big]\, \E_{\,\nu_\rho}
\Bigg(\exp\Bigg[\frac{\alpha r}{\kappa}(1-2\rho)\|K_\diag^{\kappa,T}\|_1\int_0^{t/\kappa} ds\,
\xi_s(0)\Bigg]\Bigg),
\end{aligned}
\ee
where the first equality uses the shift-invariance of $\nu_\rho$. Therefore
\be{main4-9}
\begin{aligned}
&\lim_{t,\kappa,T\to\infty}\frac{\kappa^2}{t}\log \E_{\,\nu_\rho}
\Bigg(\exp\Bigg[\frac{\alpha r}{\kappa^2}\int_0^t ds\,
\sum_{z\in\Z^d} K_\diag^{\kappa,T}(z)\Big(\xi_{\frac{s}\kappa}(z)-\rho\Big)^2\Bigg]\Bigg)\\
&\quad\leq
\lim_{\kappa,T\to\infty} \alpha r \rho^2 \|K_\diag^{\kappa,T}\|_1
+\lim_{t,\kappa,T\to\infty}\frac{\kappa^2}{t}\log \E_{\,\nu_\rho}
\bigg(\exp\bigg[\frac{\alpha r}{\kappa}(1-2\rho)\|K_\diag^{\kappa,T}\|_1
T_{\frac{t}{\kappa}}\bigg]\bigg).
\end{aligned}
\ee
Now use Lemmas \ref{Klim}--\ref{ldpart} to obtain (\ref{main4-3}).
\eprl

\bprl{Klim}
By (\ref{main3-5}), we have
\be{Klim5}
\begin{aligned}
\|K_\diag^{\kappa,T}\|_1
&= 2\sum_{\{x,y\}} \int_0^T du\, \int_u^T dv\,
\Big(p_{2du\onek}(0,y)-p_{2du\onek}(0,x)\Big)
\Big(p_{2dv\onek}(0,y)-p_{2dv\onek}(0,x)\Big)\\
&= -4 \int_0^T du\, \int_u^T dv\, \sum_{x\in\Z^d}
p_{2du\onek}(0,x)\Big[\Delta_1 p_{2dv\onek}(0,x)\Big]\\
&= -\frac{4}{\onek} \int_0^T du\, \int_u^T dv\, \sum_{x\in\Z^d}
p_{2du\onek}(0,x)\bigg[\frac{\partial}{\partial v} p_{2dv\onek}(0,x)\bigg],
\end{aligned}
\ee
where we recall the remark below (\ref{Tayl.23}). After performing the
integration w.r.t.\ the variable $v$, we get
\be{Klim7}
\begin{aligned}
\|K_\diag^{\kappa,T}\|_1
&=\frac{4}{\onek} \Bigg(
\int_0^T du\, \sum_{x\in\Z^d} p^2_{2du\onek}(0,x)
-\int_0^T du\, \sum_{x\in\Z^d} p_{2du\onek}(0,x)p_{2dT\onek}(0,x)\Bigg)\\
&=\frac{4}{\onek} \Bigg(
\int_0^T du\,\, p_{4du\onek}(0,0)
-\int_0^T du\,\, p_{2d(u+T)\onek}(0,0)\Bigg).
\end{aligned}
\ee
Hence
\be{kapnormlim}
\lim_{\kappa\to\infty} \|K_\diag^{\kappa,T}\|_1
= \frac{1}{d}\Bigg(\int_0^{2dT} du\,\,p_u(0,0) - \int_{2dT}^{4dT} du\,\,p_u(0,0)\Bigg),
\ee
which gives (\ref{Klim3}).
\eprl

\bprl{ldpart}
To derive (\ref{ldpart-3}), we use the large deviation principle for $(T_t)_{t\geq 0}$
stated in Section \ref{S2.1}. By Varadhan's Lemma we have, for all $\kappa,T>0$,
\be{ldpart-5}
\begin{aligned}
&\limsup_{t\to\infty}\frac{1}{t}\log \E_{\,\nu_\rho}
\bigg(\exp\bigg[\frac{\alpha r}{\kappa}(1-2\rho)\|K_\diag^{\kappa,T}\|_1
T_{t/\kappa}\bigg]\bigg)\\
&\qquad\qquad\qquad=
\frac1\kappa \max_{\beta\in[0,1]}\Big[
\frac{\alpha r}\kappa (1-2\rho)\|K_\diag^{\kappa,T}\|_1 \beta
-\Psi_d(\beta)\Big]\\
&\qquad\qquad\qquad\leq
\frac1\kappa \max_{\beta\in[0,1]}\Big[
\frac{\alpha r}\kappa (1-2\rho)\|K_\diag^{\kappa,T}\|_1 \beta
-\frac{1}{2G_d}\big(\sqrt{\beta}-\sqrt{\rho}\big)^2\Big],
\end{aligned}
\ee
where we use the quadratic lower bound in (\ref{rflb}). By Lemma \ref{Klim},
$(1/\kappa)\|K_\diag^{\kappa,T}\|_1 \downarrow 0$ as $\kappa\to\infty$ for any
$T>0$. Hence Lemma \ref{varform} can be applied to get (\ref{ldpart-3}).
\eprl

\subsubsection{Proof of Lemma \ref{mainLem5}}
\label{main5}

The proof of Lemma \ref{mainLem5} is based on the following two lemmas.
Recall (\ref{main3-5}). For $z_1,z_2\in\Z^d$ with $z_1\neq z_2$ and $\gamma\in\R$,
let
\be{hfun}
\begin{aligned}
h_{\gamma,\kappa}(z_1,z_2)
&=\limsup_{t\to\infty}\frac1t\log\E_{\,\nu_\rho}\bigg(
\exp\left[\frac{\gamma}{\kappa^2}\int_0^t ds\,\,
\Big(\xi_{\frac{s}\kappa}(z_1)-\rho\Big)
\Big(\xi_{\frac{s}\kappa}(z_2)-\rho\Big)\right]\bigg)\\
&=\limsup_{t\to\infty}\frac1{\kappa t}\log\E_{\,\nu_\rho}\bigg(
\exp\left[\frac{\gamma}{\kappa}\int_0^t ds\,\,
\big(\xi_s(z_1)-\rho\big)\big(\xi_s(z_2)-\rho\big)\right]\bigg).
\end{aligned}
\ee

\bl{normKoff}
For all $\kappa,T>0$,
\be{normKoff-3}
\|K_\off^{\kappa,T}\|_1 \leq 8dT^2.
\ee
\el

\bl{+bound}
For any $z_1,z_2\in\Z^d$ with $z_1\neq z_2$ and any $\gamma\in\R$,
\be{+bound-3}
\limsup_{\kappa\to\infty}\kappa^2 h_{\gamma,\kappa}(z_1,z_2)\leq 0.
\ee
\el

\noindent
Before giving the proof of Lemmas \ref{normKoff}--\ref{+bound}, we first prove Lemma
\ref{mainLem5}.

\bprl{mainLem5}
Let $K_\off^{\kappa,T;+}$ and $K_\off^{\kappa,T;-}$ denote, respectively, the positive 
and negative part of $K_\off^{\kappa,T}$. By the Cauchy-Schwarz inequality, we have
\be{main5-7}
\begin{aligned}
&\log\E_{\,\nu_\rho}
\Bigg(\exp\Bigg[\frac{\alpha r}{\kappa^2}\int_0^t ds\,
\sum_{{z_1,z_2\in\Z^d}\atop{z_1\neq z_2}} K_\off^{\kappa,T}(z_1,z_2)
\Big(\xi_{\frac{s}\kappa}(z_1)-\rho\Big)\Big(\xi_{\frac{s}\kappa}(z_2)-\rho\Big)\Bigg]\Bigg)\\
&\leq\frac12\log\E_{\,\nu_\rho}
\Bigg(\exp\Bigg[\frac{2\alpha r}{\kappa^2}
\sum_{{z_1,z_2\in\Z^d}\atop{z_1\neq z_2}} K_\off^{\kappa,T;+}(z_1,z_2) \int_0^t ds\,
\Big(\xi_{\frac{s}\kappa}(z_1)-\rho\Big)\Big(\xi_{\frac{s}\kappa}(z_2)-\rho\Big)\Bigg]\Bigg)\\
&\quad+\frac12\log\E_{\,\nu_\rho}
\Bigg(\exp\Bigg[-\frac{2\alpha r}{\kappa^2}
\sum_{{z_1,z_2\in\Z^d}\atop{z_1\neq z_2}} K_\off^{\kappa,T;-}(z_1,z_2) \int_0^t ds\,
\Big(\xi_{\frac{s}\kappa}(z_1)-\rho\Big)\Big(\xi_{\frac{s}\kappa}(z_2)-\rho\Big)\Bigg]\Bigg).
\end{aligned}
\ee
We estimate the first term in the r.h.s.\ of (\ref{main5-7}). For $R>0$, let
\be{+3}
B_R =\{(z_1,z_2)\in\Z^d\times\Z^d : \|z_1\|+\|z_2\|\leq R\}.
\ee
Then
\be{+5}
\begin{aligned}
&\E_{\,\nu_\rho}
\Bigg(\exp\Bigg[\frac{2\alpha r}{\kappa^2}\int_0^t ds\,
\sum_{{z_1,z_2\in\Z^d}\atop{z_1\neq z_2}} K_\off^{\kappa,T;+}(z_1,z_2)
\Big(\xi_{\frac{s}\kappa}(z_1)-\rho\Big)\Big(\xi_{\frac{s}\kappa}(z_2)-\rho\Big)\Bigg]\Bigg)\\
&\quad\leq \exp\Bigg[\frac{2|\alpha| r t}{\kappa^2}
\sum_{{(z_1,z_2)\in B_R^{\c}}\atop{z_1\neq z_2}}
K_\off^{\kappa,T;+}(z_1,z_2)\Bigg]\\
&\qquad\times
\E_{\,\nu_\rho}
\Bigg(\exp\Bigg[\frac{2\alpha r}{\kappa^2}\int_0^t ds\,
\sum_{{(z_1,z_2)\in B_R}\atop{z_1\neq z_2}} K_\off^{\kappa,T;+}(z_1,z_2)
\Big(\xi_{\frac{s}\kappa}(z_1)-\rho\Big)\Big(\xi_{\frac{s}\kappa}(z_2)-\rho\Big)\Bigg]\Bigg).
\end{aligned}
\ee
Applying Jensen's inequality, we get
\be{+5*}
\begin{aligned}
&\frac{\kappa^2}{t}\log\E_{\,\nu_\rho}
\Bigg(\exp\Bigg[\frac{2\alpha r}{\kappa^2}\int_0^t ds\,
\sum_{{z_1,z_2\in\Z^d}\atop{z_1\neq z_2}} K_\off^{\kappa,T;+}(z_1,z_2)
\Big(\xi_{\frac{s}\kappa}(z_1)-\rho\Big)\Big(\xi_{\frac{s}\kappa}(z_2)-\rho\Big)\Bigg]\Bigg)\\
&\leq 2|\alpha| r \sum_{{(z_1,z_2)\in B_R^{\c}}\atop{z_1\neq z_2}}
K_\off^{\kappa,T;+}(z_1,z_2)
+\frac{\kappa}{t/\kappa}\log \sum_{{(z_1,z_2)\in B_R}\atop{z_1\neq z_2}}\frac{
K_\off^{\kappa,T;+}(z_1,z_2)}{\big\|K_{\off;R}^{\kappa,T;+}\big\|_1}\\
&\qquad\qquad\qquad\qquad\qquad\quad\times \E_{\,\nu_\rho}
\Bigg(\exp\Bigg[\frac{2\alpha r}{\kappa}\big\|K_{\off;R}^{\kappa,T;+}\big\|_1
\int_0^{t/\kappa} ds\,
\Big(\xi_{s}(z_1)-\rho\Big)\Big(\xi_{s}(z_2)-\rho\Big)\Bigg]\Bigg),
\end{aligned}
\ee
where
\be{+7}
\|K_{\off;R}^{\kappa,T;+}\big\|_1
=\sum_{{(z_1,z_2)\in B_R}\atop{z_1\neq z_2}}K_{\off}^{\kappa,T;+}(z_1,z_2).
\ee
By Lemma \ref{+bound} (with $\gamma=2\alpha r \|K_{\off;R}^{\kappa,T;+}\|_1$), the second term in the 
r.h.s.\ of (\ref{+5*}) is asymptotically bounded by above by zero (as $t\to\infty$) for any $\kappa,T>0$,
$\alpha\in\R$ and $r>0$, and any $R$ finite. The first term in the r.h.s.\ of (\ref{+5*}) does not depend
on $t$ and, by Lemma \ref{normKoff}, tends to zero as $R\to\infty$. This shows that the first term in the
r.h.s\ of (\ref{main5-7}) yields a zero contribution. The same is true for the second term by the same
argument. This completes the proof of (\ref{main5-3}).
\eprl

\bprl{normKoff}
The claim follows from (\ref{main3-5}--\ref{main3-7}).
\eprl

\bprl{+bound}
The proof of Lemma \ref{+bound} is long, since it is based on three further lemmas.
Let $z_1,z_2\in\Z^d$ with $z_1\neq z_2$. Without loss of generality, we may assume that
\be{Hdef1}
z_1\in H^- \quad\text{and}\quad z_2\in H^+
\ee
with
\be{Hdef2}
H^-=\{z\in\Z^d\colon\,z^1\leq 0\} \quad\text{and}\quad H^+=\{z\in\Z^d\colon\,z^1> 0\}.
\ee
Let
\be{hfuns}
\begin{aligned}
h_{\gamma,\kappa}^{-}(z_1)&=\limsup_{t\to\infty}\frac{1}{3\kappa t}\log\E_{\,\nu_\rho}^\IRW\bigg(
\exp\left[-\frac{3\gamma}{\kappa}\rho\int_0^t ds\,\,
\tilde\xi_s^-(z_1)\right]\bigg),\\[0.3cm]
h_{\gamma,\kappa}^{+}(z_2)&=\limsup_{t\to\infty}\frac{1}{3\kappa t}\log\E_{\,\nu_\rho}^\IRW\bigg(
\exp\left[-\frac{3\gamma}{\kappa}\rho\int_0^t ds\,\,
\tilde\xi_s^+(z_2)\right]\bigg),\\[0.3cm]
h_{\gamma,\kappa}^\pm(z_1,z_2)&=\limsup_{t\to\infty}\frac{1}{3\kappa t}\log\E_{\,\nu_\rho}^\IRW\bigg(
\exp\left[\frac{3\gamma}{\kappa}\int_0^t ds\,\,
\tilde\xi_s^-(z_1)\tilde\xi_s^+(z_2)\right]\bigg),
\end{aligned}
\ee
where $(\tilde\xi_t^-)_{t\geq 0}$ and $(\tilde\xi_t^+)_{t\geq 0}$ are independent IRW's on $H^-$
and $H^+$, respectively, with transition kernels $p^-(\cdot,\cdot)$ and $p^+(\cdot,\cdot)$
corresponding to simple random walks stepping at rate $1$ such that steps outside $H^-$ and $H^+$,
respectively, are suppressed.

\bl{+bound1}
For all $\kappa>0$, $z_1\in H^-$, $z_2\in H^+$ and $\gamma\in\R$,
\be{+bound1-5}
h_{\gamma,\kappa}(z_1,z_2)
\leq \frac{\gamma}{\kappa^2}\rho^2
+h_{\gamma,\kappa}^{-}(z_1)
+h_{\gamma,\kappa}^{+}(z_2)
+h_{\gamma,\kappa}^\pm(z_1,z_2).
\ee
\el

\bl{+bound2}
For all $\gamma\in\R$,
\be{+bound2-5}
\limsup_{\kappa\to\infty} \kappa^2 \sup_{z_1\in H^{-}} h_{\gamma,\kappa}^{-}(z_1)
\leq -\gamma\rho^2
\,\,\,\,\text{and}\,\,\,\,
\limsup_{\kappa\to\infty} \kappa^2 \sup_{z_2\in H^{+}} h_{\gamma,\kappa}^{+}(z_2)
\leq -\gamma\rho^2.
\ee
\el

\bl{+bound3}
For all $\gamma\in\R$,
\be{+bound3-5}
\limsup_{\kappa\to\infty} \kappa^2 \sup_{{z_1\in H^-}\atop{z_2\in H^+}}
h_{\gamma,\kappa}^\pm(z_1,z_2)\leq \gamma\rho^2.
\ee
\el

\noindent
Combining (\ref{+bound1-5}--\ref{+bound3-5}), we get (\ref{+bound-3}).
\eprl

\bprl{+bound1}
Similarly as in the proof of Lemma \ref{spec-rep}, by cutting the bonds connecting
$H^-$ and $H^+$ in the analogue of the variational formula of Proposition \ref{RR},
we get
\be{+bound1-7}
\begin{aligned}
&h_{\gamma,\kappa}(z_1,z_2)
\leq\limsup_{t\to\infty}\frac1{\kappa t}\log\E_{\,\nu_\rho}\bigg(
\exp\left[\frac{\gamma}{\kappa}\int_0^t ds\,\,
\big(\xi_s^-(z_1)-\rho\big)
\big(\xi_s^+(z_2)-\rho\big)\right]\bigg)\\
&\qquad\quad=\limsup_{t\to\infty}\frac1{\kappa t}\log\E_{\,\nu_\rho}\bigg(
\exp\left[\frac{\gamma}{\kappa}\int_0^t ds\,\,
\Big(\rho^2 -\rho\, \xi_s^-(z_1)+\rho\, \xi_s^+(z_2)+\xi_s^-(z_1)\xi_s^+(z_2)
\Big)\right]\bigg),
\end{aligned}
\ee
where $(\xi_t^-)_{t\geq 0}$ and $(\xi_t^+)_{t\geq 0}$ are independent
exclusion processes in $H^-$ and $H^+$, respectively, obtained from $(\xi_t)_{t\geq 0}$
by suppressing jumps between $H^-$ and $H^+$. Applying H\"older's inequality in the r.h.s.\
of (\ref{+bound1-7}) to separate terms, we obtain
\be{+bound1-9}
\begin{aligned}
h_{\gamma,\kappa}(z_1,z_2)
\leq\frac{\gamma}{\kappa^2}\rho^2
&+\limsup_{t\to\infty}\frac{1}{3\kappa t}\log\E_{\,\nu_\rho}\bigg(
\exp\left[-\frac{3\gamma}{\kappa}\rho\int_0^t ds\,\,
\xi_s^-(z_1)\right]\bigg)\\
&+\limsup_{t\to\infty}\frac{1}{3\kappa t}\log\E_{\,\nu_\rho}\bigg(
\exp\left[-\frac{3\gamma}{\kappa}\rho\int_0^t ds\,\,
\xi_s^+(z_2)\right]\bigg)\\
&+\limsup_{t\to\infty}\frac{1}{3\kappa t}\log\E_{\,\nu_\rho}\bigg(
\exp\left[\frac{3\gamma}{\kappa}\int_0^t ds\,\,
\xi_s^-(z_1)\xi_s^+(z_2)\right]\bigg).
\end{aligned}
\ee
In order to get (\ref{+bound1-5}), we apply Proposition \ref{IRW-comp} to the last three terms 
in the r.h.s.\ of (\ref{+bound1-9}). For the first two terms, pick, respectively $K(z,s)= 
-(3\gamma/\kappa)\rho\, 1_{z_1}$ and $K(z,s)= -(3\gamma/\kappa)\rho\, 1_{z_2}(z)$. For the 
last term, we have to apply Proposition \ref{IRW-comp} twice, once for the exclusion process 
$(\xi_t^+)_{t\geq 0}$ on $H^+$ with $K(z,s)=-(3\gamma/\kappa)\, \xi_s^-(z_1)\, 1_{z_2}(z)$ and 
once for the exclusion process $(\xi_t^-)_{t\geq 0}$ on $H^-$ with $K(z,s)=-(3\gamma/\kappa)\,
\tilde\xi_s^-(z_2)\,1_{z_1}(z)$. Here, we in fact apply a \emph{modification} of Proposition 
\ref{IRW-comp} by considering $(\xi_t^-)_{t\geq 0}$ and $(\xi_t^+)_{t\geq 0}$ on $\Z^d$ with 
particles not moving on $H^+$ and $H^-$, respectively. See the proof of Proposition \ref{IRW-comp}
in Appendix \ref{A} to verify that this modification holds true.
\eprl

\bprl{+bound2}
We prove the second line of (\ref{+bound2-5}). The first line follows by symmetry.

Let
\be{+bound2-8}
H^+_\eta = \{x\in H^+\colon\,\eta(x)=1\},\qquad \eta\in\Omega.
\ee
Fix $z\in H^+$. Then
\be{+bound2-9}
\E_{\,\nu_\rho}^{\IRW}\bigg(
\exp\left[-\frac{3\gamma}{\kappa}\rho\int_0^t ds\,\, \tilde\xi_s^+(z)\right]\bigg)
= \int_\Omega \nu_\rho(d\eta)\, \prod_{x\in H^+_\eta}\E_{\,x}^{\RW,+}
\left(\exp\bigg[-\frac{3\gamma}{\kappa}\rho\int_0^t ds\,\,
1_{z}(Y_s^+) \bigg]\right),
\ee
where $\E_x^{\RW,+}$ is expectation w.r.t.\ simple random walk $Y^{+}=(Y^{+}_t)_{t\geq 0}$ on
$H^+$ with transition kernel $p^+(\cdot,\cdot)$ and step rate $1$ starting from $Y_0^+=x\in H^+$.
Using that $\nu_\rho$ is the Bernoulli product measure with density $\rho$, we get
\be{+bound2-10}
\begin{aligned}
&\E_{\,\nu_\rho}^{\IRW}\bigg(
\exp\left[-\frac{3\gamma}{\kappa}\rho\int_0^t ds\,\, \tilde\xi_s^+(z)\right]\bigg)\\
&\qquad = \int_\Omega \nu_\rho(d\eta) \prod_{x\in H^+}\E_{\,x}^{\RW,+}
\left(\exp\bigg[-\eta(x)\frac{3\gamma}{\kappa}\rho\int_0^t ds\,\,
1_{z}(Y_s^+) \bigg]\right)\\
&\qquad = \prod_{x\in H^+}\Big(1-\rho + \rho v(x,t)\Big)
\leq \exp\Big[\rho\sum_{x\in H^+}\big(v(x,t)-1\big)\Big]
\end{aligned}
\ee
with
\be{+bound2-11}
v(x,t)
=\E_{\,x}^{\RW,+}
\bigg(\exp\bigg[-\frac{3\gamma}{\kappa}\rho\int_0^t ds\,\,
1_{z}(Y_s^+) \bigg]\bigg).
\ee

By the Feynman-Kac formula, $v\colon\,H^+\times[0,\infty)\to \R$ is the solution of
the Cauchy problem
\be{+bound2-13}
\frac{\partial}{\partial t} v(x,t) = \frac1{2d} \Delta^+ v(x,t)
-\left\{\frac{3\gamma}{\kappa}\,\rho\,1_{z}(x)\right\}\,v(x,t),
\qquad v(\cdot,0) \equiv 1,
\ee
where
\be{+bound2-14}
\Delta^+ v(x,t)=\sum_{{y\in H^+} \atop {\|y-x\|=1}} [v(y,t)-v(x,t)], \qquad x\in H^+.
\ee
Put
\be{+bound2-15}
w(x,t)=v(x,t)-1.
\ee
Then $w\colon\,H^+\times[0,\infty)\to \R$ is the solution of the Cauchy problem
\be{+bound2-17}
\frac{\partial w}{\partial t}(x,t) = \frac1{2d} \Delta^+ w(x,t)
-\left\{\frac{3\gamma}{\kappa}\,\rho\,1_{z}(x)\right\}[w(x,t)+1],
\qquad w(\cdot,0)\equiv 0.
\ee
Since $\sum_{x\in H^+}\Delta^+f(x)=0$ for all $f\colon\,H^+\to\R$, (\ref{+bound2-17})
gives
\be{+bound2-19}
\frac{\partial}{\partial t} \sum_{x\in H^+} w(x,t)
=-\frac{3\gamma}{\kappa}\, \rho[w(z,t)+1].
\ee
After integrating (\ref{+bound2-19}) w.r.t.\ $t$, we obtain
\be{+bound2-21}
\sum_{x\in H^+} w(x,t)
=-\frac{3\gamma}{\kappa}\rho\, t
-\frac{3\gamma}{\kappa}\rho \int_0^t ds\,\, w(z,s).
\ee
Combining (\ref{hfuns}), (\ref{+bound2-10}), (\ref{+bound2-15}) and (\ref{+bound2-21}),
we arrive at
\be{+bound2-23}
h^+_{\gamma,\kappa}(z) \leq -\frac{\gamma}{\kappa^2}\rho^2
\bigg(1+\lim_{t\to\infty}\frac1t\int_0^t ds\,\, w(z,s)\bigg).
\ee
The limit in the r.h.s.\ exists since, by (\ref{+bound2-11}) and (\ref{+bound2-15}),
$w(z,t)$ is monotone in $t$.

We will complete the proof by showing that the second term in the r.h.s.\ of (\ref{+bound2-23})
tends to zero as $\kappa\to\infty$. This will rely on the following lemma, the proof of which
is deferred to the end of this section.

\bl{G+normbd}
Let $G^+(x,y)$ be the Green kernel on $H^+$ associated with $p_t^+(x,y)$. Then
$\|G^+\|_\infty\leq 2G_d<\infty$.
\el
Return to (\ref{+bound2-11}). If $\gamma>0$, then by Jensen's inequality we have
\be{+bound2-25}
1 \geq v(x,t)
\geq\exp\bigg[-\frac{3\gamma}{\kappa}\rho\int_0^t ds\,\,p_s^+(x,z) \bigg]
\geq\exp\bigg[-\frac{3\gamma}{\kappa}\rho\|G^+\|_\infty \bigg],
\ee
where $\|G^+\|_\infty <\infty$ by Lemma \ref{G+normbd}. To deal with the case 
$\gamma\leq 0$, let $\cG^+$ denote the Green operator acting on functions
$V\colon\,H^+\to [0,\infty)$ as
\be{+bound2-29}
(\cG^+ V)(x)=\sum_{y\in H^+}G^+(x,y)V(y),\qquad x\in H^+.
\ee
We have
\be{+bound2-31}
\bigg\|\cG^+\Big(\frac{3\gamma}{\kappa}\rho\, 1_{z}\Big)\bigg\|_{\infty}
\leq\frac{3|\gamma|}{\kappa}\rho\,\|G^+\|_\infty.
\ee
The r.h.s.\ tends to zero as $\kappa\to\infty$. Hence Lemma 8.2.1 in G\"artner and den
Hollander \cite{garhol04} can be applied to (\ref{+bound2-11}) for $\kappa$ large enough,
to yield
\be{+bound2-33}
1 \leq v(x,t)
\leq \frac1{1-\frac{3|\gamma|}{\kappa}\rho\,\|G^+\|_\infty} \da 1
\quad\text{as}\quad \kappa\to\infty.
\ee
Therefore, combining (\ref{+bound2-25}) and (\ref{+bound2-33}), we see that for all
$\gamma\in\R$ and $\delta\in(0,1)$ there exists $\kappa_0=\kappa_0(\gamma,\delta)$ such that
\be{+bound2-35}
\|v-1\|_\infty\leq \delta \qquad \forall\,\kappa>\kappa_0.
\ee
By (\ref{+bound2-15}--\ref{+bound2-17}), we have
\be{+bound2-37}
w(z,t)=-\frac{3\gamma}{\kappa}\rho\int_0^t ds\,\,
\E_{\,z}^{\RW,+}\bigg(1_z(Y_s^+)\, v(Y_s^+,t-s)\bigg).
\ee
Via (\ref{+bound2-35}) it therefore follows that
\be{+bound2-39}
-\frac{3\gamma}{\kappa}\rho(1\pm\delta) G^+(z,z)
\leq \lim_{t\ra\infty}\frac1t \int_0^t ds\,\, w(z,s)
\leq -\frac{3\gamma}{\kappa}\rho(1\mp\delta) G^+(z,z)
\qquad \forall\,\kappa>\kappa_0,
\ee
where the choice of $+$ or $-$ in front of $\delta$ depends on the sign of $\gamma$.
The latter shows that the second term in the r.h.s.\ of (\ref{+bound2-23}) is $O(1/\kappa)$.
This proves (\ref{+bound2-5}).
\eprl

\bprl{+bound3}
The proof is similar to that of Lemma \ref{+bound2}. Let
\be{+bound3-8}
\begin{aligned}
H^+_\eta &= \{x\in H^+\colon\, \eta(x)=1\},\quad \eta\in\Omega,\\[0.3cm]
H^-_\eta &= \{x\in H^-\colon\, \eta(x)=1\},\quad \eta\in\Omega.
\end{aligned}
\ee
Fix $z_1\in H^-$ and $z_2\in H^+$. Then
\be{+bound3-9}
\begin{aligned}
&\E_{\,\nu_\rho}^\IRW\bigg(
\exp\left[\frac{3\gamma}{\kappa}\int_0^t ds\,\,
\tilde\xi_s^-(z_1)\tilde\xi_s^+(z_2)\right]\bigg)\\
&\qquad= \int_\Omega \nu_\rho(d\eta)\, \prod_{x\in H^-_\eta}\prod_{y\in H^+_\eta}
\E_{\,x}^{\RW,-}\, \E_{\,y}^{\RW,+}\left(\exp\bigg[\frac{3\gamma}{\kappa}
\int_0^t ds\,\, 1_{(z_1,z_2)}(Y_s^-,Y_s^+) \bigg]\right),
\end{aligned}
\ee
where $Y^-$ on $H^-$ and $Y^+$ on $H^+$ are simple random walks with step rate $1$
and transition kernel $p^-(\cdot,\cdot)$ and $p^+(\cdot,\cdot)$ starting from
$Y_0^-=x\in H^-$ and $Y_0^+=y\in H^+$, respectively. Using that $\nu_\rho$
is the Bernoulli product measure with density $\rho$, we get
\be{+bound3-9*}
\begin{aligned}
&\E_{\,\nu_\rho}^\IRW\bigg(
\exp\left[\frac{3\gamma}{\kappa}\int_0^t ds\,\,
\tilde\xi_s^-(z_1)\tilde\xi_s^+(z_2)\right]\bigg)\\
&\qquad= \int_\Omega \nu_\rho(d\eta) \prod_{x\in H^-}\prod_{y\in H^+}
\E_{\,x}^{\RW,-}\, \E_{\,y}^{\RW,+}
\left(\exp\bigg[\eta(x)\eta(y)\frac{3\gamma}{\kappa}\int_0^t ds\,\,
1_{(z_1,z_2)}(Y_s^-,Y_s^+) \bigg]\right)\\
&\qquad= \prod_{x\in H^-}\prod_{y\in H^+}\Big(1-\rho^2+\rho^2 v(z_1,z_2;t)\Big)
\leq \exp\Big[\rho^2\sum_{x\in H^-}\sum_{y\in H^+}\big(v(z_1,z_2;t)-1\big)\Big]
\end{aligned}
\ee
with
\be{+bound3-11}
v(z_1,z_2;t)
=\E_{\,x}^{\RW,-}\, \E_{\,y}^{\RW,+}
\bigg(\exp\bigg[\frac{3\gamma}{\kappa}\int_0^t ds\,\,
\int_0^t ds\,\,
1_{(z_1,z_2)}(Y_s^-,Y_s^+) \bigg]\bigg).
\ee

By the Feynman-Kac formula, $v\colon\,(H^-\times H^+)\times [0,\infty)\to\R$ is the
solution of the Cauchy problem
\be{+bound3-13}
\frac{\partial}{\partial t} v(x,y;t) = \frac1{2d} \big(\Delta^- + \Delta^+\big) v(x,y;t)
+\left\{\frac{3\gamma}{\kappa}\, 1_{z_1,z_2}(x,y)\right\}\,v(x,y;t),
\qquad v(\cdot,\cdot;0) \equiv 1,
\ee
where
\be{+bound3-14}
\begin{aligned}
\Delta^- v(x;t)&=\sum_{{y\in H^-} \atop {\|y-x\|=1}} [v(y,t)-v(x,t)],
\qquad x\in H^-,\\[0.3cm]
\Delta^+ v(x;t)&=\sum_{{y\in H^+} \atop {\|y-x\|=1}} [v(y,t)-v(x,t)],
\qquad x\in H^+.
\end{aligned}
\ee
Put
\be{+bound3-15}
w(x,y;t)=v(x,y;t)-1.
\ee
Then, $w\colon\,(H^-\times H^+)\times [0,\infty)\to \R$ is the solution of the Cauchy problem
\be{+bound3-17}
\frac{\partial w}{\partial t}(x,y;t) = \frac1{2d} (\Delta^- + \Delta^+) w(x,y;t)
+\left\{\frac{3\gamma}{\kappa}\,1_{(z_1,z_2)}(x,y)\right\}\,[w(x,y;t)+1],
\qquad
w(\cdot,\cdot\,;0)\equiv 0.
\ee
By (\ref{+bound3-14}) and (\ref{+bound3-17}),
\be{+bound3-19}
\frac{\partial}{\partial t} \sum_{x\in H^-}\sum_{y\in H^+} w(x,y;t)
=\frac{3\gamma}{\kappa}[w(z_1,z_2;t)+1].
\ee
After integrating (\ref{+bound3-19}) w.r.t.\ $t$, we obtain
\be{+bound3-21}
\sum_{x\in H^-}\sum_{y\in H^+} w(x,y;t)
=\frac{3\gamma}{\kappa}\, t
+\frac{3\gamma}{\kappa} \int_0^t ds\,\, w(z_1,z_2;s).
\ee
Combining (\ref{hfuns}), (\ref{+bound3-9*}), (\ref{+bound3-15}) and (\ref{+bound3-21}),
we arrive at
\be{+bound3-23}
h_{\gamma,\kappa}^\pm(z_1,z_2) \leq \frac{\gamma}{\kappa^2}\rho^2
\bigg(1+\lim_{t\to\infty}
\frac1t\int_0^t ds\,\, w(z_1,z_2;s)\bigg).
\ee
The limit in the r.h.s.\ exists, since $w(z_1,z_2;s)$ is monotone in $s$.

We will complete the proof by showing that the second term in the r.h.s.\ of
(\ref{+bound3-23}) tends to zero as $\kappa\to\infty$. Return to (\ref{+bound3-11}).
If $\gamma\leq 0$, then by Jensen's inequality we have
\be{+bound3-25}
1 \geq v(x,y;t)
\geq\exp\bigg[-\frac{3|\gamma|}{\kappa}\int_0^\infty ds\,\,
p_s^-(x,z_1)\, p_s^+(y,z_2) \bigg]
\geq\exp\bigg[-\frac{3|\gamma|}{\kappa}\,
\Big(\|G^-\|_\infty\wedge\|G^+\|_\infty\Big) \bigg],
\ee
where $\|G^-\|_\infty,\|G^+\|_\infty<\infty$ by Lemma \ref{G+normbd}.
To deal with the case $\gamma>0$, let $\cG^\pm$ denote the Green operator acting on
functions $V\colon\,H^-\times H^+\to [0,\infty)$ as
\be{+bound3-29}
(\cG^\pm V)(x,y)=\sum_{{a\in H^-}\atop{b\in H^+}}G^\pm(x,y;a,b)V(a,b),
\qquad x\in H^-,\, y\in H^+,
\ee
where
\be{+bound3-30}
G^\pm(x,y;a,b)=\int_0^\infty ds\,\, p_s^-(x,a)\, p_s^+(y,b).
\ee
We have
\be{+bound3-31}
\bigg\|\cG^\pm\Big(\frac{3\gamma}{\kappa}\,1_{(z_1,z_2)}\Big)\bigg\|_{\infty}
\leq\frac{3\gamma}{\kappa}\, \|G^\pm\|_\infty
\leq\frac{3\gamma}{\kappa}\, \Big(\|G^-\|_\infty\wedge\|G^+\|_\infty\Big).
\ee
The r.h.s.\ tends to zero as $\kappa\to\infty$. Hence Lemma 8.2.1 in G\"artner and
den Hollander \cite{garhol04} can be applied to (\ref{+bound3-11}) for $\kappa$ large
enough, to yield
\be{+bound3-33}
1 \leq v(x,t)
\leq \frac1{1-\frac{3\gamma}{\kappa}
\Big(\|G^-\|_\infty\wedge\|G^+\|_\infty\Big)} \da 1\quad
\text{as}\quad \kappa\to\infty.
\ee
Therefore, combining (\ref{+bound3-25}) and (\ref{+bound3-33}), we see that for all
$\gamma\in\R$ and $\delta>0$ there exists $\kappa_0=\kappa_0(\gamma,\delta)$ such that
\be{+bound3-35}
\|v-1\|_\infty\leq \delta \qquad \forall\,\kappa>\kappa_0.
\ee
By (\ref{+bound3-15}--\ref{+bound3-17}), we have
\be{+bound2-37*}
w(z_1,z_2;t)=\frac{3\gamma}{\kappa}\int_0^t ds\,\,
\E_{\,z_1}^{\RW,-}\, \E_{\,z_2}^{\RW,+}\bigg(1_{(z_1,z_2)}(Y_s^-,Y_s^+)\,
v(Y_s^-,Y_s^+;t-s)\bigg).
\ee
Via (\ref{+bound3-35}) it therefore follows that for all $\kappa>\kappa_0$,
\be{+bound3-39}
\frac{3\gamma}{\kappa}(1\pm\delta) G^\pm(z_1,z_1;z_2,z_2)
\leq \lim_{t\to\infty}\frac1t \int_0^t ds\,\, w(z_1,z_2;s)
\leq \frac{3\gamma}{\kappa}(1\mp\delta) G^\pm(z_1,z_1;z_2,z_2).
\ee
Combining (\ref{+bound3-23}) and (\ref{+bound3-39}), we arrive at (\ref{+bound3-5}).
\eprl

\bprl{G+normbd}
We have $G^+(x,y)=\sum_{n=0}^\infty p_n^+(x,y)$, $x,y\in H^+$, with $p_n^+(x,y)$ the 
$n$-step transition probability of simple random walk on $H^+$ whose steps outside $H^+$ 
are suppressed (i.e., the walk pauses when it attempts to leave $H^+$). Let $p_n(x,y)$
be the $n$-step transition probability of simple random walk on $\Z^d$. Then 
\be{p+pbd}
p_n^+(x,y) \leq 2 p_n(x,y), \qquad x,y\in H^+,\, n\in\N_0.
\ee
Indeed, if we reflect simple random walk in the $(d-1)$-dimensional hyperplane between 
$H^+$ and its complement, then we obtain precisely the random walk that pauses when it
attempts to leave $H^+$. Hence, we have $p_n^+(x,y) = p_n(x,y)+p_n(x,y*)$, $x,y\in H^+$,
$n\in\N_0$, with $y*$ the reflection image of $y$. Since $p_n(x,y*) \leq p_n(x,y)$,
$x,y\in H^+$, the claim in (\ref{p+pbd}) follows. Sum on $n$, to get  $G^+(x,y) \leq 2G(x,y)$, 
$x,y\in H^+$. Now use that $G(x,y) \leq G(0,0) = G_d$, $x,y \in \Z^d$.
\eprl

\subsection{Proof of Proposition \ref{remainLem1}}
\label{S4.6}

The proof of Proposition \ref{remainLem1} is given in Section \ref{remain1}
subject to three lemmas. The latter are proved in Sections \ref{remain2}--\ref{remain4},
respectively. The first two lemmas are valid for $d\geq 3$, the third for $d\geq 4$.

\subsubsection{Proof of Proposition \ref{remainLem1}}
\label{remain1}

\bl{remainLem2}
For all $t\geq 0$, $\kappa,T>0$ and $\alpha\in\R$,
\be{remain2-3}
\E_{\,\nu_\rho,0}\bigg(\exp\bigg[
\frac{\alpha}{\kappa}\int_0^t ds\, \left(\cP_T\phi\right)(Z_s) \bigg]\bigg)
\leq\ES_{\,0}\bigg(\exp\bigg[
\frac{\alpha}{\kappa}\rho\int_0^t ds\, \sum_{x\in\Z^d} p_{2dT\onek}(X_{t-s},x)\,
w^{(t)}(x,s) \bigg]\bigg),
\ee
where $w^{(t)}\colon\,\Z^d\times[0,t)\to \R$ is the solution of the Cauchy problem
\be{remain2-5}
\frac{\partial w^{(t)}}{\partial s}(x,s) = \frac1{2d\kappa} \Delta w^{(t)}(x,s)
+ \frac\alpha\kappa\, p_{2dT\onek}(X_{t-s},x)[w^{(t)}(x,s)+1],
\qquad
w^{(t)}(\cdot,0)\equiv 0.
\ee
\el

\bl{remainLem3}
For all $t\geq 0$, $\kappa>0$, $T$ large enough and $\alpha\in\R$,
\be{remain3-3}
\begin{aligned}
&\ES_{\,0}\bigg(\exp\bigg[
\frac{\alpha}{\kappa}\rho\int_0^t ds\,\sum_{x\in\Z^d} p_{2dT\onek}(X_{t-s},x)\,
w^{(t)}(x,s) \bigg]\bigg)\\
&\qquad\leq\ES_{\,0}\bigg(\exp\bigg[
\frac{2\alpha^2}{\kappa^2}\rho\int_0^t ds\, \int_s^t du\,
p_{\frac{u-s}{\kappa}+4dT\onek}(X_u,X_s) \bigg]\bigg).
\end{aligned}
\ee
\el

\bl{remainLem4}
If $d\geq 4$, then for any $\alpha\in\R$,
\be{remain4-3}
\lim_{T,\kappa,t\to\infty}\frac{\kappa^2}{t}\log\ES_{\,0}\bigg(\exp\bigg[
\frac{2\alpha^2}{\kappa^2}\rho\int_0^t ds\, \int_s^t du\,
p_{\frac{u-s}{\kappa}+4dT\onek}(X_u,X_s) \bigg]\bigg)=0.
\ee
\el

\noindent
Lemmas \ref{remainLem2}--\ref{remainLem4} clearly imply (\ref{remain1-3}).

\subsubsection{Proof of Lemma \ref{remainLem2}}
\label{remain2}

For all $t\geq 0$, $\kappa,T>0$ and $\alpha\in\R$, let $v^{(t)}\colon\,\Z^d\times[0,t)\to \R$ 
be such that
\be{remain3-6}
v^{(t)}(x,s) = w^{(t)}(x,s)+1,
\ee
where $w^{(t)}$ is defined by (\ref{remain2-5}). Then $v^{(t)}$ is the solution of the 
Cauchy problem
\be{remain3-7**}
\frac{\partial v^{(t)}}{\partial s}(x,s) = \frac1{2d\kappa} \Delta v^{(t)}(x,s)
+ \frac{\alpha}\kappa\, p_{2dT\onek}(X_{t-s},x)\,v^{(t)}(x,s),
\qquad
v^{(t)}(\cdot,0)\equiv 1,
\ee
and has the representation
\be{remain3-9}
v^{(t)}(x,s)=\E_{\,x}^\RW\bigg(\exp\bigg[\frac{\alpha}{\kappa}\int_0^s du\,
p_{2dT\onek}\Big(X_{t-s+u},Y_{\frac{u}\kappa}\Big)\bigg]\bigg).
\ee

\bpr
By (\ref{lb2.5}) and (\ref{Gap1.73}), we have
\be{remain2-9}
\E_{\,\nu_\rho,0}\bigg(\exp\bigg[
\frac{\alpha}{\kappa}\int_0^t ds\, \left(\cP_T\phi\right)(Z_s) \bigg]\bigg)
=\E_{\,\nu_\rho,0}\bigg(\exp\bigg[
\frac{\alpha}{\kappa}\sum_{z\in\Z^d}\int_0^t ds\, p_{2dT\onek}(X_s,z)
\Big(\xi_\frac{s}{\kappa}(z)-\rho\Big)\bigg]\bigg).
\ee
Therefore, by Proposition \ref{IRW-comp} 
(with $K(z,s)=\alpha\, p_{2dT1[\kappa]}(X_{\kappa s},z)$), we get
\be{remain2-11}
\begin{aligned}
&\E_{\,\nu_\rho,0}\bigg(\exp\bigg[
\frac{\alpha}{\kappa}\int_0^t ds\, \left(\cP_T\phi\right)(Z_s) \bigg]\bigg)\\
&\qquad\leq
\ES_{\,0}\, \E_{\,\nu_\rho}^{\IRW}\bigg(\exp\bigg[
\frac{\alpha}{\kappa}\sum_{z\in\Z^d}p_{2dT\onek}(X_s,z)\int_0^t ds\,
\Big(\tilde\xi_\frac{s}{\kappa}(z)-\rho\Big)\bigg]\bigg)\\
&\qquad\leq
\exp\left[-\frac\alpha\kappa \rho t\right]\,
\ES_{\,0}\, \int_\Omega \nu_\rho(d\eta)\, \prod_{x\in A_\eta}\E_{\,x}^{\RW}\bigg(\exp\bigg[
\frac{\alpha}{\kappa}\sum_{z\in\Z^d}p_{2dT\onek}(X_s,z)\int_0^t ds\,
\delta_z\Big(Y_{\frac{s}\kappa}\Big)\bigg]\bigg)\\
&\qquad=
\exp\left[-\frac\alpha\kappa \rho t\right]\,
\ES_{\,0}\, \int_\Omega \nu_\rho(d\eta)\, \prod_{x\in A_\eta}\E_{\,x}^{\RW}\bigg(\exp\bigg[
\frac{\alpha}{\kappa}\int_0^t ds\, p_{2dT\onek}\Big(X_s,Y_{\frac{s}\kappa}\Big)
\bigg]\bigg),
\end{aligned}
\ee
where $A_\eta=\{x\in\Z^d : \eta(x)=1\}$ and $\E_x^{\RW}$ is expectation w.r.t.\ to simple 
random walk $Y=(Y_t)_{t\geq 0}$ on $\Z^d$ with step rate $1$ starting from $Y_0=x$. Using 
that $\nu_\rho$ is the Bernoulli product measure with density $\rho$, we get
\be{remain2-13}
\begin{aligned}
&\E_{\,\nu_\rho,0}\bigg(\exp\bigg[
\frac{\alpha}{\kappa}\int_0^t ds\, \left(\cP_T\phi\right)(Z_s) \bigg]\bigg)\\
&\qquad\leq
\exp\left[-\frac\alpha\kappa \rho t\right]\,
\int_\Omega \nu_\rho(d\eta)\,
\ES_{\,0}\Bigg( \prod_{x\in \Z^d}\E_{\,x}^{\RW}\, \bigg(\exp\bigg[
\eta(x)\frac{\alpha}{\kappa}\int_0^t ds\, p_{2dT\onek}\Big(X_s,Y_{\frac{s}\kappa}\Big)
\bigg]\bigg)\\
&\qquad=
\exp\left[-\frac\alpha\kappa \rho t\right]\,
\int_\Omega \nu_\rho(d\eta)\,
\ES_{\,0}\Bigg( \prod_{x\in \Z^d}\Big[1+\eta(x)\, w^{(t)}(x,t)\Big]\Bigg)\\
&\qquad=
\exp\left[-\frac\alpha\kappa \rho t\right]\,
\ES_{\,0}\bigg(\prod_{x\in\Z^d}\Big[1+\rho\, w^{(t)}(x,t)\Big]\bigg)
\leq
\exp\left[-\frac\alpha\kappa \rho t\right]\,
\ES_{\,0}\bigg(\exp\bigg[\rho\sum_{x\in\Z^d}w^{(t)}(x,t)\bigg]\bigg),
\end{aligned}
\ee
where $w^{(t)}\colon\,\Z^d\times[0,t)\to \R$ solves (\ref{remain2-5}). From
(\ref{remain2-5}) we deduce that
\be{remain2-13*}
\frac{\partial}{\partial s} \sum_{x\in\Z^d} w^{(t)}(x,s)
= \frac\alpha\kappa \sum_{x\in\Z^d} p_{2dT\onek}(X_{t-s},x)\big[1+w^{(t)}(x,s)\big].
\ee
Integrating (\ref{remain2-13*}) w.r.t.\ $s$ and inserting the result into
(\ref{remain2-13}), we get (\ref{remain2-3}).
\epr

\subsubsection{Proof of Lemma \ref{remainLem3}}
\label{remain3}

Next, we consider $v^{(t)}$ and $w^{(t)}$ as defined in (\ref{remain3-6}--\ref{remain3-9}), 
but with $|\alpha|$ instead of $\alpha$. 

\bpr
We begin by showing that, for $T$ large enough and all $x$, $s$, $t$ and $X_{(.)}$, we have 
$v^{(t)}(x,s)\leq 2$.

Do a Taylor expansion, to obtain ($s_0=0$)
\be{remain3-11}
v^{(t)}(x,s)=\sum_{n=0}^{\infty} \left(\frac{|\alpha|}{\kappa}\right)^n
\left(\prod_{l=1}^{n}\int_{s_{l-1}}^{s} ds_l\,\right)
\E_{\,x}^\RW\left(
\prod_{m=1}^{n}p_{2dT\onek}\Big(X_{t-s+s_m},Y_{\frac{s_m}\kappa}\Big)\right).
\ee
In Fourier representation the transition kernel of simple random walk with step rate
1 reads
\be{remain3-13}
p_s(x,y) = \oint dk \,\,e^{-i\, k \cdot (y-x)}\,e^{-s \widehat{\varphi}(k)},
\ee
where $\oint dk = (2\pi)^{-d}\int_{[-\pi,\pi)^d}dk$ and
\be{remain3-15}
\begin{aligned}
\widehat{\varphi}(k)
&= \frac{1}{2d}\sum_{\substack{x\in\Z^d \\ \|x\|=1}}
\left(1-e^{i\, k\cdot x}\right),
\quad k \in [-\pi,\pi)^d.
\end{aligned}
\ee
Combining (\ref{remain3-11}--\ref{remain3-13}), we get
\be{remain3-17}
\begin{aligned}
v^{(t)}(x,s)&=\sum_{n=0}^{\infty} \left(\frac{|\alpha|}{\kappa}\right)^n
\left(\prod_{l=1}^{n}\int_{s_{l-1}}^{s} ds_l\,\right)
\left(\prod_{m=1}^{n}\oint dk_m\,\right)\\
&\qquad\qquad\times
\E_{\,x}^\RW\Bigg(
\exp\bigg[i\sum_{p=1}^{n}\Big(Y_{\frac{s_p}\kappa}-X_{t-s+s_p}\Big)\cdot k_p\bigg]
\exp\bigg[-\Big(2dT\onek\Big)\sum_{q=1}^{n}\widehat{\varphi}(k_q)\bigg]\Bigg)\\
&=\sum_{n=0}^{\infty} \left(\frac{|\alpha|}{\kappa}\right)^n
\left(\prod_{l=1}^{n}\int_{s_{l-1}}^{s} ds_l\,\right)
\left(\prod_{m=1}^{n}\oint dk_m\,\right)
\exp\bigg[-i\sum_{p=1}^{n}\Big(X_{t-s+s_p}-x\Big)\cdot k_p\bigg]\\
&\qquad\qquad\times
\exp\bigg[-\Big(2dT\onek\Big)\sum_{q=1}^{n}\widehat{\varphi}(k_q)\bigg]\,
\E_{\,0}^\RW\Bigg(
\exp\bigg[i\sum_{r=1}^{n}Y_{\frac{s_r}\kappa}\cdot k_r\bigg]
\Bigg),
\end{aligned}
\ee
where in the last line we did a spatial shift of $Y$ by $x$. Because $Y$ has independent
increments, we have
\be{remain3-19}
\begin{aligned}
\E_{\,0}^\RW\Bigg(
\exp\bigg[i\sum_{r=1}^{n} Y_{\frac{s_r}\kappa}\cdot k_r\bigg]
\Bigg)
&=\E_{\,0}^\RW\Bigg(
\exp\bigg[i\sum_{r=1}^{n}(k_r+\cdots+k_n)\cdot
\Big(Y_{\frac{s_r}\kappa}-Y_{\frac{s_{r-1}}\kappa}\Big)\bigg]
\Bigg)\\
&=\prod_{r=1}^{n}\E_{\,0}^\RW\Bigg(
\exp\bigg[i(k_r+\cdots+k_n)\cdot
Y_{\frac{s_r-s_{r-1}}\kappa}\bigg]\Bigg)\\
&=\prod_{r=1}^{n}\sum_{z\in\Z^d}
p_{\frac{s_r-s_{r-1}}\kappa}(0,z)\,
\exp\big[i(k_r+\cdots+k_n)\cdot z\big]\\
&=\prod_{r=1}^{n}
\exp\bigg[-\frac{s_r-s_{r-1}}\kappa\, \widehat{\varphi}(k_r+\cdots+k_n)\bigg],
\end{aligned}
\ee
where the last line uses (\ref{remain3-13}). Since the r.h.s.\ is non-negative,
taking the modulus of the r.h.s.\ of (\ref{remain3-17}), we obtain
\be{remain3-21}
\begin{aligned}
v^{(t)}(x,s)&\leq\sum_{n=0}^{\infty} \left(\frac{|\alpha|}{\kappa}\right)^n
\left(\prod_{l=1}^{n}\int_{s_{l-1}}^{s} ds_l\,\right)
\left(\prod_{m=1}^{n}\oint dk_m\,\right)\\
&\qquad\qquad\times
\exp\bigg[-\Big(2dT\onek\Big)\sum_{q=1}^{n}\widehat{\varphi}(k_q)\bigg]\,
\E_{\,0}^\RW\Bigg(
\exp\bigg[i\sum_{r=1}^{n}Y_{\frac{s_r}\kappa}\cdot k_r\bigg]
\Bigg)\\
&=\sum_{n=0}^{\infty} \left(\frac{|\alpha|}{\kappa}\right)^n
\left(\prod_{l=1}^{n}\int_{s_{l-1}}^{s} ds_l\,\right)
\E_{\,0}^\RW\left(
\prod_{m=1}^{n}p_{2dT\onek}\Big(0,Y_{\frac{s_m}\kappa}\Big)\right),
\end{aligned}
\ee
where the last line uses (\ref{remain3-13}). Thus
\be{remain3-21*}
v^{(t)}(x,s)
\leq
\E_{\,0}^\RW\bigg(\exp\bigg[\frac{|\alpha|}{\kappa}\int_0^s du\,
p_{2dT\onek}\Big(0,Y_{\frac{u}\kappa}\Big)\bigg]\bigg)
\leq
\E_{\,0}^\RW\bigg(\exp\bigg[|\alpha|\int_0^\infty du\,
p_{2dT\onek}(0,Y_u)\bigg]\bigg).
\ee
Next, let $\cG$ denote the Green operator acting on functions $V\colon\,\Z^d\to
[0,\infty)$ as
\be{remain3-23}
(\cG V)(x)=\sum_{y\in \Z^d}G(x,y)V(y),\qquad x\in \Z^d.
\ee
With $p_t$ denoting the function $p_t(0,\cdot)$, we have
\be{remain3-25}
\Big\|\cG\Big(|\alpha|\, p_{2dT\onek}\Big)\Big\|_{\infty}
=|\alpha|\sup_{x\in\Z^d} \int_0^\infty ds\,\, \sum_{y\in\Z^d} p_s(x,y)\,
p_{2dT\onek}(0,y) \leq |\alpha| G_{2dT\onek}
\ee
with
\be{remain3-27}
G_{t}= \int_t^\infty ds\,\, p_s(0,0)
\ee
the truncated Green function at the origin. The r.h.s.\ of (\ref{remain3-25}) tends
to zero as $T\to\infty$. Hence Lemma 8.2.1 in G\"artner and den Hollander \cite{garhol04}
can be applied to the r.h.s.\ of (\ref{remain3-21*}) for $T$ large enough, to yield
\be{remain3-27*}
v^{(t)}(x,s)
\leq \frac1{1-\Big\|\cG\Big(|\alpha|\, p_{2dT\onek}\Big)\Big\|_{\infty}}
\da 1\quad \text{as}\quad T\to\infty, \quad\text{uniformly in } \kappa>0.
\ee
Thus, for $T$ large enough and all $x$, $s$, $t$, $\kappa$ and $X_{(.)}$, we have 
$v^{(t)}(x,s)\leq 2$, as claimed earlier. For such $T$, recalling (\ref{remain3-6}), 
we conclude from (\ref{remain2-5}) that $w^{(t)}\leq \bar{w}^{(t)}$, where $\bar{w}^{(t)}$ 
solves
\be{remain3-33}
\frac{\partial \bar{w}^{(t)}}{\partial s}(x,s) = \frac1{2d\kappa} \Delta \bar{w}^{(t)}(x,s)
+ \frac{2|\alpha|}\kappa\, p_{2dT\onek}(X_{t-s},x),
\qquad
\bar{w}^{(t)}(\cdot,0)\equiv 0,
\ee
The latter has the representation
\be{remain3-35}
\bar{w}^{(t)}(x,s)
=\frac{2|\alpha|}{\kappa}\int_0^s du\, \sum_{z\in\Z^d}
p_{\frac{s-u}{\kappa}}(x,z)\,p_{2dT\onek}(X_{t-u},z)
=\frac{2|\alpha|}{\kappa}\int_0^s du\,
p_{\frac{s-u}{\kappa}+2dT\onek}(x,X_{t-u}).
\ee
Hence,
\be{remain3-31}
\begin{aligned}
&\ES_{\,0}\bigg(\exp\bigg[
\frac{\alpha}{\kappa}\rho\int_0^t ds\, \sum_{x\in\Z^d} p_{2dT\onek}(X_{t-s},x)\,
w^{(t)}(x,s) \bigg]\bigg)\\
&\qquad\leq\ES_{\,0}\bigg(\exp\bigg[
\frac{|\alpha|}{\kappa}\rho\int_0^t ds\, \sum_{x\in\Z^d} p_{2dT\onek}(X_{t-s},x)\,
\bar{w}^{(t)}(x,s) \bigg]\bigg)\\
&\qquad=\ES_{\,0}\bigg(\exp\bigg[
\frac{2\alpha^2}{\kappa^2}\rho\int_0^t ds\, \int_0^s du\, 
p_{\frac{u-s}{\kappa}+4dT\onek}(X_{t-s},X_{t-u})\, \bigg]\bigg),
\end{aligned}
\ee
which proves the claim in (\ref{remain3-3}).
\epr

\subsubsection{Proof of Lemma \ref{remainLem4}}
\label{remain4}

The proof of Lemma \ref{remainLem4} is based on the following lemma. For $t\geq 0$,
$\alpha\in\R$ and $a,\kappa,T>0$, let
\be{remain4-5}
\Lambda_\alpha(t;a,\kappa,T)
=\frac{1}{2t}\log\ES_{\,0}\bigg(\exp\bigg[\frac{4\alpha^2}{\kappa^2}\rho
\int_0^t ds\, \int_s^{s+a\kappa^3} du\,\,
p_{\frac{u-s}{\kappa}+4dT\onek}(X_u,X_s)\bigg]\bigg)
\ee
and
\be{remain4-7}
\lambda_\alpha(a,\kappa,T)
=\limsup_{t\to\infty} \Lambda_\alpha(t;a,\kappa,T).
\ee

\bl{remainLem4.1}
If $d\geq 4$, then for any $\alpha\in\R$ and $a,T>0$,
\be{remain4.1-3}
\limsup_{\kappa\to\infty} \kappa^2 \lambda_\alpha(a,\kappa,T)
\leq 2\alpha^2\rho\, G_{4dT},
\ee
\el
where $G_t$ is the truncated Green function at the origin defined by (\ref{remain3-27}).
Before giving the proof of Lemma \ref{remainLem4.1}, we first prove Lemma \ref{remainLem4}.

\bprl{remainLem4}
Return to (\ref{remain4-3}). By the Cauchy-Schwarz inequality, we have
\be{remain4-9}
\begin{aligned}
&\frac{\kappa^2}{t}\log\ES_{\,0}\bigg(\exp\bigg[
\frac{2\alpha^2}{\kappa^2}\rho\int_0^t ds\, \int_s^t du\,
p_{\frac{u-s}{\kappa}+4dT\onek}(X_u,X_s) \bigg]\bigg)\\
&\quad\leq\frac{\kappa^2}{2t}\log\ES_{\,0}\bigg(\exp\bigg[\frac{4\alpha^2}{\kappa^2}\rho
\int_0^t ds\, \int_s^{s+a\kappa^3} du\,\,
p_{\frac{u-s}{\kappa}+4dT\onek}(X_u,X_s)\bigg]\bigg)\\
&\qquad+\frac{\kappa^2}{2t}\log\ES_{\,0}\bigg(\exp\bigg[\frac{4\alpha^2}{\kappa^2}\rho
\int_0^t ds\, \int_{s+a\kappa^3}^\infty du\,\,
p_{\frac{u-s}{\kappa}+4dT\onek}(X_u,X_s)\bigg]\bigg).
\end{aligned}
\ee
Moreover, by Lemma \ref{approxLaplem} and the fact that $d\geq 3$, we have
\be{remain4-11}
\begin{aligned}
\frac{1}{\kappa^2} \int_0^t ds\, \int_{s+a\kappa^3}^\infty du\,\,
p_{\frac{u-s}{\kappa}+4dT\onek}(X_u,X_s)
&\leq\frac{1}{\kappa^2} \int_0^t ds\, \int_{s+a\kappa^3}^\infty du\,\,
p_{\frac{u-s}{\kappa}}(0,0)\\
&\leq\frac{C}{\kappa}\, t \int_{a\kappa^2}^\infty du\,\,
\frac{1}{(1+u)^{\frac{d}{2}}}
\leq \frac{\tilde{C}}{a^\frac{1}{2}\kappa^2}\, t
\end{aligned}
\ee
with $C,\tilde{C}>0$. Combining (\ref{remain4-9}--\ref{remain4-11}) and Lemma \ref{remainLem4.1},
and letting $a\ra\infty$, we get (\ref{remain4-3}).
\eprl

The proof of Lemma \ref{remainLem4.1} is based on one further lemma. For $\gamma\geq 0$
and $a,\kappa,T>0$, let
\be{remain4.2-3}
\Lambda_\gamma(a,T) = \limsup_{\kappa\to\infty}
\frac{1}{a\kappa} \log
\ES_{\,0}\bigg(\frac\gamma{\kappa^2}\int_0^{a\kappa^3} ds\, \int_s^\infty du\,\,
p_{\frac{u-s}{\kappa}+4dT\onek}(X_s,X_u)\bigg).
\ee

\bl{remainLem4.2}
If $d\geq 4$, then for any $\gamma\geq 0$ and $a,T>0$,
\be{remain4.2-5}
\Lambda_\gamma(a,T) \leq \gamma\, G_{4dT}.
\ee
\el

\noindent
Before giving the proof of Lemma \ref{remainLem4.2}, we first prove
Lemma \ref{remainLem4.1}.

\bprl{remainLem4.1}
Split the integral in the exponent in the r.h.s.\ of (\ref{remain4-5}) as follows:
\be{remain4.1-7}
\begin{aligned}
&\int_0^t ds\, \int_s^{s+a\kappa^3} du\,\,
p_{\frac{u-s}{\kappa}+4dT\onek}(X_u,X_s)\\
&\qquad\leq
\left(\sum_{{k=1}\atop{\text{even}}}^{\lceil t/a\kappa^3\rceil}
+\sum_{{k=1}\atop{\text{odd}}}^{\lceil t/a\kappa^3\rceil}\right)
\int_{(k-1)a\kappa^3}^{k a\kappa^3} ds\, \int_s^{s+a\kappa^3} du\,\,
p_{\frac{u-s}{\kappa}+4dT\onek}(X_u,X_s).
\end{aligned}
\ee
Note that in each of the two sums, the summands are i.i.d. Hence, substituting
(\ref{remain4.1-7}) into (\ref{remain4-5}) and applying the Cauchy-Schwarz inequality,
we get
\be{remain4.9}
\Lambda_\alpha(t;a,\kappa,T)\leq\frac{\lceil t/a\kappa^3\rceil}{4t}\log\ES_{\,0}\bigg(
\exp\bigg[\frac{8\alpha^2}{\kappa^2}\rho
\int_0^{a\kappa^3} ds\, \int_s^{s+a\kappa^3} du\,\,
p_{\frac{u-s}{\kappa}+4dT\onek}(X_u,X_s)\bigg]\bigg).
\ee
Letting $t\to\infty$ and recalling (\ref{remain4-7}), we arrive at
\be{remain4.9*}
\lambda_\alpha(a,\kappa,T)
\leq\frac{1}{4a\kappa^3}\log\ES_{\,0}\bigg(\exp\bigg[
\frac{8\alpha^2}{\kappa^2}\rho
\int_0^{a\kappa^3} ds\, \int_s^{s+a\kappa^3} du\,\,
p_{\frac{u-s}{\kappa}+4dT\onek}(X_u,X_s)\bigg]\bigg).
\ee
Combining this with Lemma \ref{remainLem4.2} (with $\gamma=8\alpha^2\rho$), we obtain
(\ref{remain4.1-3}).
\eprl

The proof of Lemma \ref{remainLem4.2} is based on two further lemmas.

\bl{remainLem4.3}
For any $\beta>0$ and $M\in\N$,
\be{remain4.3.3}
\begin{aligned}
&\ES_{\,0}\Bigg(\exp\Bigg[\beta\sum_{k=1}^{M}\int_0^\infty ds\,\,
p_{\frac{s}\kappa+4dT\onek}\big(U_{k-1}(0),U_{k-1}(s)\big)\Bigg]\Bigg)\\
&\qquad\leq \prod_{k=1}^{M}\max_{y_1,\cdots,y_{k-1}\in\Z^d}\ES_{\,0}\Bigg(
\exp\Bigg[\beta\sum_{l=0}^{k-1}\int_0^\infty ds\,\,
p_{\frac{a\kappa^2}{M}l+\frac{s}\kappa +4dT\onek}(0,X_s+y_l)\Bigg]\Bigg),
\end{aligned}
\ee
where $U_k(t)=X(\frac{k}{M}a\kappa^3+s)$, $k\in\N_0$ and $y_0=0$.
\el

\bl{remainLem4.4}
For any $\beta>0$, $M\in\N$, $k\in\N_0$, and $y_0,\cdots,y_k\in\Z^d$,
\be{remain4.4-3}
\ES_{\,0}\Bigg(
\exp\Bigg[\beta\sum_{l=0}^{k-1}\int_0^\infty ds\,\,
p_{\frac{a\kappa^2}{M}l+\frac{s}\kappa +4dT\onek}(0,X_s+y_l)\Bigg]\Bigg)
\leq \exp\Bigg[
\frac{\beta\sum_{l=0}^{k}G_{\frac{a\kappa^2}{M}l+4dT\onek}}
{1-\beta\sum_{l=0}^{k}G_{\frac{a\kappa^2}{M}l+4dT\onek}}
\Bigg],
\ee
(recall {\rm (\ref{remain3-27})}), provided that
\be{remain4.4-5}
\beta\sum_{l=0}^{k}G_{\frac{a\kappa^2}{M}l+4dT\onek}<1.
\ee
\el

\noindent
The proofs of Lemmas \ref{remainLem4.3}--\ref{remainLem4.4} are similar to those
of Lemmas 6.3.1--6.3.2 in G\"artner and den Hollander \cite{garhol04}. We refrain
from spelling out the details. We conclude by proving Lemma \ref{remainLem4.2}.

\bprl{remainLem4.2}
As in the proof of Lemma 6.2.1 in G\"artner and den Hollander \cite{garhol04},
using Lemmas \ref{remainLem4.3}--\ref{remainLem4.4} we obtain
\be{remain4.2-7}
\begin{aligned}
&\frac{1}{a\kappa}\log\ES_{\,0}\Bigg(\exp\Bigg[\frac{\gamma}{\kappa^2}
\int_0^{a\kappa^3}ds\, \int_s^\infty du\,\,
p_{\frac{u-s}\kappa+4dT\onek}(X_s,X_u)\Bigg]\Bigg)\\
&\qquad\leq \frac{\gamma\sum_{l=0}^{M-1}G_{\frac{a\kappa^2}{M}l+4dT\onek}}
{1-\gamma\frac{a\kappa}{M}\sum_{l=0}^{M-1}G_{\frac{a\kappa^2}{M}l+4dT\onek}},
\end{aligned}
\ee
provided that
\be{remain4.2-9}
\gamma\frac{a\kappa}{M}\sum_{l=0}^{M-1}G_{\frac{a\kappa^2}{M}l+4dT\onek}<1.
\ee
But (recall (\ref{remain3-27}))
\be{remain4.2-9*}
\sum_{l=0}^{M-1} G_{\frac{a\kappa^2}{M}l+4dT\onek}
\leq G_{4dT\onek} + \sum_{l=1}^{M-1} G_{\frac{a\kappa^2}{M}l}.
\ee
{}From Lemma \ref{approxLaplem} we get $G_t \leq C/t^{\frac{d}{2}-1}$. Therefore
\be{Gsumsbds}
\frac{\kappa}{M} \sum_{l=1}^{M-1} G_{\frac{a\kappa^2}{M}l}
\leq \left\{\ba{ll}
\frac{C_3}{a^\frac12} &\mbox{if } d=3,\\[0.3cm]
\frac{C_4}{a} \frac{1}{\kappa}\log M &\mbox{if } d=4,\\[0.3cm]
\frac{C_d}{a^{\frac{d}{2}-1}} \frac{M^{\frac{d}{2}-2}}{\kappa^{d-3}} 
&\mbox{if } d\geq 5,
\ea
\right.
\ee
for some $C_d>0$, $d\geq 3$. Hence, picking $1 \ll M \leq C\kappa^2$, 
(\ref{remain4.2-9}) holds for $\kappa$ large enough when $d\geq 4$, and so
the claim (\ref{remain4.2-5}) follows from (\ref{remain4.2-7}) and 
(\ref{remain4.2-9*}--\ref{Gsumsbds}).
\eprl

\subsection{Extension to arbitrary $p$}
\label{S4.7}

In Sections \ref{S4.1}--\ref{S4.6} we proved Theorem \ref{Lyahighlim} for $p=1$. We
briefly indicate how the proof can be extended to arbitrary $p$.

As in (\ref{lscal}), after time rescaling we have, for any $p\in\N$,
\be{extp-3}
\lambda_p^*(\kappa)=\lim_{t\to\infty}\Lambda_p^*(\kappa;t)
\quad\text{with}\quad
\Lambda_p^*(\kappa;t)=\frac1t\log\E_{\,\nu_\rho,0,\cdots,0}\bigg(\exp\bigg[
\frac1\kappa\int_0^t ds\,\,\sum_{k=1}^{p}
\xi_{\frac{s}\kappa}\big(X_k(s)\big)\bigg]\bigg).
\ee
We are interested in the quantity
\be{extp-5}
\lambda_p^*(\kappa)-\frac\rho\kappa
=\lim_{t\to\infty}\frac1t\log\E_{\,\nu_\rho,0,\cdots,0}\bigg(\exp\bigg[
\frac1\kappa\int_0^t ds\,\,\sum_{k=1}^{p}
\Big(\xi_{\frac{s}\kappa}\big(X_k(s)\big)-\rho\Big)\bigg]\bigg).
\ee
As in (\ref{Gap1.57}), for $T>0$ let $\psi_p\colon\,\Omega\times(\Z^d)^p$
be defined by
\be{extp-7}
\psi(\eta,x_1,\cdots,x_p)
=\int_{0}^{T}ds \left(\cP_s^{(p)} \phi_p\right)(\eta,x_1,\dots,x_p)
\quad \mbox{with} \quad \phi_p(\eta,x_1,\cdots,x_p)=\sum_{k=1}^p\big[\eta(x_k)-\rho\big],
\ee
where $(\cP_s^{(p)})_{s\geq 0}$ is the semigroup with generator (compare with (\ref{lb2.7}))
\be{extp-9}
\cA^{(p)}=\frac1\kappa L+\sum_{k=1}^p \Delta_k.
\ee
Using (\ref{graph}), we obtain the representation (compare with (\ref{Gap1.69}))
\be{extp-11}
\psi_p(\eta,x_1,\cdots,x_p)
= \int_0^{T} ds\, \sum_{z\in\Z^d}\sum_{k=1}^{p}
p_{2ds\onek}(z,x_k)\big[\eta(z)-\rho\big]
=\sum_{k=1}^p \psi(\eta,x_k).
\ee
Let (compare with (\ref{lb2.5}))
\be{extp-15}
Z_s^{(p)}=\big(\xi_\frac{s}\kappa,X_1(s),\cdots,X_p(s)\big).
\ee

First, we have the analogue of Proposition \ref{lbholder}:
\bp{lbholder*}
For any  $p\in\N$, $\kappa,T>0$,
\be{Gap1.27*}
\begin{aligned}
\lambda_p^\ast(\kappa)-\frac{\rho}{\kappa}
&\,\,{\leq \atop \geq}\,\,
\frac1{2q}\limsup_{t\to\infty}\frac{1}{t}\log\E_{\,\nu_\rho,0}\bigg(\exp\bigg[
\frac{2q}r\int_0^t ds\, \left[
\left(e^{-\frac{r}{\kappa}\psi_p} \cA e^{\frac{r}{\kappa}\psi_p}\right)
-\cA\left(\frac{r}{\kappa}\psi_p\right)\right](Z_s^{(p)})\bigg]\bigg)\\
&\quad+\frac1{4q}\limsup_{t\to\infty}\frac{1}{t}\log\E_{\,\nu_\rho,0}\bigg(\exp\bigg[
\frac{4q}{\kappa}\int_0^t ds\, \left(\cP_T^{(p)}\phi_p\right)(Z_s^{(p)}) \bigg]\bigg),
\end{aligned}
\ee
where $1/r + 1/q=1$, for any $r,q>1$ in the first inequality and any $q<0<r<1$
in the second inequality.
\ep

Next, using (\ref{extp-11}), the bound
\be{extp-17}
\Big(\psi_p\big(\eta^{a,b},x_1,\cdots,x_p\big)
-\psi_p(\eta,x_1,\cdots,x_p)\Big)^2
\leq p\sum_{k=1}^{p}\Big(\psi\big(\eta^{a,b},x_k\big)
-\psi(\eta,x_k)\Big)^2,
\ee
and the estimate in (\ref{fdifest1}), we also have the analogue of Lemma \ref{mainLem2}:

\bl{mainLem2*}
Uniformly in $\eta\in\Omega$ and $x_1,\cdots,x_p\in\Z^d$,
\be{main2-3*}
\begin{aligned}
&\left[\left(e^{-\frac{r}{\kappa} \psi_p} \cA\,
e^{\frac{r}{\kappa}\psi_p}\right)
-\cA\,\left(\frac{r}{\kappa}\psi_p\right)\right](\eta,x_1,\cdots,x_p)\\
&\qquad= \frac{r^2}{2\kappa^2}\sum_{k=1}^{p}\sum_{e\colon\,\|e\|=1}
\Big(\psi(\eta,x_k+e)-\psi(\eta,x_k)\Big)^2
+ O\bigg(\Big(\frac{1}{\kappa}\Big)^3\bigg).
\end{aligned}
\ee
\el
Using H\"older's inequality to separate terms, we may therefore reduce to the case $p=1$
and deal with the first term in the r.h.s.\ of (\ref{Gap1.27*}) to get the analogue of
Proposition \ref{mainLem1}.

For the second term in (\ref{Gap1.27*}), using (\ref{graph}) we have
\be{extp-19}
\big(\cP_T^{(p)}\phi_p\big)(\eta,x_1,\cdots,x_p)
=\sum_{k=1}^{p}\sum_{z\in\Z^d} p_{2dT\onek}(z,x_k)\big[\eta(z)-\rho\big]
=\sum_{k=1}^{p}\big(\cP_T\phi\big)(\eta,x_k).
\ee
Using H\"older's inequality to separate terms, we may therefore again reduce to the
case $p=1$ and deal with the second term in the r.h.s.\ of (\ref{Gap1.27*}) to get
the analogue of Proposition \ref{remainLem1}.


\appendix

\section{Appendix}
\label{A}

In this appendix we give the proof of Proposition \ref{IRW-comp}.

\bpr
Fix $t\geq0$, $\eta\in\Omega$ and $K\colon\,\Z^d\times[0,\infty)\to \R$ such that
$S=\sum_{z\in\Z^d}\int_0^t ds\, |K(z,s)|<\infty$. First consider the
case $K\geq 0$. 
Since the $\xi$-process and the $\tilde\xi$-process are both monotone in their
initial configuration (as is evident from the graphical representation described
in Section \ref{S1.2}), it suffices to show that
\be{er2-8}
\E_{\,\eta}\Bigg(\exp\Bigg[
\sum_{z\in \Z^d}\int_0^t ds\,\, K(z,s)\, \xi_s(z)
\Bigg]\Bigg)
\leq
\E_{\,\eta}^{\IRW}\Bigg(\exp\Bigg[
\sum_{z\in \Z^d}\int_0^t ds\,\, K(z,s)\, \tilde\xi_s(z)
\Bigg]\Bigg),
\ee
for all $\eta\in\Omega$ such that $|\{x\in\Z^d : \eta(x)=1\}|<\infty$. This goes as follows.

Since $\xi_s(z)\in\{0,1\}$, we may write for any $r\in\R\setminus\{0\}$,
\be{er2-9}
\E_{\,\eta}\Bigg(\exp\Bigg[
\sum_{z\in \Z^d}\int_0^t ds\,\, K(z,s)\, \xi_s(z)
\Bigg]\Bigg)
=\E_{\,\eta}\Bigg(\exp\Bigg[
\sum_{z\in \Z^d}\int_0^t ds\,\, K(z,s)\, \frac{e^{r\, \xi_s(z)}-1}{e^r-1}
\Bigg]\Bigg).
\ee
By Taylor expansion, we get
\be{er2-10}
\begin{aligned}
&\E_{\,\eta}\Bigg(\exp\Bigg[
\sum_{z\in \Z^d}\int_0^t ds\,\, K(z,s)\,
\frac{e^{r\, \xi_s(z)}-1}{e^r-1}
\Bigg]\Bigg)\\
&=\exp\bigg[\frac{-t}{e^r-1}\,S\bigg]\E_{\,\eta}\Bigg(
\exp\Bigg[\sum_{z\in \Z^d}\int_0^t ds\,\,
K(z,s)\, \frac{e^{r\, \xi_s(z)}}{e^r-1}
\Bigg]\Bigg)\\
&=\exp\bigg[\frac{-t}{e^r-1}\,S\bigg]\\
&\quad\times
\sum_{n=0}^\infty\Bigg(
\frac{1}{e^r-1}\Bigg)^n\, \frac{1}{n!}
\Bigg(\prod_{j=1}^{n}\int_0^t ds_j\, \sum_{z_j\in \Z^d}\Bigg)
\Bigg(\prod_{j=1}^{n} K(z_j,s_j)\Bigg)
\E_{\,\eta}\Bigg(\exp\Bigg[
r\sum_{j=1}^{n}\xi_{s_j}(z_j)\Bigg]\Bigg).
\end{aligned}
\ee
According to Lemma 4.1 in Landim \cite{lan92}, we have for any $r\in\R$,
\be{er2-13}
\E_{\,\eta}\Bigg(\exp\Bigg[r
\sum_{j=1}^{n}\xi_{s_j}(z_j)\Bigg]\Bigg)
\leq\E_{\,\eta}^{\IRW}\Bigg(\exp\Bigg[r
\sum_{j=1}^{n}\tilde\xi_{s_j}(z_j)\Bigg]\Bigg).
\ee
Picking $r\geq 0$, combining (\ref{er2-9}--\ref{er2-13}), and using the analogue of 
(\ref{er2-10}) for $(\tilde\xi_t)_{t\geq 0}$, we obtain
\be{er2-15}
\E_{\,\eta}\Bigg(\exp\Bigg[
\sum_{z\in \Z^d}\int_0^t ds\,\, K(z,s)\, \xi_{s}(z)
\Bigg]\Bigg)
\leq\E_{\,\eta}^{\IRW}\Bigg(\exp\Bigg[
\sum_{z\in \Z^d}\int_0^t ds\,\, K(z,s)\,
\frac{e^{r\, \tilde\xi_{s}(z)}-1}{e^r-1}
\Bigg]\Bigg).
\ee
Now let $r\da 0$ and use the dominated convergence theorem to arrive at (\ref{er2-8}).

For the case $K\leq 0$ we can use the same argument with
\be{er2-17}
-\xi_s = \frac{e^{-r\xi_s}-1}{1-e^{-r}}.
\ee
\epr



\begin{thebibliography}{99}

\bibitem{arr85}
R.\ Arratia,
Symmetric exclusion processes: a comparison inequality and a large deviation result,
Ann.\ Probab.\ 13 (1985) 53--61.

\bibitem{chalanlee04}
C.-C.\ Chang, C.\ Landim and T.-Y.\ Lee,
Occupation time large deviations of two-dimensional symmetric simple exclusion process,
Ann.\ Probab.\ 32 (2004) 661--691.

\bibitem{deustr89}
J.\ D.\ Deuschel and D.\ W.\ Stroock,
\emph{Large Deviations},
Academic Press, London, 1989.

\bibitem{garhol04}
J.\ G\"artner and F.\ den Hollander,
Intermittency in a catalytic random medium,
EURANDOM Report 2004--019.
To appear in Ann. Prob.

\bibitem{garkon04}
J.\ G\"artner and W.\ K\"onig, The parabolic Anderson model, in: \emph{Interacting
Stochastic Systems} (J.-D.\ Deuschel and A.\ Greven, eds.), Springer, Berlin,
2005, pp.\ 153--179.

\bibitem{garmol90}
J.\ G\"artner and S.A.\ Molchanov, Parabolic problems for the Anderson Hamiltonian.
I.\ Intermittency and related topics. Commun.\ Math.\ Phys.\ 132 (1990) 613--655.

\bibitem{hol00}
F.\ den Hollander,
\emph{Large Deviations},
Fields Institute Monographs 14, American Mathematical Society, Providence, RI, 2000.

\bibitem{ka76}
T.\ Kato,
\emph{Perturbation Theory for Linear Operators} (2nd.\ ed.),
Springer, New York, 1976.

\bibitem{kessid03}
H.\ Kesten and V.\ Sidoravicius,
Branching random walk with catalysts,
Electr.\ J.\ Prob.\ 8 (2003) Paper no.\ 5, pp.\ 1--51.

\bibitem{kip87}
C.\ Kipnis,
Fluctuations des temps d'occupation d'un site dans l'exclusion simple sym\'etrique,
Ann.\ Inst.\ H.\ Poincar\'e\ Probab.\ Statist.\, 23 (1987) 21--35.

\bibitem{lan92}
C.\ Landim,
Occupation time large deviations for the symmetric simple exclusion process,
Ann.\ Probab.\ 20 (1992) 206--231.

\bibitem{lig85}
T.M.\ Liggett,
\emph{Interacting Particle Systems},
Grundlehren der Mathematischen Wissenschaften 276, Springer, New York, 1985.

\bibitem{sp76}
F.\ Spitzer,
\emph{Principles of Random Walk} (2nd.\ ed.),
Springer, Berlin, 1976.

\end{thebibliography}
\end{document}